\documentclass[a4paper, reqno]{amsart}
\usepackage[margin=1in]{geometry}
\allowdisplaybreaks

\title[]{On the structure of W-algebras in type A}

\author{Thomas Creutzig}
\address[T.C.]{Department Mathematik, FAU Erlangen–Nürnberg, Cauerstraße 11, 91058, Erlangen, Germany}
\email{creutzigt@math.fau.de}

\author[J.F.]{Justine Fasquel}
\address[Justine Fasquel]{School of Mathematics and Statistics, University of Melbourne, Parkville, Australia, 3010}
\email{justine.fasquel@unimelb.edu.au}

\author{Andrew R. Linshaw}
\address[A.L.]{Department of Mathematics, University of Denver, C. M. Knudson Hall, 2390 S. York St. Denver, CO 80210}
\email{andrew.linshaw@du.edu}

\author[Shigenori Nakatsuka]{Shigenori Nakatsuka}
\address[S.N.]{Department Mathematik, FAU Erlangen–Nürnberg, Cauerstraße 11, 91058, Erlangen, Germany}
\email{shigenori.nakatsuka@fau.de}

\usepackage{latexsym}
\usepackage{amsmath}
\usepackage{amsfonts}
\usepackage{amssymb,nccmath}
\usepackage{mathtools}
\usepackage{amsxtra}
\usepackage{mathrsfs}
\usepackage{eucal}
\usepackage[all]{xy}
\usepackage{color}
\usepackage{braket}
\usepackage{graphics, setspace}
\usepackage{etoolbox, textcomp}
\usepackage{comment}
\usepackage{changepage}
\usepackage{xcolor}
\definecolor{rouge}{rgb}{0.85,0.1,.4}
\definecolor{bleu}{rgb}{0.1,0.2,0.9}
\definecolor{violet}{rgb}{0.7,0,0.8}
\usepackage[colorlinks=true,linkcolor=bleu,urlcolor=violet,citecolor=rouge]{hyperref}
\usepackage{cleveref}
\usepackage{lscape}
\usepackage{caption}
\usepackage{subcaption}
\usepackage{enumitem}
\usepackage{graphicx}
\usepackage{arydshln}
\usepackage{wrapfig}
\usepackage[colorinlistoftodos]{todonotes}

\usepackage{cancel}

\usepackage[aligntableaux=center,smalltableaux=true]{ytableau}

\newtheorem{definition}{Definition}[section]
\newtheorem{proposition}[definition]{Proposition}
\newtheorem{theorem}[definition]{Theorem}
\newtheorem{ConjLetter}{Conjecture}

\newtheorem{ThmLetter}{Theorem}
\newtheorem{corollary}[definition]{Corollary}
\newtheorem{lemma}[definition]{Lemma}
\newtheorem{conjecture}[definition]{Conjecture}
\theoremstyle{remark}
\newtheorem{remark}[definition]{Remark}
\newtheorem{example}[definition]{Example}

\numberwithin{equation}{section}

\newcommand{\Z}{\mathbb{Z}}
\newcommand{\Q}{\mathbb{Q}}
\newcommand{\C}{\mathbb{C}}

\newcommand{\W}{\mathcal{W}}

\newcommand{\cO}{\mathcal{O}}

\newcommand{\g}{\mathfrak{g}}

\newcommand{\h}{\mathfrak{h}}
\newcommand{\m}{\mathfrak{m}}

\newcommand{\sll}{\mathfrak{sl}}
\newcommand{\gl}{\mathfrak{gl}}
\newcommand{\SL}{\mathrm{SL}}
\newcommand{\BP}{\mathrm{BP}}

\newcommand{\vir}[1]{\mathcal{L}(#1)}
\newcommand{\sL}{\mathcal{L}}

\newcommand{\ad}{\operatorname{ad}}

\newcommand{\ch}{\operatorname{ch}}
\newcommand{\Proj}{\pi_{P/Q}}

\newcommand{\Span}{\operatorname{Span}}

\newcommand{\pd}{\partial}
\newcommand{\ok}{$\checkmark$}

\newcommand{\osp}{\mathfrak{osp}}

\renewcommand{\mod}{\operatorname{-mod}}

\newcommand{\x}{\widehat{x}}
\newcommand{\y}{\widehat{y}}

\newcommand{\weyl}{\mathbb{V}}

\newcommand{\lbr}[2]{ {[} {#1} {}_\Lambda {#2} {]}}
\newcommand{\lm}[1]{\Lambda^{({#1})}}
\newcommand{\hwt}[1]{\mathrm{e}^{{#1}}}
\newcommand{\Com}[1]{\mathrm{Com}\left({#1}\right)}
\newcommand{\com}{\mathrm{Com}}
\newcommand{\tot}{\mathrm{tot}}

\newcommand{\ssqrt}[1]{\operatorname{\sqrt{\smash[b]{#1}}}}

\newcommand{\tA}{\texorpdfstring{$A$}{W}}
\newcommand{\tW}{\texorpdfstring{$\W$}{W}}
\newcommand{\tWinf}{\texorpdfstring{$\W_\infty$}{W}}
\newcommand{\tWinfcl}{\texorpdfstring{$\W_\infty[c,\lambda]$}{W}}
\newcommand{\tslN}{\texorpdfstring{$\sll_N$}{W}}
\newcommand{\tWtype}[2]{\texorpdfstring{$\W^k(\sll_{#1},f_{#2})$}{W}}

\newcommand\doi[2]{\href{http://dx.doi.org/#1}{#2}}
\newcommand{\arxiv}[2]{\href{https://arxiv.org/abs/#1}{#2}}

\begin{document}

\begin{abstract}
We formulate and prove examples of a conjecture which describes the $\W$-algebras in type $A$ as successive quantum Hamiltonian reductions of affine vertex algebras associated with several hook-type nilpotent orbits. This implies that the affine coset subalgebras of hook-type $\W$-algebras are building blocks of the $\W$-algebras in type $A$.
In the rational case, it turns out that the building blocks for the simple quotients are provided by the minimal series of the regular $\W$-algebras. 
In contrast, they are provided by singlet-type extensions of $\W$-algebras at collapsing levels which are irrational. 
In the latter case, several new sporadic isomorphisms between different $\W$-algebras are established.
\end{abstract}

\maketitle

\section{Introduction}
\subsection{Overview}

Let $\g$ be a simple Lie algebra or, more generally, a basic classical simple Lie superalgebra over the complex numbers $\C$. For each nilpotent element $f\in \g$ of even parity and $k \in \C$, one associates the $\W$-(super)algebra $\W^k(\g,f)$ at level $k$ as the quantum Hamiltonian reduction of the affine vertex (super)algebra $V^k(\g)$ associated with $\g$ at level $k$ \cite{FF90b, KRW03}. 

The $\W$-algebras are important objects for the study of the Whittaker models of the affine vertex algebras \cite{R} and also for 
the study of superconformal field theories in physics as they correspond to their symmetry algebras \cite{KW04}.
Moreover, they have appeared prominently in various areas such as integrable hierarchies \cite{DJKV, DS, FF90b}, moduli spaces of coherent sheaves \cite{AGT, BFN, B}, the quantum geometric Langlands program \cite{AF, Fr, FG20, Gai16} in mathematics and the 4d/2d duality \cite{SXY, WX, XYY}, the $\mathcal{N}=4$ super Yang-Mills gauge theories \cite{CG, GR, PR} in physics.

The regular $\W$-algebras in type $A$ (i.e., $\g\simeq\sll_N$), associated with the regular nilpotent elements, can be interpolated by a universal object called the universal $\W_\infty$-algebra.
It is a two-parameter vertex algebra freely and strongly generated by elements of conformal weights $2,3,4,\dots$ unique under some mild assumption.
It has been believed to exist since the 90s in physics and was constructed uniquely by one of the authors \cite{L}.
It turned out that the universal $\W_\infty$-algebra also interpolates certain subalgebras -- called the affine coset subalgebras -- of the $\W$-algebras associated with hook-type nilpotent elements. 
This object significantly improves our understanding of hook-type $\W$-algebras with major results such as the Feigin--Frenkel type duality, the Goddard--Kent--Olive type coset construction, and the reconstruction theorem \cite{CL1, CLNS, GR}.
 
The present paper aims to extend the scope to all the $\W$-algebras in type $A$ using building blocks obtained from the universal $\W_\infty$-algebra.
The description is achieved by formulating a general conjecture on the reduction by stages for the $\W$-algebras using the hook-type $\W$-algebras that is proven in small-rank examples.
As explained in the following, the conjecture lies at the center of various important applications, including the rational and non-rational representation theory of the $\W$-superalgebras, and the theory of closely related algebras such as the finite $\W$-algebras and the affine Yangians. 

\subsection{Reduction by stages}

\subsubsection{Main conjectures and results}
The $\W$-algebras $\W^k(\sll_N,f)$ in type $A$ are parametrized by the Jordan classes of nilpotent elements and thus by the partitions $\lambda$ of $N$ (we write $\lambda \vdash N$). Denote $f_\lambda\in\sll_N$ the corresponding representative. 
Accordingly, let us denote by $H_{f_\lambda}$ the functor which assigns the quantum Hamiltonian reduction $H_{f_\lambda}(M)$ to the modules $M$ over the affine vertex algebra $V^k(\g)$:
\begin{align*}
    H_{f_\lambda}\colon V^k(\g)\text{-}\mathrm{mod}\rightarrow \W^k(\sll_N,f_{\lambda})\text{-}\mathrm{mod}.
\end{align*}
It is exact on the Kazhdan--Lusztig category $\mathrm{KL}^k(\g)\subset V^k(\g)\text{-}\mathrm{mod}$ where the notion of $q$-characters makes sense \cite{A00}.

On the other hand, consider the sequence of nilpotent elements associated with partitions  
$\widehat{\lambda}_{1}, \widehat{\lambda}_{2}, \dots , \widehat{\lambda}_{n}$
with $\widehat{\lambda}_i=(1,\dots,1,\lambda_i)\vdash N_i:=N-\sum_{j>i}\lambda_{j}$.
They are called hook-type because of the shape of their associated Young diagrams.
It is then possible to define 
the following sequence of successive reductions
\begin{align*}
    H_{f_{\widehat{\lambda}_{1}}}H_{f_{\widehat{\lambda}_{2}}}\dots H_{f_{\widehat{\lambda}_{n}}}\left(V^k(\sll_N)\right).
\end{align*}
Indeed, at each step there is an affine subalgebra $V^{k_i^\sharp}({\sll}_{N_i})$, at certain level $k_i^\sharp$, sitting inside the vertex algebra $H_{f_{\widehat{\lambda}_{i+1}}}\dots H_{f_{\widehat{\lambda}_{n}}}\left(V^k(\sll_N)\right)$ to which we can apply the reduction functor $H_{f_{\widehat{\lambda}_{i}}}$.
It is then natural to conjecture that the resulting vertex algebra recovers the original $\W$-algebra 
$\W^k(\sll_N,f_\lambda)$. 
More strongly, we conjecture the following.

\begin{ConjLetter}\label{conj:successive_reductions}\hspace{0mm}
\begin{enumerate}[wide, labelindent=0pt]
\item For a partition $\lambda\vdash N$, there is an isomorphism of vertex algebras
\begin{equation*}
    \W^k(\sll_N,f_{\lambda})\simeq  H_{f_{\widehat{\lambda}_{1}}}H_{f_{\widehat{\lambda}_{2}}}\dots H_{f_{\widehat{\lambda}_{n}}}\left(V^k(\sll_N)\right).
\end{equation*}
\item 
Moreover, the functors 
\begin{align*}
    H_{f_\lambda}, H_{f_{\widehat{\lambda}_{1}}}H_{f_{\widehat{\lambda}_{2}}}\dots H_{f_{\widehat{\lambda}_{n}}}\colon \mathrm{KL}^k(\sll_N)\rightarrow \W^k(\sll_N,f_{\lambda})\text{-}\mathrm{mod}
\end{align*}
restricted to the Kazhdan--Lusztig category $\mathrm{KL}^k(\sll_N)$ are naturally isomorphic.
\end{enumerate}
\end{ConjLetter}
\noindent
This conjecture seems to be noticed by several independent groups of researchers simultaneously (see, for instance, \cite{GJ}). 
The name ``reduction by stages" is borrowed from papers \cite{GJ, Mor} that deal with an analogous construction in the \emph{finite $\W$-algebras} setting.

Our first result confirms Conjecture \ref{conj:successive_reductions} (1) for lower-rank cases:
\begin{ThmLetter}[Theorem \ref{thm: the case 2,2}/\ref{thm: the case 3,2}/\ref{thm: the case 2,2,1}]\label{thm:main lower rank case}
For all levels $k$, there are isomorphisms of vertex algebras 
\begin{align*}
        &H_{f_{2}}(\W^k(\sll_4,f_{1^2,2}))\simeq \W^k(\sll_4,f_{2,2}),\quad H_{f_{2}}(\W^k(\sll_5,f_{1^2,3}))\simeq \W^k(\sll_5,f_{2,3}),\\
        &H_{f_{1,2}}(\W^k(\sll_5,f_{1^3,2}))\simeq \W^k(\sll_5,f_{1,2,2}).
    \end{align*}
\end{ThmLetter}

The proof relies on the explicit description of the $\W$-algebras using the OPEs. The OPEs of $\W$-algberas are quite difficult to obtain in general.
One of the novelties of the present paper is the ``uniformly-organized" OPEs for $\W$-algebras of type $A$ up to rank four presented in \S \ref{appx1}. 
They have been obtained by fixing the cohomology classes of the BRST complexes \cite{KRW03, KW04} so that they have the maximal number of mutually commuting Virasoro elements and primary vectors.
Each Virasoro element is associated with a \emph{building block} (see \S\ref{sec:intro_building_blocks}), simplifying considerably the description of the $\W$-algebras structure.

Moreover, we give strong evidence for the second part of the conjecture by showing the coincidence of some important invariants, namely the central charges of vertex algebras and the $q$-characters of their modules.
\begin{ThmLetter}[Theorem \ref{Conj at the level of q}]\label{thm:main_characters}
The vertex algebras $\W^k(\sll_N,f_{\lambda})$ and  
$H_{f_{\widehat{\lambda}_{1}},f_{\widehat{\lambda}_{2}},\dots,f_{\widehat{\lambda}_{n}}}(V^k(\sll_N))$ have the same central charge and the $q$-characters. Moreover, the equality of the $q$-characters 
$$\ch[H_{f_\lambda}(M)](q) =\ch[H_{f_{\widehat{\lambda}_{1}},f_{\widehat{\lambda}_{2}},\dots,f_{\widehat{\lambda}_{n}}}(M)](q).$$ 
holds for $V^k(\g)$-modules $M$ in $\mathrm{KL}^k(\sll_N)$.
\end{ThmLetter}

The concept of reduction by stages and Conjecture \ref{conj:successive_reductions} have numerous important applications. We give several in the following and illustrate them with the examples in small ranks appearing in Theorem \ref{thm:main lower rank case}.

\subsubsection{Motivations from Whittaker models}
It is widely believed that the representation category of the $\W$-algebras $\W^k(\g,f)$ captures the representation category of the affine vertex algebra $V^k(\g)$ satisfying a certain nilpotency condition, called Whittaker models -- see \cite{R} for the case of regular nilpotent elements.

The finite $\W$-algebras are associative algebras which serve as the finite version of the $\W$-algebras $\W^k(\g,f)$. They are defined as the (finite) quantum Hamiltonian reduction of the enveloping algebras $\mathscr{U}(\g)$ and denoted by $\mathscr{U}(\g,f)$ \cite{Pre02}.
They are also realized as a quotient of $\W^k(\g,f)$ called the Zhu algebra \cite{DK06,Zhu}. 
Since Zhu's functor behave nicely with the Hamiltonian reduction \cite{Ara07}, Conjecture \ref{conj:successive_reductions} implies the following isomorphism
\begin{align}\label{intro: finite version}
    \mathscr{U}(\sll_N,f_\lambda)\simeq H_{f_{\widehat{\lambda}_{1}}}H_{f_{\widehat{\lambda}_{2}}}\dots H_{f_{\widehat{\lambda}_{n}}}(\mathscr{U}(\sll_N)).
\end{align}

After the first version of the present article was submitted, 
we were informed by H. Nakajima in a private communication that the classical analogue of this finite setting already appears in his paper \cite{NT17} with Y. Takayama as Hanany--Witten transitions of bow varieties. 
Interestingly, the work has a physics background \cite{HW} very closely related to our motivation explained in the next subsection. According to him, the quantization \eqref{intro: finite version} and the vertex algebraic analogue, Conjecture \ref{conj:successive_reductions} (1), were expected. 
As further evidence of \eqref{intro: finite version}, note that it implies a chain of embeddings
\begin{align*}
    \mathscr{U}(\sll_{N_1},f_{\lambda_1}) \subset \mathscr{U}(\sll_{N_2},f_{\lambda_1,\lambda_2}) \subset\cdots \subset \mathscr{U}(\sll_N,f_\lambda).
\end{align*}
This is already known to be true through their realizations as truncated shifted Yangians and has been used to study the representation theory \cite{BK, FRZ}.
The conjecture \eqref{intro: finite version} itself is consistent with a conjecture in Morgan's thesis \cite{Mor}, and several cases have been established in \cite{GJ}. 

The Skryabin's equivalence gives a correspondence between
the category of $\mathscr{U}(\g,f)$-modules and the category of Whittaker modules over $\mathscr{U}(\g)$. 
Compared with the vertex algebraic setting, it is the finite analog of the usual quantum Hamiltonian reduction.
It is of independent interest to find a natural subcategory of $\mathscr{U}(\g)$-modules on which the isomorphism \eqref{intro: finite version} for modules holds. 

Finally, it is worth mentioning that the analog of Conjecture \ref{conj:successive_reductions} has been established in the area of Whittaker models for the general linear group over local fields by Gomez, Gourevitch, Sahi \cite{GGS}. 
This suggests that there should be a common phenomenon of ``reduction by stages" for Whittaker models, which might be understood deeper from the framework for the quantum (local) geometric Langlands program \cite{Gai16, Winter school}.

\subsubsection{Motivations from physics}\label{sec: Motivations from physics}
The physical origin of Conjecture \ref{conj:successive_reductions} comes from the topologically twisted $\mathcal{N} = 4$ super Yang-Mills theory in four dimensions \cite{GR, PR}. In this theory, we consider various unitary gauge groups placed in four dimensions glued along a family of half-BPS interfaces $\mathcal{B}_{p,q}$ parametrized by integers $(p,q)$. 
We usually describe them by trivalent graphs, called $(p,q)$-webs, as in Figure \ref{fig: boundary conditions} below.

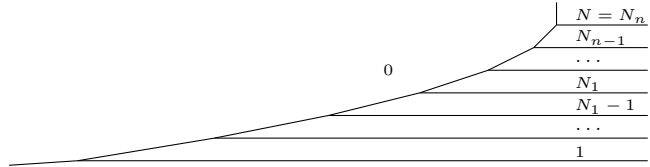
\begin{figure}[htbp]
  \centering
\begin{tikzpicture}[scale=0.6]
\draw (1,4) -- (1,4.5);
\draw (0.5,3.5) -- (1,4);
\draw (-0.5,3) -- (0.5,3.5);
\draw (-2,2.5) -- (-0.5,3);
\draw (-4,2) -- (-2,2.5);
\draw (-6.5,1.5) -- (-4,2);
\draw (-6.5,1.5) -- (-9.5,1);
\draw (1,4) -- (3,4);
\draw (0.5,3.5) -- (3,3.5);
\draw (-0.5,3) -- (3,3);
\draw (-2,2.5) -- (3,2.5);
\draw (-4,2) -- (3,2);
\draw (-6.5,1.5) -- (3,1.5);
\draw (-9.5,1) -- (3,1);
\draw (-11,0.9) -- (-9.5,1);
\node[right] at (-3,3){\tiny $0$ };
\node[right] at (1.2,4.2){\tiny$N=N_n$};
\node[right] at (1.2,3.7){\tiny{$N_{n-1}$}};
\node[right,font=\scriptsize] at (1.2,3.2){$\dots$};
\node[right] at (1.2,2.7){\tiny{$N_1$}};
\node[right] at (1.2,2.2){\tiny{$N_1-1$}};
\node[right,font=\scriptsize] at (1.2,1.7){$\dots$};
\node[right] at (1.2,1.2){\tiny$1$};
\end{tikzpicture}
 \caption{One associates to the $(p,q)$-web diagram the vertex algebra obtained from $V^k(\gl_N)$ by reduction by stages. The gauge group on each face is given by the unitary group $U(a)$ associated with the number, and the half-BPS interface $\mathcal{B}_{p,q}$ on each line is given by the title angle $p/q$.}
 \label{fig: boundary conditions}
 \end{figure}
\setlength{\unitlength}{1mm}

Vertex (super)algebras are assigned to $(p,q)$-webs in the following manner \cite{GR, PR}. 
    We start with the affine vertex algebra associated with the largest gauge group. In Figure \ref{fig: boundary conditions}, it is $V^k(\gl_N)$ where the level $k$ is the parameter for the topological twist of the theory. 
    Then we apply the quantum Hamiltonian reductions associated with hook-type partitions, corresponding to the fact that the gauge groups are reduced to smaller ones when passing across horizontal lines.
    Finally, we take the subalgebra of the resulting vertex algebra, called the affine coset, in terms of the affine vertex algebra associated with the remaining gauge group.
In that way, the vertex algebra associated to the $(p,q)$-web in Figure \ref{fig: boundary conditions} is 
    $H_{f_{\widehat{\lambda}_{1}}}H_{f_{\widehat{\lambda}_{2}}}\dots H_{f_{\widehat{\lambda}_{n}}}\left(V^k(\gl_N)\right)$.
    
On the other hand, the associated vertex algebra is also described from the local structure of the $(p,q)$-web in the following manner: it is an extension of the tensor product of vertex algebras associated with each trivalent (sub)graphs by tensor representations corresponding to the edges connecting vertices, which will be discussed in the next section. 
Based on this description, the vertex algebra is expected to be $\W^k(\gl_N,f_\lambda)$ and we obtained the general linear version of Conjecture \ref{conj:successive_reductions} (1):
\begin{align}\label{intro: gl version of reduction by stages}
    \W^k(\gl_N,f_\lambda)\simeq H_{f_{\widehat{\lambda}_{1}}}H_{f_{\widehat{\lambda}_{2}}}\dots H_{f_{\widehat{\lambda}_{n}}}\left(V^k(\gl_N)\right).
\end{align}
Since the decomposition $\gl_N=\sll_N\oplus \C$ induces the decomposition $\W^k(\gl_N,f_\lambda)\simeq \W^k(\sll_N,f_\lambda)\otimes \pi$ by using the rank-one Heisenberg vertex algebra $\pi$, \eqref{intro: gl version of reduction by stages} is equivalent to Conjecture \ref{conj:successive_reductions} (1).

\subsection{Universal objects and building blocks}\label{sec:intro_building_blocks}

The reduction by stages in Conjecture \ref{conj:successive_reductions} suggests that the quantum Hamiltonian reductions associated with hook-type partitions play a special role among all the reductions. 
For the beauty of the formulation, let us explain it in the general linear version \eqref{intro: gl version of reduction by stages}.

Given a vertex algebra $V$ and a subalgebra $W\subset V$, the coset subalgebra $\Com{W,V}$ is defined as 
\begin{align*}
    \Com{W,V}:=\{u\in V; [Y(u,z),Y(v,w)]=0,\ \forall v\in W\}\subset V
\end{align*}
by using the state-field correspondence $u\mapsto Y(u,z)$ of vertex algebras. 
In particular, for the hook-type $\W$-algebras $\W^k(\gl_N,f_{1^n,m})$ and their affine subalgebras $V^{k^\sharp}(\gl_{m})$ at level $k^\sharp$ which depends linearly on $k$, 
the coset subalgebras 
\begin{align*}
  C^k(\gl_N,f_{1^m,n}):=\Com{V^{k^\sharp}(\gl_{m}), \W^k(\gl_N,f_{1^n,m}))}
\end{align*}
are called the affine cosets of $\W^k(\gl_N,f_{1^n,m})$.
A remarkable property of these algebras is that they are all realized as quotients of the same universal object, called the $\W_{1+\infty}$-algebra 
$$\W_{1 +\infty}[c,\lambda]=\pi\otimes \W_\infty[c,\lambda].$$
Here 
$\W_\infty[c,\lambda]$ is the vertex algebra constructed in \cite{L} which depends on the central charge $c$ and another parameter $\lambda$ appearing in the OPEs of its strong generators. 
The vertex algebra $\W_{1 +\infty}[c,\lambda]$ admits is strongly generated by generators of conformal weights $1,2,3,\dots$ and admits a PBW base for these generators.
It is the unique vertex algebra object satisfying these properties under mild conditions.

The affine cosets $C^k(\gl_N,f_{1^m,n})$ are realized as truncations of $\W_{1 +\infty}[c,\lambda]$ for appropriate choices of parameters $(c,\lambda)$.
Hence there exist surjective maps 
\begin{align}\label{surjection from the universal objects}
    \W_{1 +\infty}[c,\lambda] \twoheadrightarrow C^k(\gl_N,f_{1^m,n}).
\end{align}
Even if the affine cosets $C^k(\gl_N,f_{1^m,n})$ admit strong generators of conformal weights $1,2,3,\dots, r(m,n)$, they do not have the PBW base property when $m\geq1$ which makes them difficult to study with a direct approach. 
Note that the affine cosets $C^k(\gl_N,f_{1^m,n})$ are simple for generic levels $k$ and thus are characterized only by the values $(c,\lambda)$ which themselves depend only on $m,n$.

The affine cosets $C^k(\gl_N,f_{1^m,n})$ and their modules serve as building blocks to reconstruct the $\W$-algebras $\W^k(\gl_N,f)$ in type $A$. 
Indeed, since we have the natural embeddings
\begin{align*}
    C^k(\gl_N,f_{1^m,n}) \otimes V^{k^\sharp}(\gl_{m})\hookrightarrow \W^k(\gl_N,f_{1^m,n}),
\end{align*}
Conjecture \ref{conj:successive_reductions} implies the following conformal embeddings, i.e. embeddings that preserve the Virasoro vector and thus the conformal structure.
\begin{ConjLetter}\label{conj:decomposition}
For generic levels $k$, there is a conformal embedding
\begin{equation*}
    \bigotimes_{i=1}^n C^{k_i^\sharp}(\gl_{N_i},f_{\widehat{\lambda}_i})\hookrightarrow \W^k(\gl_N,f_{\lambda})
\end{equation*}
for some levels $k_i^\sharp$ (see \eqref{eq: precise levels} for a precise formula).
\end{ConjLetter}

Although this conjecture is weaker than Conjecture \ref{conj:successive_reductions}, 
it is still important because of its application to the representation theory of $\W$-algebras (explained in \S \ref{intro: rational case} below) and its relation to the theory of affine Yangians.

The universal $\W_{1 +\infty}[c,\lambda]$ is intimately connected to another class of algebras called affine Yangians or, more generally, shifted quiver Yangians \cite{GLY}, which are another main players in the physics context reviewed in \S \ref{sec: Motivations from physics} \cite{GR, PR}.
More precisely, the universal enveloping algebra of $\W_{1 +\infty}[c,\lambda]$ is expected to be isomorphic to the affine Yangian $Y_{\hbar,\varepsilon}(\widehat{\gl}_1)$ after suitable completion \cite{GGLP}. The embedding in Conjecture \ref{conj:decomposition} is a shadow of a standard morphism for affine Yangians  
\begin{align*}
    \bigotimes_{i=1}^n Y_{\hbar_i,\varepsilon_i}(\widehat{\gl}_1) \hookrightarrow Y_{\hbar,\varepsilon}^\Xi(\widehat{\gl}_n)
\end{align*}
where $\Xi$ denotes a certain ``shift" of $Y_{\hbar,\varepsilon}(\widehat{\gl}_n)$.
Although $Y_{\hbar,\varepsilon}^\Xi(\widehat{\gl}_n)$ has not yet been rigorously constructed and studied in the literature, it suggests a new class of ``universal" $\W$-algebras which correspond to the shifted affine Yangians through the above-mentioned procedure. In \S \ref{sec: Universal objects}, we propose their potential constructions on the vertex algebra side.
We refer \cite{BR, KU, Ra20, U1, U2} and references therein for more details in this direction.

\subsection{Applications to the representation theory: rational case} \label{intro: rational case}

\subsubsection{Exceptional \tW-algebras}\label{sec: Exceptional Walgebras}

One of the main problems in the theory of vertex algebras is to classify strongly rational vertex algebras, whose associated variety is a point (lisse) and representation category is semisimple (rational). Based on the modular invariance of the admissible representations over the affine vertex algebras \cite{KW89}, it was conjectured that the so-called exceptional $\W$-algebras are strongly rational \cite{AvE22, KW08}. 
The exceptional $\W$-algebras were proven to be lisse for all Lie types by Arakawa \cite{A00} and rational in type $A$ by Arakawa-van Ekeren \cite{AvE22}. The general case was proven by McRae \cite{McR21}.

The description of $\W$-algebras through building blocks (Conjecture \ref{conj:decomposition}) provides a new and numerical approach to this problem.
Under mild assumptions, the conformal embedding in Conjecture \ref{conj:decomposition} implies the embedding for their simple quotients:
\begin{align*}
\begin{array}{lcc}
    \bigotimes_{i=1}^n C_{k_i^\sharp}(\gl_{N_i},f_{\widehat{\lambda}_i}) &\hookrightarrow& \W_k(\gl_N,f_{\lambda})\\
    \hspace{2cm}\cup &&\cup \\
    \bigotimes_{i=1}^n C_{k_i^\sharp}(\sll_{N_i},f_{\widehat{\lambda}_i})\otimes \pi^{\otimes (n-1)}\ &\hookrightarrow& \W_k(\sll_N,f_{\lambda}).
    \end{array}
\end{align*}
As the affine cosets $C_{k}(\sll_{N},f_{1^n,m})$ (and thus $C_{k}(\gl_{N},f_{1^n,m})$) are determined by the parameters $(c,\lambda)$, so if the parameters correspond to those of a strongly rational vertex algebra, then $\W_k(\sll_N,f_{\lambda})$ is also strongly rational by the general theory of vertex algebra extensions \cite{CMSY}.
We rediscover the strong rationality for the exceptional $\W$-algebras in type $A$ in this way. 
The $\W$-algebras $\W_k(\sll_N,f_\lambda)$ are exceptional for
\begin{align*}
    \lambda=(r,q^s)\quad (0\leq r<q)\quad\text{and}\quad k=-N+\frac{p}{q} 
\end{align*}
with $k$ admissible, i.e. $p\geq N$ and $\mathrm{gcd}(p,q)=1$.
The key observation is the conformal embedding for the hook-type $\W$-algebras for these levels (Theorem \ref{thm:branching_rules_hook_type}):
\begin{align}\label{intro: embedding at exceptional levels}
\W_{-(p-N)+\frac{p}{p-q}}(\sll_{p-N})\otimes L_{-(N-q)+\frac{p}{q}}(\sll_{N-q})\otimes \pi \hookrightarrow \W_{k}(\sll_{N},f_{1^{N-q},q})
\end{align}
where $\W_{-(p-N)+\frac{p}{p-q}}(\sll_{p-N})$ is a strongly rational (simple) $\W$-algebra associated with the regular nilpotent orbits and $L_{-(N-q)+\frac{p}{q}}(\sll_{N-q})$ is an admissible affine vertex algebra on which a new reduction can be performed.
Thereby, Conjecture \ref{conj:decomposition} implies the following.
\begin{ConjLetter}\label{conj:decomposition_rational}
For the admissible levels $k=-N+p/q$ and $f=f_{r,q^s}$ with $0\leq r<q$, there is a conformal embedding
\begin{equation*}
\left(\bigotimes_{\ell=0}^{s-1}\W_{-(p-N)+\frac{p-\ell q}{p-(\ell+1)q}}(\gl_{p-N})\right)\otimes \W_{-r+\frac{(p-sq)}{q}}(\gl_{r}) \hookrightarrow \W_k(\gl_N,f).
\end{equation*}
\end{ConjLetter}

The conjecture gives a new explanation for the strong rationality of $\W_k(\sll_N, 
f_{r,q^s})$. 
This method is a large generalization of some extremal cases appearing in the literature (see for instance \cite{ACL2, ACL, ALY, CL1, CL3}).
The conjecture is verified for the lower-rank cases thanks to Theorem \ref{thm:main lower rank case}. 
Interestingly, at these admissible levels one has 
\begin{equation}\label{eq: extended KRW conjecture}
    \W_k(\sll_N,f_{\lambda}) \simeq H_{f_{\widehat{\lambda}_{1}}}H_{f_{\widehat{\lambda}_{2}}}\dots H_{f_{\widehat{\lambda}_{n}}}\left(L_k(\sll_N)\right).
\end{equation}
It is a particular example of  Conjecture \ref{conj:successive_reductions} (2) as $\W_k(\sll_N,f_{\lambda})\simeq H_{f_{\lambda}}(L_k(\sll_N))$ holds \cite{AvE22}.
It is also crucial to explore the representation theory of the admissible affine vertex algebras $L_k(\sll_N)$ {beyond the admissible representations}. We will return to this subject later in \S \ref{sec: App: non-rational case}.

\subsubsection{Regular \tW-superalgebras}
In contrast to the $\W$-algebras, the nature of $\W$-superalgebras associated with simple Lie \emph{super}algebras
remains largely mysterious.
The main obstacle lies in the difficulty of developing the theory of admissible representations over affine Lie superalgebras. 
Conjecture \ref{conj:successive_reductions} offers an alternative to study $\W$-superalgebras through the webs of $\W$-algebras mentioned in \S \ref{sec: Motivations from physics} bypassing this difficulty.

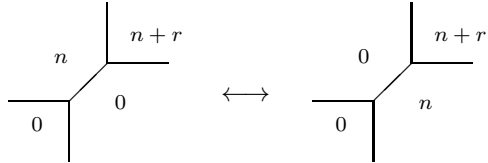
\begin{figure}[htbp]
\centering
\begin{picture}(80,22)(0,0)
\put(20,13){\line(1,0){8}}
\put(20,13){\line(0,1){8}}
\put(20,13){\line(-1,-1){5}}
\put(15,8){\line(-1,0){8}}
\put(15,8){\line(0,-1){8}}
\put(23,16){\footnotesize$n+r$}
\put(21,7){\footnotesize$0$}
\put(13,13){\footnotesize$n$}
\put(10,4){\footnotesize$0$}
\put(35,8){$\longleftrightarrow$}
\put(60,13){\line(1,0){8}}
\put(60,13){\line(0,1){8}}
\put(60,13){\line(-1,-1){5}}
\put(55,8){\line(-1,0){8}}
\put(55,8){\line(0,-1){8}}
\put(63,16){\footnotesize$n+r$}
\put(61,7){\footnotesize$n$}
\put(53,13){\footnotesize$0$}
\put(50,4){\footnotesize$0$}
\end{picture}
\caption{$(p,q)$-web for $\W^k(\sll_{n+r|n}, f_{n+r|n})$ and its flip.}
\label{fig: (p,q)-web for the regular Wsalg.}
\end{figure}

The regular $\W$-algebra $\W^k(\sll_{n+r|n}, f_{n+r|n})$ is associated with the left web in Figure \ref{fig: (p,q)-web for the regular Wsalg.}. 
Remarkably, the web has a $\Z_2$-symmetry, which flips to the web on the right. 
Under the flips, the associated vertex superalgebras are expected to be preserved. This gives a vast generalization of the celebrated Feigin--Frenkel duality for the regular $\W$-algebras \cite{FF91}. 
In the Figure \ref{fig: (p,q)-web for the regular Wsalg.}, the right web corresponds to the affine coset subalgebra inside the hook-type $\W$-algebra $\W^{\ell}(\sll_{n+r}, f_{1^n,r})$ tensored with a certain free field algebra \cite{PR}. 
Using the decomposition \eqref{intro: embedding at exceptional levels} in \S \ref{sec: Exceptional Walgebras}, one can use the rationality of exceptional $\W$-algebras to prove some new rationality results.
A precise statement is as follows.
\begin{ConjLetter}\label{conjecture for regular Wsalg}\hspace{0mm}
  \begin{enumerate}[wide, labelindent=0pt]
    \item There is an isomorphism of vertex superalgebras
    $$\W^k(\sll_{n+r|n}, f_{n+r|n}) \simeq \begin{dcases}\Com{V^{\ell + r}(\gl_n), \W^{\ell}(\sll_{n+r}, f_{1^n,r}) \otimes V_{\Z^n}}, & (r\geq 1),\\ 
    \Com{V^{\ell }(\gl_n), V^{\ell}(\sll_{n})\otimes \beta\gamma^{\otimes n}\otimes V_{\Z^n}}, & (r=0),\end{dcases}$$
    for generic levels $k,\ell$ satisfying $(k+r)(\ell + n+r) = 1$.
    \item The regular $\W$-superalgebras $\W_k(\sll_{n+r|n},f_{n+r|n})$ are strongly rational at levels 
    \begin{align*}
    k=-r+\frac{r}{p},\quad (p\geq n+r,\  r>n).
\end{align*}
\end{enumerate} 
\end{ConjLetter}

We note that the first statement unifies (1) the Feigin--Frenkel duality \cite{FF91} in type $A$, (2) the Kazama--Suzuki duality \cite{CGN, CL1}, and (3) Ito conjecture \cite{GL, Ito}.
As for the second statement, the cases $n=0,1$ are established in \cite{Ar2, CGN}
but the general conjecture is new as far as we know.

\subsection{Applications to the representation theory: non-rational case}
\label{sec: App: non-rational case}

\subsubsection{Inverse Hamiltonian reduction}
In parallel with the concept of reduction by stages, the idea of inverse Hamiltonian reduction has emerged.
Inverse Hamiltonian reductions are realized as conformal embeddings between $\W$-algebras associated with different nilpotent orbits of the same Lie algebra up to 
tensor products of the $\beta\gamma$-system, its localization $\Pi[0]$ and extensions $\Pi^\frac{1}{m}[0]$ (Remark \ref{CDO of GL1}).

The first example of inverse Hamiltonian reduction was given by Semikhatov \cite{Sem94} in string theory and relates the affine vertex algebra $V^k(\sll_2)$ to the Virasoro vertex algebra (see also \cite{Ad}):
\begin{align*}
    V^{k}(\sll_2) \hookrightarrow \W^{k}(\sll_2)\otimes \Pi[0].
\end{align*}
Despite recent efforts, the construction of the inverse Hamiltonian reductions is still limited and the results are mainly for hook-type $\W$-algebras e.g.\ \cite{ACG,AKR,FN23, Feh23, Feh23a}.
Conjectures~\ref{conj:successive_reductions} and \ref{conj:decomposition} suggest a large family of inverse Hamiltonian reductions providing paths from hook-type $\W$-algebras to all $\W$-algebras in type $A$.
The construction of inverse Hamiltonian reductions is the main ingredient for the proof of Theorem~\ref{thm:main lower rank case} and relies on the following isomorphisms.

\begin{ThmLetter}[Theorem \ref{thm: the case 2,2}/\ref{thm: the case 3,2}/\ref{thm: the case 2,2,1}]\label{thm:main of iHR for lower rank case }
For all levels $k$, there are isomorphisms of vertex algebras 
\begin{align*}
        &\W^k(\sll_4,f_{1^2,2})\simeq (\W^k(\sll_4,f_{2,2})\otimes \Pi^{\frac{1}{2}}[0])^{\SL_2},\quad \W^k(\sll_5,f_{1^2,3})\simeq (\W^k(\sll_5,f_{2,3})\otimes \Pi^{\frac{1}{2}}[0])^{\SL_2},\\
        &\W^k(\sll_5,f_{1^3,2})\simeq (\W^k(\sll_5,f_{1,2^2})\otimes \Pi^{\frac{1}{3}}[0]\otimes \beta\gamma)^{\SL_3}.
    \end{align*}
    Here $M^{\SL_n}$ denotes the maximal $\sll_n$-integrable submodule of $M$ with respect to an appropriate $\sll_n$-action. 
\end{ThmLetter}

To prove the isomorphisms, we use the uniformly-organized OPEs for $\W$-algebras presented in \S \ref{appx1}.
For instance, the first isomorphism in Theorem \ref{thm:main of iHR for lower rank case } is an extension of the original inverse Hamiltonian reduction for $\sll_2$ \cite{Sem94} as illustrated on the following diagram
\begin{align*}
	\begin{array}{ccc}
		\W^k(\sll_4,f_{1^2,2}) &\hookrightarrow& \W^k(\sll_4,f_{2,2})\otimes \Pi^{\frac{1}{2}}[0]\\
		\cup &&\cup \\
		V^{k+1}(\sll_2) &\hookrightarrow& \W^{k+1}(\sll_2)\otimes \Pi[0].
	\end{array}
\end{align*}

In general, hook-type $\W$-algebras were proven to be connected by inverse Hamiltonian reductions by Fehily in \cite{Feh23}. More precisely, there exist conformal embeddings
$$\W^k(\sll_{n+m},f_{1^m,n})\hookrightarrow \W^k(\sll_{n+m},f_{1^{m-1},n+1})\otimes\Pi[0]\otimes \beta\gamma^{\otimes(m-1)}.$$
By combining it with the building blocks decomposition in Conjecture \ref{conj:decomposition}, it gives strong evidence for the existence of a large family of inverse Hamiltonian reductions outside the hook-types. More precisely, we arrive at the following conjecture.
\begin{ConjLetter}\label{conj:embedding}
There exists an embedding of vertex algebras 
$$\W^k(\sll_N,f_{\sigma(\lambda)})\hookrightarrow \W^k(\sll_N,f_{\lambda})\otimes \Pi^{{\frac{1}{m}}}[0]\otimes \beta\gamma^{\otimes(m-1)},\quad (m=\lambda_1+\dots+\lambda_{a-1}=a-1)$$
where $1\leq a\leq n$ is the {unique} integer satisfying {$\lambda_{a-1}=1$} and $\lambda_{a}>1$,
$\sigma(\lambda)=(1^{a},\lambda_{a}-1,\lambda_{a+1},\dots,\lambda_n)$.
\end{ConjLetter}

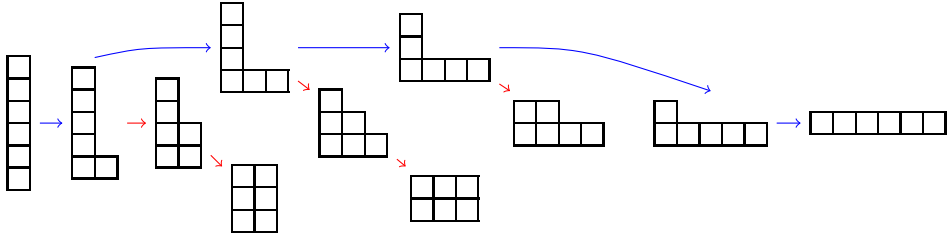
\begin{figure}[htbp]
    \centering
    \begin{tikzpicture}
			\node at (-0.1, 0) (aff)  {\ydiagram{1,1,1,1,1,1}};
			\node at (0.9, 0) (min) {\ydiagram{1,1,1,1,2}};
			\node at (2, 0) (2211) {\ydiagram{1,1,2,2}};
			\node at (3,-1) (222) {\ydiagram{2,2,2}};
			\node at (3,1) (3111) {\ydiagram{1,1,1,3}};
			\node at (4.3,0) (321) {\ydiagram{1,2,3}};
			\node at (5.5,-1) (33) {\ydiagram{3,3}};
			\node at (5.5,1) (411) {\ydiagram{1,1,4}};
			\node at (7,0) (42) {\ydiagram{2,4}};
			\node at (9,0) (sub) {\ydiagram{1,5}};
			\node at (11.2,0) (reg) {\ydiagram{6}};
			
			\draw[->,blue](aff) -- (min);
			\draw[->,blue] (sub) -- (reg);
			\draw[->,blue] (min.north) .. controls (1.5,1) .. (3111.west);
			\draw[->,blue] (3111) -- (411);
			\draw[->,blue] (411.east) .. controls (7.3,1) .. (sub.north);

			\draw[->,red](min)--(2211);
			\draw[->,red](2211)--(222);
			\draw[->,red](3111)--(321);
			\draw[->,red](321)--(33);
			\draw[->,red](411)--(42);
		\end{tikzpicture}
    \caption{The Young diagrams represent the $\W$-algebras of $\sll_6$. They are sorted according to the closure relations of nilpotent orbits with the smallest (corresponding to $V^k(\sll_6)$) on the left and the biggest ($\W^k(\sll_6,f_6)$) on the right.
    The arrows in blue describe inverse Hamiltonian reductions for hook-type $\W$-algebras given by \cite{Feh23}. The arrows in red describe the expected embeddings given by Conjectures~\ref{conj:embedding}.}
    \label{fig:IHR_sl6}
\end{figure}

The embeddings appearing in the conjecture significantly improve the global picture that inverse Hamitonian reductions exist along the closure relations of nilpotent orbits (see Figure~\ref{fig:IHR_sl6}). 
We remark that such embeddings are also expected in the relation of $\W$-algebras and vertex algebras arising from divisors on Calabi--Yau threefolds \cite{B}. 

\subsubsection{Relaxed modules and modularity}
The original motivation of the inverse Hamiltonian reduction is 
the study of the representation theory of admissible but non-exceptional $\W$-algebras with
the construction of the so-called relaxed highest weight modules \cite{FST}. 

At non-negative integer levels, the simple affine vertex algebra $L_k(\g)$ is known to be strongly rational \cite{FZ}.
Strongly rational vertex algebras behave very nicely from the perspective of representation theory. 
All their modules are ordinary and $\Z_\geq0$-graded. 
They form a (semisimple) modular tensor category \cite{Huang}.
Unfortunately for admissible non-rational $\W$-algebra, ordinary modules, and more generally highest weight modules, do not define a modular tensor category.

One of the most blatant examples of this failure is the ill-behaving Verlinde formula which sometimes returns negative coefficients for the decomposition of the fusion product in $L_k(\sll_2)$ at admissible non-integer levels \cite{KS88}.
It appeared later that the category $\mathcal{O}$ is not even closed under the fusion product \cite{Gab01} and the characters do not carry an action of the modular group \cite{R1}.
Considering instead the category of relaxed modules seems to remedy 
the second issue. 
Indeed, the relaxed modules of $L_k(\sll_2)$ at admissible non-integer levels were classified in \cite{AM95}
and the modularity of their characters verified in \cite{CR1, CR2}.

Nevertheless, the direct approach, still possible for $\sll_2$, cannot be generalized realistically.
In this context, the inverse Hamiltonian reduction provides a new efficient method to reconstruct the characters of the relaxed modules and verify their modularity.
Adamovi\'c \cite{Ad} showed that Semikhatov's embedding \cite{Sem94} 
still holds when considering simple quotients at non-integer admissible levels.
At these levels, the simple Virasoro vertex algebra is rational and the embedding reflects the modular tensor category given by its highest weight modules onto  
the nonfinite nonsemisimple category defined by relaxed modules of $L_k(\sll_2)$ (see for instance \cite{C2}).
In an attempt to develop this new approach, $\W$-algebras associated with $\sll_3$ have been studied in detail \cite{ACK, FRR, FR, KRW}.
Inverse Hamiltonian reduction appears then as a promising tool that carries back the modularity structure of the highest-weight modules of a rational $\W$-algebra on the set of relaxed modules of a non-rational $\W$-algebra.

The decomposition \eqref{intro: embedding at exceptional levels} implies the existence of inverse Hamiltonian reductions for the simple quotients of the $\W$-algebras appearing in Theorem \ref{thm:main of iHR for lower rank case } at specific levels.
Precisely
\begin{align*}
        &\W_k(\sll_4,f_{1^2,2})\simeq (\W_k(\sll_4,f_{2,2})\otimes \Pi^{\frac{1}{2}}[0])^{\SL_2},\quad \W_k(\sll_5,f_{1^2,3})\simeq (\W_k(\sll_5,f_{2,3})\otimes \Pi^{\frac{1}{2}}[0])^{\SL_2},\\
        &\W_k(\sll_5,f_{1^3,2})\simeq (\W_k(\sll_5,f_{1,2^2})\otimes \Pi^{\frac{1}{3}}[0]\otimes \beta\gamma)^{\SL_3}
\end{align*}
hold for exceptional levels $k$ defined in Cororally \ref{iHR for simple quotients}.
This implies that the representation category of the hook-type $\W$-algebras can be studied by combining results on the exceptional $\W$-algebras 
and the inverse Hamiltonian reductions for $L_{\ell}(\sll_2)$ and $L_{\ell}(\sll_3)$. 

\subsection{Collapsing levels and singlet \tW-algebras}
The representation theory of the simple $\W$-algebras $\W_k(\g,f)$ is quite wild for generic levels as $\W_k(\g,f)=\W^k(\g,f)$.
On the other extreme, the representation theory for strongly rational $\W$-algebras (e.g. at exceptional levels) is always semisimple. 
The collapsing levels lie in the intermediate and form an interesting class where the representation theory is often nonsemisimple but still accessible.
Collapsing historically referred to levels where the simple $\W$-algebra $\W_k(\g,f)$ is isomorphic to a simple affine vertex algebra \cite{AKMPP1}.
They have been tracked down from the representation theory \cite{AKMPP1, AKFPP, AMP, AFP0, AFP}, the geometry of nilpotent Slodowy slices \cite{AvEM}, and the 4D/2D duality in physics \cite{XY}.
In this paper, generalizing the initial definition, a level $k$ is called collapsing as soon as the strong generating type of $\W_k(\g,f)$ is strictly smaller than that of $\W^k(\g,f)$. 
In this wider sense of collapsing levels, the following coincidences between different $\W$-algebras are established.
\begin{ThmLetter}\label{Main thm C}
There exist isomorphisms of vertex algebras as follows.
  \begin{enumerate}[wide, labelindent=0pt]
    \item {\textup{[Theorem \ref{thm: collapsing for higherrank hooks}]}} For $m\geq 1$ and $n\geq 3$,
\begin{align*}
    \W_{-(n+m)+(m+1)}(\sll_{n+m},f_{1^m,n})\simeq \pi\otimes \W_{-(n-1)+(m+1)}(\sll_{n-1}).
\end{align*}
    \item {\textup{[Proposition \ref{prop:collaps_sl5_32}/Theorem \ref{thm: collapsing for hook type in cases}/Corollary \ref{cor: 3.2 v.s. 2.1}]}}
    \begin{align*}
    &\W_{-4+4/3}(\sll_4,f_{1^2,2})\simeq\W_{-5+5/3}(\sll_5,f_{1,2^2}), && \W_{-4+4/5}(\sll_4,f_{1,3})\simeq \W_{-3+3/5}(\mathfrak{sp}_4,f_{\mathrm{sub}}),\\
    &\W_{-5+5/4}(\sll_5,f_{2,3}) \simeq \W_{-3+3/4}(\sll_3,f_{1,2}), &&\W_{-5+5/2}(\sll_5,f_{2,3})\simeq \W_{-3+3/1}(\sll_3,f_{1,2}).
\end{align*}
\end{enumerate}
\end{ThmLetter} 

Some of the isomorphisms in the second case confirm in lower-rank cases a conjecture proposed by Xie and Yan \cite[\S3.1]{XY}.
The proof is rather interesting as we derive certain realizations 
using singlet-type extensions of $\W$-algebras \cite{Ad2, ACGY, ACKR, CMY3, CNS, CRW}. 
Indeed, we find that many of the $\W$-algebras at collapsing levels can be decomposed in terms of various singlet-type algebras. For example, we establish the decompositions
\begin{align*}
    &\W_{-5+5/2}(\sll_5,f_{2,3})\simeq \bigoplus_{n\in \Z}S_{2;n}\otimes \pi_{n}, && \W_{-5+4/3}(\sll_5,f_{2,3})\simeq \bigoplus_{n\in \Z}S_{2;n}\otimes S_{3;n}\otimes \pi_{\ssqrt{-5/2}n},
\end{align*}
in terms of the singlet algebras $S_{p;0}$, which are new as far as we know.
See \S\ref{sec: collapsing levels} for details and more examples.
Moreover, the isomorphism at level $-5+4/3$ seems to extend into the following family.
\begin{ConjLetter} 
We have the following decomposition
\begin{align*}
    \W_{-(p+2)+\frac{p+1}{p}}(\sll_{p+2},f_{2,p})\simeq \bigoplus_{n\in \Z}S_{2;n}\otimes S_{p;n}\otimes \pi_{\ssqrt{-\frac{p+2}{2}}n}.
\end{align*}
\end{ConjLetter}

These decompositions imply that singlet-type extensions of $\W$-algebras might be regarded as building blocks at collapsing levels. 
This conjecture is especially interesting from the representation theory point of view since one might describe the representation categories of $\W$-algebras at non-admissible levels, which are hard to study in general, in terms of those of singlet algebras, which are well-understood (see for instance \cite{Ad2, CMY, CMY2}), by using the theory of vertex algebra extensions \cite{CKM2, CLR, HKL}.

\subsection*{Organization of the paper} 
The rest of the paper is organized as follows. 
In \S \ref{sec: Iterated Quantum Hamiltonian reductions in type A}, we give the strong generating type of $\W$-algebras in type $A$ and explain the concept of reduction by stages. Then we prove Theorem \ref{thm:main_characters}.
In \S \ref{sec: Universal objects}, we collect some useful facts concerning the $\W_\infty$-algebra and integral form of $\W$-algebras including the large level limit.
We discuss possible ways to construct more universal objects and formulate the rationality conjecture of $\W$-superalgebras in type $A$. 
In \S \ref{sec: Non-trivial examples in low ranks}, we prove Theorems \ref{thm:main lower rank case} and \ref{thm:main of iHR for lower rank case }. 
In \S \ref{sec: collapsing levels}, we study the collapsing levels of $\W$-algebras and establish various decomposition theorems, including Theorem \ref{Main thm C}.  
Finally in \S \ref{sec: exceptional levels}, we study the exceptional $\W$-algebras and establish the decomposition theorem for hook-type nilpotent orbits (Theorem \ref{thm:branching_rules_hook_type}).
In \S \ref{appx1}, we present the OPEs for all the $\W$-algebras in type $A$ up to rank four. 
In \S \ref{Winfinity algebra}, we recollect the OPEs of the universal $\W_\infty$-algebra from \cite{L}.

\subsection*{Acknowledgements} 
S.N. would like to thank Dražen Adamovi\'c for the exciting discussions that triggered this project. 
We thank Hiraku Nakajima for explaining his work \cite{NT17}
and his expectation on Conjecture A (1) and its finite analogue.
S.N. thanks Gurbir Dhillon, Thibault Juillard, Dmytro Matvieievskyi, and Jinwei Yang for interesting correspondence during the work.
J.F. thanks David Ridout and Zac Fehily for their useful feedback.
Part of the work was done while S.N. and J.F. visited the Universities of Alberta and Denver. They are grateful for their hospitality. 
T.C. is supported by NSERC Grant Number RES0048511.
J.F. is supported by a University of Melbourne Establishment Grant and Andrew Sisson Support Package 2023.
A.L. is supported by NSF Grants DMS-2001484 and DMS-2401382 and  Simons Foundation Grant MPS-TSM-00007694.
S.N. is supported by JSPS Overseas Research Fellowships Grant Number 202260077.

\section{Quantum Hamiltonian reductions in type \tA}\label{sec: Iterated Quantum Hamiltonian reductions in type A}

We start recalling the general construction of $\W$-algebras using quantum Hamiltonian reductions due to Kac--Roan--Wakimoto \cite{KRW03}. Then we restrict to type $A$ and construct vertex algebras applying the so-called reduction by stages. We show that the character of the resulting vertex algebra coincides with the character of a well-chosen $\W$-algebra.
This proves Conjecture~\ref{conj:successive_reductions} at level of characters.

\subsection{\tW-algebras}\label{sec:BRST}
Let $\g$ be a finite-dimensional simple Lie algebra over $\C$ with Cartan subalgebra $\h$ and let $\Delta$ denote the set of roots.
In the following, we fix an element $x\in\h$ which induces an $\frac{1}{2}\Z$-grading $\Gamma_x$ on $\g$ given by the eigenspaces of $x$.
For $i\in\frac{1}{2}\Z$, we define 
$$\Delta_i^{x}:=\{\alpha \in \Delta\mid \alpha(x)=i\}\subset \Delta,$$ 
and $\Delta_+^x=\bigcup_{i>0}\Delta_i^{x}$ that is the subset consisting of roots whose root subspaces have positive $\Gamma_x$-grading.

In the following, we fix a nilpotent element $f\in\g$ that we embed into a $\sll_2$-triple $(e,h=2x,f)$ whose existence is guaranteed by the Jacobson--Morosov theorem. We assume that $\Gamma_x$ is a good grading (see for instance \cite{BG,EK}).
Then, to each $\alpha\in \Delta_+^{x}$, we associate a copy of a $bc$-system $\mathcal{F}_{\alpha}[x]$ generated by the odd fields $\varphi_\alpha,\varphi_\alpha^*$ satisfying the $\Lambda$-bracket
\begin{align*}
    \lbr{\varphi_\alpha}{\varphi_\alpha}=\lbr{\varphi^*_\alpha}{\varphi^*_\alpha}=0,\quad \lbr{\varphi_\alpha}{\varphi^*_\alpha}=1.
\end{align*}
The $bc$-system $\mathcal{F}_{\alpha}[x]$ is equipped with a conformal vector 
\begin{equation}\label{eq:alphax}
    L_\alpha^x=(1-\alpha(x)) (\partial \varphi_\alpha^*) \varphi_\alpha - \alpha(x) \varphi_\alpha^*(\partial\varphi_\alpha),
\end{equation}
whose central charge is given by
\begin{align}\label{cc for charged fermion}
    c_\alpha^x=-c[\alpha(x)],\quad c[t]:=12t(t-1)+2.
\end{align}
Moreover the fields $\varphi_\alpha,\varphi_\alpha^*$ have conformal weights $\Delta(\varphi_\alpha)=1-\alpha(x)$ and $\Delta(\varphi_\alpha^*)=\alpha(x)$.

On the other hand, we associate to the vector space $\g_{1/2}$ a product of $\beta\gamma$-systems, denoted by $\Phi(\g_{1/2})$ generated by the even fields $\Phi_\alpha$ ($\alpha \in \Delta^x_{1/2}$) satisfying the $\Lambda$-bracket
\begin{align*}
    \lbr{\Phi_\alpha}{\Phi_\beta}=\langle e_\alpha,e_\beta\rangle
\end{align*}
where $e_\alpha$ is a fixed root vector associated to $\alpha\in\Delta_{1/2}^x$ and 
$\langle x,y\rangle=(f,[x,y])$, for all $x,y \in \g_{1/2}$. Here $(\text{-},\text{-})$ is a normalized invariant bilinear symmetric form on $\g$.
We equip $\Phi(\g_{1/2})$ with the conformal vector 
$$L^{x}_{\Phi} = \frac{1}{2} \sum_{\alpha\in \Delta_{1/2}^x}(\partial \Phi_\alpha^*) \Phi_{\alpha}$$
where $\Phi_{\alpha}^*$ is uniquely determined by $\lbr{\Phi_\beta}{\Phi_\alpha^*}=\delta_{\alpha,\beta}$. Later, we will shift this conformal vector as 
\begin{align}\label{eq:phi_xa}
   L^{x+a}_{\Phi}=L^{x}_{\Phi}+\partial a^{\Phi}, \quad a^{\Phi}=-\frac{1}{2}\sum_{\alpha\in \Delta_{1/2}^{x}}\alpha(a)\,\Phi_\alpha \Phi_{\alpha}^*.
\end{align}
for $a\in\h$ satisfying $[a,f]=0$.
The conformal vector $L^{x+a}_{\Phi}$ has central charge 
\begin{align}\label{cc for neutral fermion}
    c_{\Phi}^{x+a}=\frac{1}{2}\sum_{\alpha \in \Delta_{1/2}^{x}}c[\alpha(x+a)]
\end{align}
and $\Phi_\alpha$ has the conformal weight 
$\Delta(\Phi_\alpha)=1-\alpha(x+a)$.

The $\W$-algebra associated with $f$ at level $k\in\C$ is the BRST reduction of the affine vertex algebra $V^k(\g)$.
More concretely, it is defined as the cohomology
\begin{equation}\label{eq:complex}
\W^k(\g,f)=H_{f}^0(V^k(\g)):=H^0(C_{f}(V^k(\g)),d),\quad C_{f}(V^k(\g)):=V^k(\g)\otimes \mathcal{F}(\g, f)
\end{equation}
where 
\begin{align*}
\mathcal{F}(\g, f)=\Phi(\g_{1/2})\otimes \bigotimes_{\alpha \in \Delta_{+}^x}\mathcal{F}_\alpha[x],\quad d=\int Y(Q,z)\mathrm{d}z, 
\end{align*}
with 
$$Q=\sum_{\alpha\in \Delta^{x}_+} (e_\alpha+(f,e_\alpha))\otimes \varphi_\alpha^*-\frac{1}{2}\sum_{ \alpha,\beta,\gamma \in\Delta^{x}_+}c_{\alpha,\beta}^\gamma \varphi_\gamma\varphi_\alpha^*\varphi_\beta^*
+\sum_{\alpha\in \Delta^{x}_{1/2}}\Phi_{\alpha}\otimes\varphi_\alpha^*$$
where $c_{\alpha,\beta}^\gamma$ is given by $[e_\alpha,e_\beta]=\sum c_{\alpha,\beta}^\gamma e_\gamma$. 
Despite appearances, the $\W$-algebra $\W^k(\g,f)$ depends on the nilpotent orbit of $f$ rather than $f$ itself.
In addition, at non-critical levels $k\neq-h^\vee$, where $h^\vee$ is the dual Coxeter number of $\g$, $\W^k(\g,f)$ inherits a conformal structure from $V^k(\g)$. The conformal vector is given by
\begin{align}\label{KRW}
    L_{f}=\left( L^{\g}_{\mathrm{sug}} + \partial x \right)+ \sum_{\alpha\in \Delta_{+}^x}L_\alpha^x+ L_{\Phi}^x
\end{align}
where $L^{\g}_{\mathrm{sug}}$ is the Sugawara conformal vector of $V^k(\g)$. 
The central charge of $L_{f}$ is  
$$\mathrm{c}(\W^k(\g,{f})) = c_{\mathrm{aff}}({f})+ c_{\mathrm{gh}}({f})$$
where 
$$c_{\mathrm{aff}}({f})=\frac{k \dim (\g)}{k + h^\vee} -12k(x, x),\quad c_{\mathrm{gh}}({f}) =c_{\Phi}^x+ \sum_{\alpha\in \Delta_{+}^x} c_\alpha^x.$$
Since the $L_{{f},0}$-action is semisimple on $\W^k(\g,{f})$ with finite-dimensional eigenspaces, the $q$-character
\begin{align*}
    \ch[\W^k(\g,f)](q)=\mathrm{tr}_{\W^k(\g,f)}q^{L_{{f},0}}
\end{align*}
is well-defined.
Thanks to the cohomology vanishing $H_{f}^{\neq0}(V^k(\g))=0$ \cite{KW04},
it agrees with the Euler--Poincar\'{e} character of the BRST complex:
\begin{align*}
    \ch[\W^k(\g,f)](q)=\mathrm{str}_{V^k(\g)\otimes \mathcal{F}(\g, f)}q^{L_{f,0}}.
\end{align*}
We note that the naive computation of the right-hand side is not well-defined at first sight but can still proceed using, for instance, the limit procedure described in \cite{KRW03}.  

Finally, given a $V^k(\g)$-module $M$, we define the BRST complex $C_f(M)$ as in \eqref{eq:complex}. The cohomology $H^0_f(M):=H^0(C_f(M),d)$ is a module over $\W^k(\g,f)$.

\subsection{\tW-algebras of type \tslN}\label{sec:WslN}
From now on, $\g$ is the Lie algebra $\sll_N$ otherwise specified.
As mentioned previously, the $\W$-algebra $\W^k(\g,f)$ depends on the nilpotent element $f$ and more precisely its nilpotent orbit.
In $\sll_N$, nilpotent orbits are parametrized by the partitions of $N$. 
In what follows, we describe the strong generating type -- the conformal weights of the strong generators -- of the $\W$-algebra 
$$\W^k(\gl_N,f_\lambda)\simeq \pi\otimes \W^k(\sll_N,f_\lambda)$$
corresponding to the nilpotent orbit associated with the partition 
$\lambda=(\lambda_1\leq\dots \leq\lambda_n)\vdash N$. 
The nilpotent element $f_\lambda$ is a representative of the orbit. 
The vertex algebra $\pi$ is the Heisenberg vertex algebra or equivalently the affine vertex algebra associated with $\gl_1$.
Recall that a vertex (operator) algebra $V$ is said to have strong generating type $\W(a_1,\dots,a_\ell)$ if $V$ has a set of generators $w_1,\dots,w_\ell$ of conformal weights $a_1,\dots,a_\ell$ respectively so that their differential monomials $\partial^{j_1}w_{i_1}\dots \partial^{j_s}w_{i_s}$ span $V$.

\begin{figure}[htbp]
    \centering
	\begin{minipage}[l]{0.3\linewidth}
		\begin{center} 
			\begin{tikzpicture}[node distance = 0pt,every node/.style = {draw, minimum size=7mm, inner sep=0pt, outer sep=0pt}]
\node (n11)  {$1$};
\node (n12) [right=of n11]  {$2$};
\node (n21) [below left=of n11.south]   {$3$};
\node (n22) [right=of n21]  {$4$};
\node (n23) [right=of n22]  {$5$};
\node (n31) [below left=of n21.south]   {$6$};
\node (n32) [right=of n31]              {$7$};
\node (n33) [right=of n32]              {$8$};
\node (n34) [right=of n33]              {$9$};
            \end{tikzpicture}
		\end{center}	
	\end{minipage}
	\begin{minipage}[c]{0.6\linewidth}
		\begin{center}
			$$\Gamma_x=
        \begin{tabular}{|cc|ccc|cccc|}
            \hline
             0&1 & -1/2&1/2&3/2 & -1&0&1&2 \\
             -1&0 & -3/2&-1/2&1/2 & -2&-1&0&1 \\\hline
             1/2&3/2 & 0&1&2 & -1/2&1/2&3/2&5/2 \\
             -1/2&1/2 & -1&0&1 & -3/2&-1/2&1/2&3/2 \\
             -3/2&-1/2 & -2&-1&0 & -5/2&-3/2&-1/2&1/2 \\\hline
             1&2 & 1/2&3/2&5/2 & 0&1&2&3 \\
             0&1 & -1/2&1/2&3/2 & -1&0&1&2 \\
             -1&0 & -3/2&-1/2&1/2 & -2&-1&0&1 \\
             -2&-1 & -5/2&-3/2&-1/2 & -3&-2&-1&0 \\
             \hline
        \end{tabular}$$
		\end{center}
	\end{minipage}
    \caption{Symmetric pyramid and grading associated with the partition $(2,3,4)$}
    \label{fig:sym_grading_gl9}
\end{figure}
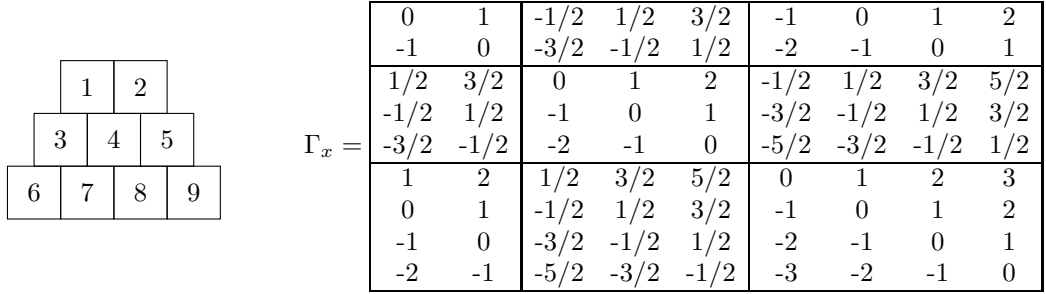

For $\lambda\vdash N$, the Jordan classification gives a natural choice of $f_\lambda$ to be
\begin{align}\label{nilpotent element}
    f_\lambda=f_{\lambda_1}+\dots + f_{\lambda_n},\quad 
\end{align}
where, in terms of the $N$-square elementary matrices $e_{i,j}$,
$$f_{\lambda_i}=e_{\mu_i+2,\mu_i+1}+\dots +e_{\mu_i+\lambda_i,\mu_i+\lambda_i-1},\qquad \mu_i=\lambda_1+\dots+\lambda_{i-1}.$$
Briefly, $f_{\lambda}$ is the matrix with $1$ on the subdiagonals of all its diagonal blocks.
For example, if $\lambda=(2,3,4)$ we have  
\begin{align*}
    f_{(2,3,4)}=f_{2}+f_3+f_4,\quad f_2=e_{2,1},\ f_3=e_{4,3}+e_{5,4},\ f_4=e_{7,6}+e_{8,7}+e_{9,8}.
\end{align*}
We may realize the grading $\Gamma_x$ by using the unique symmetric numbered pyramid associated with \eqref{nilpotent element}, see \cite{BG,EK} for the precise construction.
Then we have 
\begin{align}\label{semisimple element}
    x=x_{1}+\dots +x_{n}
\end{align}
where $x_{i}$ is the Weyl vector
\begin{align*}
    x_{i}=\frac{\lambda_i{-}1}{2}e_{\mu_i+1,\mu_i+1}+\frac{\lambda_i-{3}}{2}e_{\mu_i+2,\mu_i+2}+\dots +\frac{-(\lambda_i{-}1)}{2}e_{\mu_i+\lambda_i,\mu_i+\lambda_i}.
\end{align*}
In the case $f_\lambda=f_{(2,3,4)}$ of the above example, the pyramid and the good grading of $\sll_9$ (and $\gl_9$) are presented in Figure \ref{fig:sym_grading_gl9}.

By \cite{KW04}, the $\W$-algebra $\W^k(\gl_N,f_\lambda)$ has a PBW basis in terms of the strong generators represented by a homogeneous basis of the centralizer $\gl_N^{f_\lambda}=\{x\in\gl_N\mid [f_\lambda,x]=0\}$ with respect to $\Gamma_x$. Each strong generator $A$ corresponding to $a\in\gl_N^{f_\lambda}$ has conformal weight $-\Gamma_x(a)+1$.
We decompose $\gl_N$ by blocks of size $\lambda_i\times\lambda_j$ ($1\leq i,j\leq n$):
\begin{align*}
\gl_N=\bigoplus_{i,j=1}^n \gl_{N}^{i,j}.
\end{align*}
The blocks are separated by lines in the grading $\Gamma_x$ in Figure \ref{fig:sym_grading_gl9}.
It is straightforward to see that each block $\gl_{N}^{i,j}$ is stable under $\mathrm{ad}_{f_\lambda}$ and thus the centralizer of $f_\lambda$ decomposes into
$$\gl_N^{f_\lambda}=\bigoplus_{i,j}\gl^{f_\lambda}_{N;i,j},\quad \gl^{f_\lambda}_{N;i,j}:=\gl_N^{f_\lambda}\cap \gl_{N}^{i,j}.$$
Moreover, for each block $\gl_{N}^{i,j}$, the centralizer $\gl^{f_\lambda}_{N;i,j}$ has bases 
\begin{align*}
&{\sum_{0\leq \alpha < \lambda_j}e_{\mu_{i+1}-\alpha,\mu_{j+1}-\alpha}\ ,\ \sum_{1\leq \alpha < \lambda_j}e_{\mu_{i+1}+1-\alpha,\mu_{j+1}-\alpha} ,\dots, e_{\mu_{i+1},\mu_{j+1}-\lambda_j+1}},\quad (i\geq j),\\
&\sum_{0< \alpha \leq \lambda_i}e_{\mu_i+\alpha,\mu_j+\alpha}\ ,\ \sum_{0< \alpha \leq \lambda_i-1}e_{\mu_i+1+\alpha,\mu_j+\alpha} ,\dots, e_{\mu_{i}+\lambda_i,\mu_j+1} ,\quad (i\leq j),
\end{align*}
corresponding to strong generators in $\W^k(\gl_N,f_\lambda)$ with conformal weights
\begin{align}\label{coformal weights}
\frac{\lambda_j-\lambda_i}{2}+1,\dots,\frac{\lambda_j+\lambda_i}{2}
\end{align}
in both cases.
It is convenient to introduce the set
\begin{align*}
[m]_p=\{m,m+1,\dots,m+p-1 \},\quad(m\in\frac{1}{2}\Z)
\end{align*}
where $p\in\Z_{>0}\cup\{\infty\}$, to represent the sequence \eqref{coformal weights} as $[\frac{\lambda_j-\lambda_i}{2}+1]_{\lambda_i}$.
Then, the following is clear.
\begin{proposition}\label{strong generating type}
    The $\W$-algebra $\W^k(\gl_N,f_\lambda)$ has strong generating type 
    \begin{align*}
\W\left(\bigsqcup_{i,j=1}^n\left[\tfrac{|\lambda_j-\lambda_i|}{2}+1\right]_{\mathrm{min}\{\lambda_i,\lambda_j\}}\right).
\end{align*}
The strong generating type of the $\W$-algebra $\W^k(\sll_N,f_\lambda)$ is the same as the one of $\W^k(\gl_N,f_\lambda)$ after removing the strong generator of conformal weight $1$ corresponding to $\pi$.
\end{proposition}
\noindent For instance, for $f_{2,3,4}$, we have
\begin{align*}
    \W^k(\gl_9,f_{2,3,4})=\W\left(1^3,(3/2)^4,2^5,(5/2)^4,3^4,(7/2)^2,4\right)
\end{align*}
and
\begin{align*}
    \W^k(\sll_9,f_{2,3,4})=\W\left(1^2,(3/2)^4,2^5,(5/2)^4,3^4,(7/2)^2,4\right).
\end{align*}


\subsection{Reduction by stages}
In this section, we explain how we apply the reduction by stages. We then give evidence that the vertex algebra obtained after the successive reductions is a certain $\W$-algebra associated with $\sll_N$.
Consider $V$ a vertex algebra extension of $V^k(\g)$ with the following property:
\begin{itemize}
    \item[(P)] $V$ admits a $\frac{1}{2}\Z$-grading $V=\bigoplus_{\Delta}V_\Delta$ bounded from below such that $\dim V_\Delta <\infty$ for all $\Delta$ and compatible with the action of $V^k(\g)$, that is
    $a_{(n)}\colon V_\Delta\subset V_{\Delta-n}$ for all $a\in \g$ and $n\in\Z$. 
\end{itemize}
For a nilpotent element $f\in\g$, we define the BRST complex $C_{f}(V):=V\otimes \mathcal{F}(\g, f)$ where $\mathcal{F}(\g, f)$ and the differential $d$ are defined as in \S\ref{sec:BRST}.
We have the cohomology vanishing $H_{f}^{\neq0}(V)= 0$ by \cite{A00}.
Therefore, $H_{f}(V)\simeq H_{f}^0(V)$ is a vertex algebra extension of $\W^k(\g,{f})$ and has the $q$-character
\begin{align*}
    \ch[H_{f}(V)](q)=\mathrm{str}_{V\otimes \mathcal{F}(\g, {f})}q^{L_{{f},0}}.
\end{align*}

In the following, we want to consider $V$ to be the $\W$-algebra $\W^k(\g,f)$. Indeed, $\W$-algebras are extensions over the affine vertex subalgebra setting inside them. This affine vertex subalgebra is strongly generated by vectors represented by a basis of the centralizer $\g^\sharp:=\g_0\cap\g^f$. It is quite direct to see that $\W$-algebras satisfy the property (P).

Given a partition $\lambda=(\lambda_1\leq\dots \leq\lambda_n)\vdash N$, we set 
\begin{align*}
    \widehat{\lambda}_i:=(\underbrace{1,\dots,1}_{N_{i-1}},\lambda_i )\vdash N_i,\qquad N_i:=N-\sum_{j>i}{\lambda_j}.
\end{align*}
The partitions $\widehat{\lambda}_i$ are called hook-type as their Young tableaux look like hooks. 
For $i=1,\dots n$, the $\W$-algebra $\W^{k_i^\sharp}(\sll_{N_i},f_{\widehat{\lambda}_i})$ contains the affine vertex subalgebra $V^{k_{i-1}^\sharp}(\gl_{N_{i-1}})$
where the levels $k_i^\sharp$ are defined recursively by
\begin{align}\label{eq: precise levels}
  k_{i-1}^\sharp:=k_i^\sharp+\lambda_i-1,\quad k_{n}^\sharp:=k. 
\end{align}
As $\W^{k_i^\sharp}(\sll_{N_i},f_{\widehat{\lambda}_i})$ is a vertex algebra extension of $V^{k_{i-1}^\sharp}(\gl_{N_{i-1}})$ -- and thus of $V^{k_{i-1}^\sharp}(\sll_{N_{i-1}})$ -- satisfying (P), we can apply the BRST cohomology $H_{\widehat{\lambda}_{i-1}}(\cdot):=H_{f_{\widehat{\lambda}_{i-1}}}(\cdot)$. 
By induction, we define the reduction by stages 
\begin{align*}
\W_{[i]}^k(\sll_N,{f_\lambda}):=H_{f_{\widehat{\lambda}_{i}}}H_{f_{\widehat{\lambda}_{i+1}}}\dots\, H_{f_{\widehat{\lambda}_{n}}}\left(V^k(\sll_N)\right), \quad(i=1,\dots,n).
\end{align*}
Again we have the cohomology vanishing 
\begin{align*}
   \W_{[i]}^k(\sll_N,{f_\lambda})\simeq H^0_{f_{\widehat{\lambda}_{i}}}H^0_{f_{\widehat{\lambda}_{i+1}}}\dots H^0_{f_{\widehat{\lambda}_{n}}}\left(V^k(\sll_N)\right).
\end{align*}
Therefore, the $q$-character of $\W_{[i]}^k(\sll_N,{f_\lambda})$ is given by the Euler--Poincar\'{e} character of the total complex
\begin{align*}
 C_{f_\lambda}^{[i]}(V^k(\sll_N)):=V^k(\sll_N)\otimes \mathcal{F}_{{f_\lambda}}^{[i]}(V^k(\sll_N)),\quad  \mathcal{F}_{{f_\lambda}}^{[i]}(V^k(\sll_N)):=\bigotimes_{j=i}^n\mathcal{F}(\sll_{N_i},f_{\widehat{\lambda}_i}).
\end{align*}
For $i=1,\dots,n$ and $\spadesuit\in\{+,1/2\}$, $\Delta_{\spadesuit}^{x_i}$ is a set of roots in $\sll_{N_{i}}$ but 
the roots in $\Delta_{\spadesuit}^{x_i}$ can naturally be viewed as roots for $\sll_{N}$.
For $i\neq j$, the subsets $\Delta_{\spadesuit}^{x_i}$ and $\Delta_{\spadesuit}^{x_j}$ are disjoint and 
we set
\begin{align}\label{eq:delta_xbullet}
    \Delta_{\spadesuit}^{x_\bullet}=\Delta_{\spadesuit}^{x_1}\sqcup\dots \sqcup \Delta_{\spadesuit}^{x_n},\quad (\spadesuit\in\{+,1/2\}).
\end{align}
Indeed, realizing $\sll_{N}$ as the set of $N\!\times\!N$-matrices with zero trace as in \S\ref{sec:WslN}, the roots in $\Delta_+^{x_i}$ -- and so of $\Delta_{1/2}^{x_i}$ -- correspond to root vectors that are represented by certain elementary matrices $e_{m,n}$ ($1\leq m,n\leq N$) satisfying
\begin{itemize}
    \item $m<n$ and $N_{i-1}+1\leq n\leq N_i$ or
    \item $m>n$ and $N_{i-1}+1\leq m\leq N_i$.
\end{itemize}
The space $(\sll_{N_i})_{1/2}$ can also be viewed as a subspace of $\sll_N$ and we fix $\Phi[i]:=\Phi((\sll_{N_i})_{1/2})$. Then 
\begin{equation}\label{eq:partial_ghosts}
    \mathcal{F}(\sll_{N_i},f_{\widehat{\lambda}_i})=\Phi[i]\otimes\bigotimes_{\alpha\in\Delta^{x_i}_+}\mathcal{F}_\alpha[x_i].
\end{equation}
In particular, for $i=n$, we recover the definition of the usual complex appearing in the definition of the hook-type $\W$-algebra $\W^k(\sll_N,f_{\widehat{\lambda}_{n}})$, that is $C_{f_\lambda}^{[n]}(V^k(\sll_N))=C_{f_{\widehat{\lambda}_{n}}}(V^k(\sll_N))$, and for $i=n-1$, $C_{f_\lambda}^{[n-1]}(V^k(\sll_N))=C_{f_{\widehat{\lambda}_{n-1}}}(\W^{k}(\sll_N,f_{\widehat{\lambda}_{n}}))$.

Except for finitely many levels -- corresponding to the successive critical levels $k^\sharp_i\neq-h^\vee_{\sll_{N_i}}$ -- the vertex algebra $\W_{[i]}^k(\sll_N,{f_\lambda})$ admits a conformal vector $L^{[i]}_{{f_\lambda}}$ which is defined inductively. 
We firstly set $L^{[n]}_{{f_\lambda}}:=L_{f_{\widehat{\lambda}_n}}$ as defined in \eqref{KRW}.
For $i=1,\dots,n$, since the conformal vector $L_{f_{\widehat{\lambda}_i}}$ of $\W^{k_i^\sharp}(\sll_{N_i},f_{\widehat{\lambda}_i})$ satisfies the $\Lambda$-brackets
\begin{align*}
    \lbr{L_{f_{\widehat{\lambda}_n}}}{a}=\partial a+a\Lambda,\quad (a\in \gl_{N_{i-1}}\subset V^{k_{i-1}^\sharp}(\gl_{N_{i-1}})),
\end{align*}
by \cite[Thm 2.4]{KW04}, we can decompose 
\begin{align*}
L_{f_{\widehat{\lambda}_i}}=L_{f_{\widehat{\lambda}_i}, \perp}+L_{\mathrm{sug}}^{\sll_{N_{i-1}}},
\end{align*}
With this decomposition, the field $L_{f_{\widehat{\lambda}_i}, \perp}$ satisfies the Virasoro commutation relations and commutes with $V^{k_{i-1}^\sharp}(\sll_{N_{i-1}})$, and thus $\lbr{L_{f_{\widehat{\lambda}_i}, \perp}}{L_{\mathrm{sug}}^{\sll_{N_{i-1}}}}=0$.
Accordingly, one has the decomposition 
\begin{align}\label{eq:decompo_conf_vector}
    L_{f_\lambda}^{[i]}=L^{[i]}_{{f_\lambda},\perp}+L_{\mathrm{sug}}^{\sll_{N_{i-1}}},\quad [L^{[i]}_{{f_\lambda},\perp}{}_\Lambda L_{\mathrm{sug}}^{\sll_{N_{i-1}}}]=0.
\end{align}
Then we define
\begin{equation*}
    L^{[i-1]}_{f_\lambda}:=L^{[i]}_{{f_\lambda},\perp}+L_{f_{\widehat{\lambda}_{i-1}}},\quad (i=2,\dots,n).
\end{equation*}

As in \S\ref{sec:BRST}, we define similarly the complexes $C_{f_\lambda}^{[i]}(M)$ for any $V^k(\g)$-module $M$ satisfying the property $(P)$. The reduction by stages
$H_{f_{\widehat{\lambda}_{1}},f_{\widehat{\lambda}_{2}},\dots,f_{\widehat{\lambda}_{n}}}(M)$ is a module over $H_{f_{\widehat{\lambda}_{1}},f_{\widehat{\lambda}_{2}},\dots,f_{\widehat{\lambda}_{n}}}(V^k(\sll_N))$.

\begin{theorem}\label{Conj at the level of q}
The vertex algebras $\W^k(\sll_N,f_{\lambda})$ and  
$H_{f_{\widehat{\lambda}_{1}},f_{\widehat{\lambda}_{2}},\dots,f_{\widehat{\lambda}_{n}}}(V^k(\sll_N))$ have the same central charge and the $q$-characters. Moreover, the equality of the $q$-characters 
$$\ch[H_{f_\lambda}(M)](q) =\ch[H_{f_{\widehat{\lambda}_{1}},f_{\widehat{\lambda}_{2}},\dots,f_{\widehat{\lambda}_{n}}}(M)](q).$$ 
holds for any $V^k(\g)$-module $M$ satisfying the property $(P)$.
\end{theorem}

\subsection{Proof of Theorem \ref{Conj at the level of q}}\label{sec: Proof of Conj A at the level of q}


We will compare the complexes $C_{f_\lambda}(M)$ and $C_{f_\lambda}^{[1]}(M)$ and the structure defined by the respective conformal vectors.
To this end, we derive an explicit formula of the conformal vector $L^{[1]}_{{f_\lambda}}$ of $\W_{[1]}^k(\sll_N,{f_\lambda})$ by induction. 
We consider first the case $\W_{[n-1]}^k(\sll_N,{f_\lambda})$.
By \cite[Thm 2.4]{KRW03}, the subalgebra $V^{k^\sharp_{n-1}}(\sll_{N_{n-1}})\subset \W_{[n-1]}^k(\sll_N,{f_\lambda})$ is generated by the fields in $C_{{f_\lambda}}^{[n]}(V^k(\sll_N))$
\begin{align}\label{wide a first case}
    \widehat{a}=a+a^{\Phi[n]} +\sum_{\beta,\gamma\in \Delta_{+}^{x_n}} c_{a,\beta}^\gamma \varphi_\gamma \varphi_\beta^*,\quad (a\in \sll_{N_{n-1}}),
\end{align}
where $\Phi[i]$ is defined in \eqref{eq:partial_ghosts} and $a^{\Phi[n]}$ in \eqref{eq:phi_xa}. Using \eqref{KRW} and \eqref{eq:decompo_conf_vector}, we have
\begin{align*}
    &L_{f_\lambda}^{[n-1]}
    =L^{[n]}_{{f_\lambda},\perp}+ L_{f_{\widehat{\lambda}_{n-1}}}
    =\left(L_{f_{\widehat{\lambda}_{n}}}-L_{\mathrm{sug}}^{\sll_{N_{n-1}}}\right)+ L_{f_{\widehat{\lambda}_{n-1}}}\\
    &=\left( L^{\sll_{N}}_{\mathrm{sug}} + \partial x_{n}+ \sum_{\alpha\in \Delta_{+}^{x_n}}L_\alpha^{x_n}+L_{\Phi[n]}^{x_n}- L_{\mathrm{sug}}^{\sll_{N_{n-1}}}\right)+\left(L_{\mathrm{sug}}^{\sll_{N_{n-1}}}+\partial \widehat{x}_{n-1}+\sum_{\alpha\in \Delta_{+}^{x_{n-1}}}L_\alpha^{x_{n-1}}+L_{\Phi[n-1]}^{x_{n-1}} \right)\\
    &= L^{\sll_{N}}_{\mathrm{sug}} + \partial (x_{n-1}+x_{n})+L_{\Phi[n-1]}^{x_{n-1}}+L_{\Phi[n]}^{x_{n-1}+x_n}+ \sum_{\alpha\in \Delta_{+}^{x_{n-1}}\sqcup\Delta_{+}^{x_n}}L_\alpha^{x_{n-1}+x_{n}}.
\end{align*}
In the last equation, we used that $L_{\Phi[n]}^{x_n}+\partial x_{n-1}^{\Phi[n]}=L_{\Phi[n]}^{x_{n-1}+x_n}$ (see \eqref{eq:phi_xa}). Moreover, for $\alpha\in \Delta_{+}^{x_{n-1}}$, $\alpha(x_{n-1}+x_{n})=\alpha(x_{n-1})$ by restriction, hence $L_\alpha^{x_{n-1}}=L_\alpha^{x_{n-1}+x_{n}}$.
It is straightforward to generalize the formula \eqref{wide a first case} and to show by induction that the subalgebra $V^{k^\sharp_{i}}(\sll_{N_{i}})\subset \W_{[i+1]}^k(\g,{f_\lambda})$ is generated by the fields in $C_{{f_\lambda}}^{[i+1]}(V^k(\sll_N))$
\begin{align*}
    \widehat{a}=a+\sum_{j=i+1}^na^{\Phi[j]} +\sum_{j=i+1}^n\sum_{\beta,\gamma\in \Delta_{+}^{x_j}} c_{a,\beta}^\gamma \varphi_\gamma \varphi_\beta^*,\quad (a \in \sll_{N_{i}})
\end{align*}
Hence, we have
\begin{equation}\label{eq:L1f}
\begin{aligned}
    L^{[1]}_{f_\lambda}
    &= L^{\sll_{N}}_{\mathrm{sug}} + \partial(x_1+\dots+ x_{n})+L_{\Phi[1]}^{x_1}+L_{\Phi[2]}^{x_1+x_2}+\dots +L_{\Phi[n]}^{x_1+\dots+ x_n}  \\
    &\hspace{1cm}+\sum_{\alpha\in \Delta_{+}^{x_{1}}}L_\alpha^{x_{1}}+\sum_{\alpha\in \Delta_{+}^{x_{2}}}L_\alpha^{x_{1}+x_2}+\dots+ \sum_{\alpha\in \Delta_{+}^{x_n}}L_\alpha^{x_1+x_2+\dots +x_{n}}\\
    &= L^{\sll_{N}}_{\mathrm{sug}}+\partial x+ \sum_{i=1}^n L_{\Phi[i]}^x+ \sum_{\alpha\in \Delta^{x_\bullet}_+}L_\alpha^x
\end{aligned} 
\end{equation}
where $\Delta^{x_\bullet}_+$ was defined in \eqref{eq:delta_xbullet}. 

It suffices to consider the central charges and the $q$-characters of the ghost parts of the complexes $C_{f_\lambda}(M)$ and $C_{f_\lambda}^{[1]}(M)$, namely $\mathcal{F}(\g, f_\lambda)$ and $\mathcal{F}^{[1]}(\g, f_\lambda)$.
The decomposition of $\sll_N$ into blocks of size $\lambda_i\times\lambda_j$ ($1\leq i,j\leq n$) induces the decomposition 
\begin{align*}
    \Delta_+^{x}=\bigsqcup_{i,j=1}^n\Delta_+^{x}[i,j],\quad  \Delta_+^{x_\bullet}=\bigsqcup_{i,j=1}^n\Delta_+^{x_\bullet}[i,j],
\end{align*}
and similarly for $\Delta_{1/2}^{x}$ and $\Delta_{1/2}^{x_\bullet}$.
Accordingly, we decompose the ghost parts of the complexes
\begin{align*}
&\mathcal{F}(\g, f_\lambda)\simeq \bigotimes_{i,j=1}^n \mathcal{F}_{i,j},& \mathcal{F}_{i,j}=\bigotimes_{\alpha\in \Delta_+^x[i,j]} \mathcal{F}_\alpha^x \otimes \bigotimes_{\alpha\in \Delta_{1/2}^x[i,j]} \langle \Phi_\alpha \rangle,\\
&\mathcal{F}^{[1]}(\g, f_\lambda)\simeq \bigotimes_{i,j=1}^n\mathcal{F}^{[1]}_{i,j} & \mathcal{F}^{[1]}_{i,j}=\bigotimes_{\alpha\in \Delta_+^{x_\bullet}[i,j]} \mathcal{F}_\alpha^x \otimes \bigotimes_{\alpha\in \Delta_{1/2}^{x_\bullet}[i,j]} \langle \Phi_\alpha \rangle,
\end{align*}
where $\langle \Phi_\alpha \rangle$ denotes the subalgebra generated by the field $\Phi_\alpha$.
We compute the action of the conformal vectors $L_{f_\lambda}$ and $L_{f_\lambda}^{[1]}$ restricted to the tensor product components for the $(i,j)$-block, that are $\mathcal{F}_{i,j}$ and $\mathcal{F}^{[1]}_{i,j}$ respectively.
First, by \eqref{KRW}, we have
\begin{align*}
    (L_{f_\lambda})_{| \mathcal{F}_{i,j}}= \sum_{\alpha\in \Delta_+^{x}[i,j]}(L^x_\alpha)_{| \mathcal{F}_{i,j}}+(L^x_\Phi)_{| \mathcal{F}_{i,j}}.
\end{align*}
Note that for $\alpha\in\Delta_+^{x}[i,j]$, $\alpha(x)=\alpha(x_i+x_j)$ if $i\neq j$ and $\alpha(x)=\alpha(x_i)$ if $i=j$. Then it follows from \eqref{eq:alphax} that $(L^{x}_\alpha)_{| \mathcal{F}_{i,j}}=(L^{x_i+x_j}_\alpha)_{| \mathcal{F}_{i,j}}$.
Similarly, using \eqref{eq:phi_xa}, we obtain
\begin{equation*}
    (L^x_\Phi)_{|\mathcal{F}_{i,j}}=(L^{x_i+x_j}_\Phi)_{| \mathcal{F}_{i,j}}-\frac{1}{2}\sum_{\substack{1\leq m\leq n\\m\neq i,j}}\sum_{\alpha\in\Delta^x_{1/2}[i,j]}\underbrace{\alpha(x_m)}_{0}\partial\left(\Phi_\alpha\Phi_\alpha^*\right)_{| \mathcal{F}_{i,j}}.
\end{equation*}
Hence, $(L_{f_\lambda})_{| \mathcal{F}_{i,j}}=(L_{\Phi}^{x_i+x_j}+\sum_{\alpha\in \Delta_{+}^x[i,j]}L_\alpha^{x_i+x_j})_{| \mathcal{F}_{i,j}}$.
Applying the same argument to the decomposition \eqref{eq:L1f}, we have
\begin{align*}
    (L_{f_\lambda}^{[1]})_{|\mathcal{F}^{[1]}_{i,j}}=\left(\sum_{m=1}^n L_{\Phi[m]}^{x_i+x_j}+\sum_{\alpha\in \Delta^{x_\bullet}_+[i,j]}L_\alpha^{x_i+x_j}\right)_{|\mathcal{F}^{[1]}_{i,j}}=\left(L_{\Phi[i]}^{x_i+x_j}+\sum_{\alpha\in \Delta^{x_\bullet}_+[i,j]}L_\alpha^{x_i+x_j}\right)_{|\mathcal{F}^{[1]}_{i,j}}.
\end{align*}
It follows that the two actions restricted to the $(i,j)$-block depend only on $\lambda_i$ and $\lambda_j$.
Thus the proof reduces to the case $n=2$, corresponding to the partition $(\lambda_i, \lambda_j)$.

We show the case $\lambda=(\lambda_1,\lambda_2)$ with $x=x_1+x_2$.
By \eqref{cc for charged fermion} and \eqref{cc for neutral fermion}, the central charges of $(L_{f_\lambda})_{|\mathcal{F}_{i,j}}$ and $(L_{f_\lambda}^{[1]})_{|\mathcal{F}^{[1]}_{i,j}}$ are given respectively by
\begin{equation}\label{eq:central_charges}
    \frac{1}{2}\sum_{\alpha\in \Delta^x_{1/2}[i,j]}c[\alpha(x)]\ {-}\sum_{\alpha\in \Delta_+^x[i,j]} c[\alpha(x)],\quad
    \frac{1}{2}\sum_{\alpha\in \Delta^{x_\bullet}_{1/2}[i,j]}c[\alpha(x)]\ {-}\sum_{\alpha\in \Delta_+^{x_\bullet}[i,j]}c[\alpha(x)]
\end{equation}
for $i,j\in \{1,2\}$, . The $q$-characters of $\mathcal{F}_{i,j}$ and $\mathcal{F}^{[1]}_{i,j}$ are given by 
\begin{equation}\label{eq:characters}
    \frac{\prod_{\alpha\in \Delta_+^x[i,j]}(q^{1-\alpha(x)},q^{\alpha(x)};q)_\infty}{\prod_{\alpha\in \Delta_{1/2}^x[i,j]} (q^{1-\alpha(x)};q)_\infty},\quad 
    \frac{\prod_{\alpha\in \Delta_+^{x_\bullet}[i,j]}(q^{1-\alpha(x)},q^{\alpha(x)};q)_\infty}{\prod_{\alpha\in \Delta_{1/2}^{x_\bullet}[i,j]} (q^{1-\alpha(x)};q)_\infty},
\end{equation}
respectively, where $(a_1,\dots,a_n;q)_\infty$ is the $q$-Pochhammer symbol:
\begin{equation*}
    (a_1,\dots,a_n;q)_\infty=\prod_{i=1}^n(a_i;q)_\infty=\prod_{i=1}^n\prod_{k\geq0}(1-a_i q^k).
\end{equation*}

When $i=j$, $\Delta_+^x[i,j]=\Delta_+^{x_\bullet}[i,j]$ and $\Delta_{1/2}^x[i,j]=\Delta_{1/2}^{x_\bullet}[i,j]$ so the central charges and the $q$-characters of $(L_{f_\lambda})_{|\mathcal{F}_{i,j}}$ and $(L_{f_\lambda}^{[1]})_{|\mathcal{F}^{[1]}_{i,j}}$ are identical. 
Assume $i\neq j$, then swapping the role of $i$ and $j$, we can restrict to the case $(i,j)=(1,2)$.
We consider this case in the following.
Then $\Delta_{\spadesuit}^{x_\bullet}{[1,2]}=\Delta_{\spadesuit}^{x_2}{[1,2]}$ for $\spadesuit=+,1/2$.

\begin{figure}[htbp]
    \centering
\begin{align*}\label{two gradings}
\Gamma_x\colon {
        \begin{tabular}{|cccc|}
            \hline
              $\frac{a-b}{2}$ & $\dots$ &$\frac{a+b-4}{2}$ &$\frac{a+b-2}{2}$ \\
              $\vdots$ & \ &\ &$\frac{a+b-4}{2}$ \\
              $\frac{-a-b+4}{2}$ & \ &\ &$\vdots$ \\
              $\frac{-a-b+2}{2}$ & $\frac{-a-b+4}{2}$ &$\dots$ &$\frac{-a+b}{2}$ \\
             \hline
        \end{tabular}}\ ,\quad 
\Gamma_{x_2}\colon 
        \begin{tabular}{|cccc|}
            \hline
              \hspace{2mm}$\frac{-b+1}{2}$\hspace{2mm} & $\dots$ & \hspace{4mm}$\frac{b-3}{2}$ &\hspace{4mm}$\frac{b-1}{2}$\hspace{2mm} \\
              \hspace{2mm}$\frac{-b+1}{2}$\hspace{2mm} & \ &\hspace{4mm}$\frac{b-3}{2}$ &\hspace{4mm}$\frac{b-1}{2}$\hspace{2mm} \\
              $\vdots$ & \ &\hspace{4mm}$\vdots$ &\hspace{4mm}$\vdots$\hspace{2mm} \\
              \hspace{2mm}$\frac{-b+1}{2}$\hspace{2mm} & $\dots$ &\hspace{4mm}$\frac{b-3}{2}$ &\hspace{4mm}$\frac{b-1}{2}$\hspace{2mm} \\
             \hline
        \end{tabular}
\end{align*}
\begin{align*}
\Delta_{+}^{x}[1,2]\colon {
        \begin{tabular}{|cccccccc|}
            \hline
              \ & \ &$\ast$&$\ast$ &$\dots$ &$\ast$&$\dots$ &$\ast$ \\
              \ & \ &\ &$\ast$ &$\ddots$&$\vdots$ &\ &$\vdots$ \\
             \ &\ & \ & \ &$\ddots$&$\ast$&$\dots$&$\ast$ \\
             \ &\ & \ & \ & \ &$\ast$&$\dots$ &$\ast$ \\
             \hline
        \end{tabular}}\ ,\quad 
\Delta_{+}^{x_2}[1,2]\colon {
        \begin{tabular}{|cccccccc|}
            \hline
              \hspace{4mm} & \ &\ &\ &$\ast$ &$\dots$&$\ast$ &$\ast$ \\
              \ & \ &\ &\ &$\vdots$&\ &$\vdots$ &$\vdots$ \\
             \ &\ & \ & \ &$\ast$&$\dots$ &$\ast$ & $\ast$ \\
             \ &\ & \ & \ & $\ast$ &$\dots$&$\ast$ &$\ast$ \\
             \hline
        \end{tabular}}        
\end{align*}
    \caption{Grading and positives roots for $x$ and $x_2$.}
    \label{fig:gradings_roots}
\end{figure}

Remark that in the $q$-characters, two factors on the denominator cancel with one factor on the numerator corresponding to a root in $\Delta_{1/2}^{\clubsuit}[1,2]$ ($\clubsuit=x,x_2$). Roughly speaking, two roots in $\Delta_{1/2}^{\clubsuit}[1,2]$ in the denominator ``kill'' one root in $\Delta_{+}^{\clubsuit}[1,2]$ in the numerator.
Similar cancellation happens for the central charges.
Thus, we will construct the sets $\overline{\Delta}_{\spadesuit}^{\clubsuit}[1,2]$ removing certain roots of grading $\frac{1}{2}$ in the sets $\Delta_{\spadesuit}^{\clubsuit}[1,2]$ ($\spadesuit=+,1/2$) to avoid the redundancy in the computations. 
The basic idea is to remove a maximal even number of roots in ${\Delta}_{1/2}^{\clubsuit}[1,2]$, says $2m$, so that $\overline{\Delta}_{1/2}^{\clubsuit}[1,2]$ has cardinal $0$ or $1$, and remove $m$ roots of grading $\frac{1}{2}$ in ${\Delta}_{+}^{\clubsuit}[1,2]$.
The precise definition of the resulting sets $\overline{\Delta}_{\spadesuit}^{\clubsuit}[1,2]$ depends on the parity of $\lambda_1,\lambda_2$ and is explained below.
We write the symbol $\xrightarrow[\mathrm{q.isom}]{\sim}$ to indicate this restriction:
\begin{align*}
\Delta_+^{\clubsuit}[1,2]\sqcup \Delta_{1/2}^{\clubsuit}[1,2]\xrightarrow[\mathrm{q.isom}]{\sim} \overline{\Delta}_+^{\clubsuit}[1,2]\sqcup\overline{\Delta}_{1/2}^{\clubsuit}[1,2].
\end{align*}

Then we construct an identification 
\begin{align*}
    \psi\colon \overline{\Delta}_+^{x}[1,2]\sqcup\overline{\Delta}_{1/2}^{x}[1,2] \xrightarrow{\sim} \overline{\Delta}_+^{x_2}[1,2]\sqcup\overline{\Delta}_{1/2}^{x_2}[1,2],
\end{align*}
satisfying $\psi(\alpha)(x)=1-\alpha(x)$ for all $\alpha$ such that $\psi(\alpha)\neq \alpha$.
This transformation preserves the $q$-characters obviously as well as the central charges since $c[1-t]=c[t]$ for all $t$.

Below, we describe these constructions case by case depending on $(\lambda_1,\lambda_2)=(a,b)$ modulo 2.
The gradings $\Gamma_x$ and $\Gamma_{x_2}$ on the $(1,2)$-block and the roots in $\Delta_{+}^{x}[1,2]$ and $\Delta_{+}^{x_2}[1,2]$ appear in Figure~\ref{fig:gradings_roots}.
Note, when they exist, roots in $\Delta_{1/2}^{x}[1,2]$ and $\Delta_{1/2}^{x_2}[1,2]$ (and thus in the subsets $\overline{\Delta}_{1/2}^{x}[1,2]$, $\overline{\Delta}_{1/2}^{x_2}[1,2]$) only appear along the left boundary. 
In this case, they are indicated by $\overline{\ast}$. 
If ``$\xrightarrow[\mathrm{q.isom}]{\sim}$" is omitted, we do not need to take reduced sets.

\begin{itemize}
    \item $(a,b)\equiv (1,1)$
        
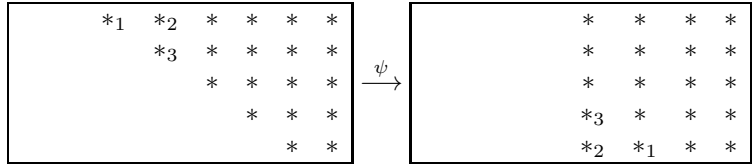
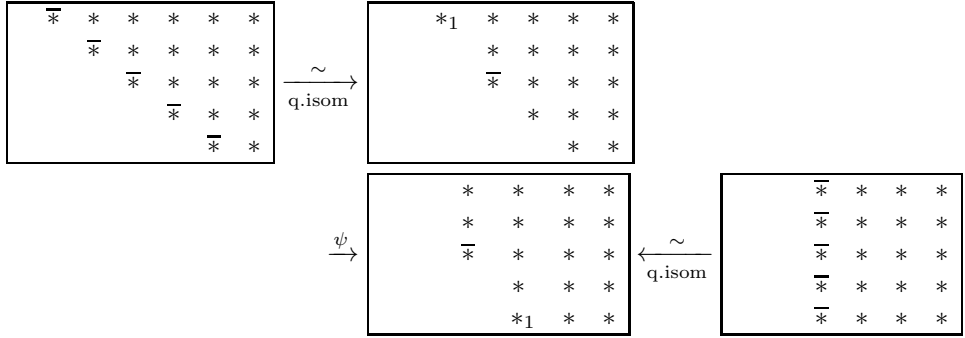
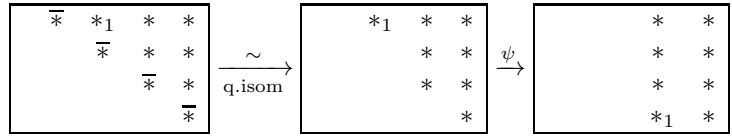
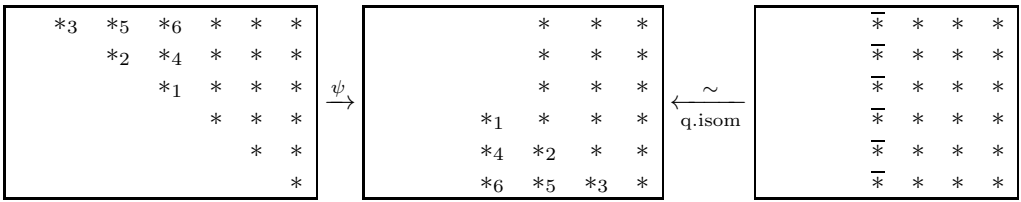
\begin{align*}
        \begin{tabular}{|ccccccccc|}
            \hline
             \ & \ & \ &$*_1$&$*_2$ &$\ast$&$\ast$&$\ast$&$\ast$ \\
             \ & \ & \ &\ &$*_3$ &$\ast$&$\ast$&$\ast$&$\ast$ \\
             \ &\ &\ & \ & \ &$\ast$&$\ast$&$\ast$&$\ast$ \\
             \ &\ &\ & \ & \ & \ &$\ast$&$\ast$ &$\ast$ \\
             \ &\ &\ & \ & \ & \  & &$\ast$&$\ast$  \\
             \hline
        \end{tabular} \overset{\psi}{\longrightarrow}
        \begin{tabular}{|ccccccccc|}
            \hline
            \hspace{3mm} & \ & \ &\ & \ &$\ast_{\hspace{1mm}}$&$\ast_{\hspace{1mm}}$&$\ast$&$\ast$ \\
             \ & \ & \ &\ &\ &$\ast_{\hspace{1mm}}$&$\ast_{\hspace{1mm}}$&$\ast$&$\ast$ \\
             \ &\ &\ & \ & \ &$\ast_{\hspace{1mm}}$& $\ast_{\hspace{1mm}}$&$\ast$&$\ast$ \\
             \ &\ &\ & \ & \ & $*_3$ &$\ast_{\hspace{1mm}}$&$\ast$ &$\ast$ \\
             \ &\ &\ & \ & \ & $*_2$ &$*_1$ &$\ast$&$\ast$  \\
             \hline
        \end{tabular}        
    \end{align*}
    \item $(a,b)\equiv (1,0)$
        \begin{align*}
        \begin{tabular}{|ccccccc|}
            \hline
             \ &$\overline{\ast}$&$\ast$ &$\ast$&$\ast$&$\ast$&$\ast$ \\
              \ &\ &$\overline{\ast}$ &$\ast$&$\ast$&$\ast$&$\ast$ \\
             \ & \ & \ &$\overline{\ast}$&$\ast$&$\ast$&$\ast$ \\
             \ & \ & \ & \ &$\overline{\ast}$&$\ast$ &$\ast$ \\
             \ & \ & \ & \  & &$\overline{\ast}$&$\ast$  \\
             \hline
        \end{tabular}\xrightarrow[\mathrm{q.isom}]{\sim}\
                &\begin{tabular}{|ccccccc|}
            \hline
               \ &\ &$*_1$ &$\ast$&$\ast$&$\ast$&$\ast$ \\
              \ &\ &\ &$\ast$&$\ast$&$\ast$&$\ast$ \\
             \ & \ & \ &$\overline{\ast}$&$\ast$&$\ast$&$\ast$ \\
             \ & \ & \ & \ &$\ast$&$\ast$ &$\ast$ \\
             \ & \ & \ & \  & &$\ast$&$\ast$  \\
             \hline
        \end{tabular}\\
        \xrightarrow{\psi}\ &
        \begin{tabular}{|ccccccc|}
            \hline
            \ &\ & \ &$\ast_{\hspace{1mm}}$&$\ast_{\hspace{1mm}}$&$\ast$&$\ast$ \\
               \ &\ &\ &$\ast_{\hspace{1mm}}$&$\ast_{\hspace{1mm}}$&$\ast$&$\ast$ \\
             \ & \ & \ &$\overline{\ast}_{\hspace{1mm}}$& $\ast_{\hspace{1mm}}$&$\ast$&$\ast$ \\
             \ & \ & \ & \ &$\ast_{\hspace{1mm}}$&$\ast$ &$\ast$ \\
             \ & \ & \ & \ &$*_1$ &$\ast$&$\ast$  \\
             \hline
        \end{tabular}\xleftarrow[\mathrm{q.isom}]{\sim} 
                \begin{tabular}{|ccccccc|}
            \hline
             \ &\ & \ &$\overline{\ast}$&$\ast$&$\ast$&$\ast$ \\
               \ &\ &\ &$\overline{\ast}$&$\ast$&$\ast$&$\ast$ \\
             \ & \ & \ &$\overline{\ast}$& $\ast$&$\ast$&$\ast$ \\
             \ & \ & \ & $\overline{*}$ &$\ast$&$\ast$ &$\ast$ \\
             \ & \ & \ & $\overline{*}$ &$\ast$ &$\ast$&$\ast$  \\
             \hline
        \end{tabular}
    \end{align*}
    \item $(a,b)\equiv (0,1)$
\begin{align*}
        \begin{tabular}{|ccccc|}
            \hline
            \ &$\overline{*}$&$*_1$ &$\ast$&$\ast$ \\
            \ &\ &$\overline{*}$ &$\ast$&$\ast$ \\ 
            \ & \ & \ &$\overline{\ast}$&$\ast$ \\
            \ & \ & \ & \ &$\overline{\ast}$\\
             \hline
        \end{tabular} \xrightarrow[\mathrm{q.isom}]{\sim}
        \begin{tabular}{|ccccc|}
            \hline
            \ &\ &$*_1$ &$\ast$&$\ast$ \\
            \ &\ &\  &$\ast$&$\ast$ \\ 
            \ & \ & \ &$\ast$&$\ast$ \\
            \ & \ & \ & \ &$\ast$\\
             \hline
        \end{tabular}
        \xrightarrow{\psi}
       \begin{tabular}{|ccccc|}
            \hline
            \hspace{3mm} & \ & \ &$\ast_{\hspace{1mm}}$&$\ast$ \\
            \ &\ & \ &$\ast_{\hspace{1mm}}$&$\ast$ \\ 
            \ & \ & \ &$\ast_{\hspace{1mm}}$&$\ast$ \\
            \ & \ & \ & $*_1$ &$\ast$\\
             \hline
        \end{tabular}        
    \end{align*}
    
\item $(a,b)\equiv (0,0)$
    \begin{align*}
        \begin{tabular}{|ccccccc|}
            \hline
            \ &$*_3$ &$*_5$ &$*_6$ &$\ast$ &$\ast$ &$\ast$ \\
           \ &\ &$*_2$ &$*_4$ &$\ast$ &$\ast$ &$\ast$ \\
            \ \ &\ &\ &$*_1$ &$\ast$ &$\ast$ &$\ast$ \\
            \ &\ &\ &\  &$\ast$ &$\ast$ &$\ast$ \\
            \  &\ &\ &\ &\ &$\ast$ &$\ast$ \\
            \ &\ &\ &\ &\ &\ &$\ast$ \\
             \hline
        \end{tabular} \xrightarrow{\psi} 
       \begin{tabular}{|ccccccc|}
            \hline
            \hspace{3mm}  & \ & \ & \ &$\ast$ &$\ast$ &$\ast$ \\
             \ & \ & \ & \ &$\ast$ &$\ast$ &$\ast$ \\
             \ & \ & \ & \ &$\ast$ &$\ast$ &$\ast$ \\
             \ & \ & \ & $*_1$ &$\ast$ &$\ast$ &$\ast$ \\
             \ & \ & \ & $*_4$ &$*_2$ &$\ast$ &$\ast$ \\
             \ & \ & \ & $*_6$ &$*_5$ &$*_3$ &$\ast$ \\
             \hline
        \end{tabular} 
        \xleftarrow[\mathrm{q.isom}]{\sim}
       \begin{tabular}{|ccccccc|}
            \hline
            \hspace{3mm} &  \ & \ & $\overline{\ast}$ &$\ast$ &$\ast$ &$\ast$ \\
            \ & \ & \ & $\overline{\ast}$ &$\ast$ &$\ast$ &$\ast$ \\
             \ & \ & \ & $\overline{\ast}$ &$\ast$ &$\ast$ &$\ast$ \\
             \ & \ & \ & $\overline{\ast}$ &$\ast$ &$\ast$ &$\ast$ \\
            \ & \ & \ & $\overline{\ast}$ &$\ast$ &$\ast$ &$\ast$ \\
             \ & \ & \ & $\overline{\ast}$ &$\ast$ &$\ast$ &$\ast$ \\
             \hline
        \end{tabular}        
    \end{align*}
\end{itemize}
For each block, the identification of the roots guarantees the equality of the central charges \eqref{eq:central_charges} and $q$-characters \eqref{eq:characters}.
This completes the proof.

\bigskip
For readers' convenience, we detail the construction of the reduced sets for $\lambda=(5,6)$, which corresponds to the case $(a,b)\equiv(1,0)$.
The gradings $\Gamma_x$ and $\Gamma_{x_2}$ restricted to the $(1,2)$-block are given by
        $$\Gamma_x=
        \begin{tabular}{|cccccc|}
            \hline
             -1/2&1/2&3/2&5/2&7/2&9/2 \\
             -3/2&-1/2&1/2&3/2&5/2&7/2 \\
             -5/2&-3/2&-1/2&1/2&3/2&5/2 \\
             -7/2&-5/2&-3/2&-1/2&1/2&3/2 \\
             -9/2&-7/2&-5/2&-3/2&-1/2&1/2 \\
             \hline
        \end{tabular}\ ,\qquad
        \Gamma_{x_2}=
        \begin{tabular}{|cccccc|}
            \hline
             -5/2&-3/2&-1/2&1/2&3/2&5/2\\
             -5/2&-3/2&-1/2&1/2&3/2&5/2 \\
             -5/2&-3/2&-1/2&1/2&3/2&5/2 \\
             -5/2&-3/2&-1/2&1/2&3/2&5/2 \\
             -5/2&-3/2&-1/2&1/2&3/2&5/2 \\
             \hline
        \end{tabular}\ .$$
Accordingly, we have
\begin{align*}
  &\Delta_+^x[1,2]=\{\alpha_{i,j},\, 1\leq i\leq 10,\, i+5\leq j\leq 10\},
  &&\Delta_{1/2}^x[1,2]=\{\alpha_{i,i+5},\, 1\leq i\leq 5\},\\
  &\Delta_+^{x_{\bullet}}[1,2]=\{\alpha_{i,j},\, 1\leq i\leq 5,\, 8\leq j\leq 10\},
  &&\Delta_{1/2}^{x_{\bullet}}[1,2]=\{\alpha_{i,8},\, 1\leq i\leq 5\},
\end{align*}
where $\alpha_{i,j}=\alpha_i+\dots+\alpha_j$ with $i\leq j$.
We take the reduced sets
\begin{align*}
    &\overline{\Delta}_+^{x}[1,2]=\Delta_+^x[1,2]\backslash\{\alpha_{1,6},\alpha_{2,7}\},
    &\overline{\Delta}_{1/2}^{x}[1,2]=\{\alpha_{3,8}\},\\
    &\overline{\Delta}_+^{x_\bullet}[1,2]=\Delta_+^{x_{\bullet}}[1,2]\backslash\{\alpha_{4,8},\alpha_{5,8}\},
    &\overline{\Delta}_{1/2}^{x_\bullet}[1,2]=\{\alpha_{3,8}\}
\end{align*}
where we remove from ${\Delta}_+^{x}[1,2]$, the roots in ${\Delta}_{1/2}^{x}[1,2]$ which are not in ${\Delta}^{x_\bullet}_+[1,2]$ and conversely.

Comparing the reduced sets, we note that the only root of $\overline{\Delta}_+^{x}[1,2]$ which is not in $\overline{\Delta}_+^{x_\bullet}[1,2]$ is $\alpha_{1,7}$ whereas $\alpha_{5,9}$ is the only root in $\overline{\Delta}_+^{x_\bullet}[1,2]$ not in $\overline{\Delta}_+^{x}[1,2]$.
Their corresponding root vectors $e_{1,8}$ and $e_{5,10}$ have respective grading $\frac{3}{2}$ and $-\frac{1}{2}=1-\frac{3}{2}$ in $\Gamma_{x_2}$. Hence, we swap these two roots.
As a consequence, we define the isomorphism $\psi$ such that
\begin{equation*}
    \psi(\alpha)=\left\{\begin{aligned}
        \alpha_{5,9},&\quad\text{if }\alpha=\alpha_{1,7},\\
        \alpha,&\quad\text{otherwise}.
    \end{aligned}\right.
\end{equation*}

\begin{remark}
In the above proof, the relation between the two gradings in Figure~\ref{fig:gradings_roots} is the key structure. 
From the $\sll_2$-representation point of view, it corresponds to the weight structures of the tensor product $\C^a\otimes \C^b$ where $\sll_2$ acts by coproduct or acts only on the second component. On the other hand, the identification of roots up to cancellation works as $a\leq b$ holds. This implies that a similar argument works in a more general setting.
Let $\g$ be a basic classical Lie superalgebra or a simple Lie algebra. Suppose that we have an embedding 
   $\mathfrak{a}_1\oplus \mathfrak{a}_2 \dots \oplus \mathfrak{a}_n \hookrightarrow \g$ 
so that the $\mathfrak{a}_i$'s are all simple Lie algebras and set 
$f=f_1+f_2+\dots+f_n$
where the $f_i$'s are regular nilpotent elements of $\mathfrak{a}_i$'s. 
Fix $\sll_2$-triples $\sll_2^{(i)}\subset \mathfrak{a}_i$ containing $f_i$.
Let 
$$\g\simeq (\mathfrak{a}_1\oplus \mathfrak{a}_2 \dots \oplus \mathfrak{a}_n)\oplus M$$  
as $\mathfrak{a}_1\oplus \mathfrak{a}_2 \dots \oplus \mathfrak{a}_n$-modules.
By using two $\sll_2$-triples, which are  $\sll_2=\sll_2^{(i)}$ and the diagonal embedding $\sll_2\hookrightarrow \bigoplus_{j>i}\sll_2^{(j)}$,
we decompose 
$$M\simeq \bigoplus_{A}\C^{d_A^1}\otimes \C^{d_A^2}$$
as an $\sll_2$-bimodule. Then, we expect that Theorem \ref{Conj at the level of q} holds when $d_A^1\leq d_A^2$ holds for all $A$ and $i$.
\end{remark}

\section{Universal objects for cosets of \tW-algebras in type \tA}\label{sec: Universal objects}  
In this section, we recall the universal $\W_\infty$-algebra $\W_\infty[c,\lambda]$ constructed in \cite{L} and discuss some expected construction and properties of generalizations based on our conjectures. 
{For simplicity we consider $k\in$ to be generic in this section.}

\subsection{\tWinf-algebra} \label{sec: W_infty-algebra}
Recall that the regular $\W$-algebras $\W^k(\sll_N,f_N)$ form a family of vertex algebras obtained as quotients of a universal object, called the $\W_\infty$-algebra $\W_\infty[c,\lambda]$ \cite{L}. The vertex algebra $\W_\infty[c,\lambda]$ is uniquely defined, free over the polynomial ring $\C[c,\lambda]$ and weakly generated by the fields $L, W_3$, completed into the strong generating set $\{W_2=L, W_3, W_4,\dots\}$, which gives a PBW base.
By definition, the weak generators satisfy the OPEs
\begin{align*}
    &\lbr{L}{L}=\frac{c}{2}\lm{3}+2L\lm{1} +L',\qquad \lbr{L}{W_3}=3 W_3\lm{1} +W_3',
\end{align*}
and the remaining strong generators are defined as 
\begin{align*}
    W_{n+1}=W_3{}_{(1)}W_n\qquad (n\geq 3).
\end{align*}
The complete algebraic structure on $\W_\infty[c,\lambda]$ is uniquely determined by the following properties:
\begin{itemize}
    \item (Normalization) $\lbr{W_{3}}{W_3}=\tfrac{c}{3}\lm{5}+(\text{lower order terms})$.
    \item (Involution) The assignment $W_n\mapsto (-1)^nW_n$ $(n\geq 2)$ extends to an automorphism of $\W_\infty[c,\lambda]$.
\end{itemize}
For example, we have the following OPEs (see \S \ref{Winfinity algebra} for more OPEs in lower conformal weights)
\begin{align*}
    &\lbr{W_3}{W_3}=\tfrac{c}{3}\lm{5}+2L\lm{3}+L' \lm{2}+W_4 \lm{1} +(\tfrac{1}{2}W_4'-\tfrac{1}{12}L'''),\\
    &\lbr{L}{W_5}=(185-80\lambda(c+2))W_3\lm{3}+(55-16\lambda(c+2))W_3'\lm{2}+5W_5\lm{1}+W_5'.
\end{align*}

Let us introduce the following specializations 
\begin{align}
    &c_{N,m}(k)=-\frac{(n\psi-N-1)(n\psi-\psi-N+1)(n\psi+\psi-N)}{(\psi-1)\psi}\label{eq:c_parameter}\\
    &\lambda_{N,m}(k)=-\frac{(\psi-1)\psi}{(n\psi-N-2)(n\psi-2\psi-N+2)(n\psi+2\psi-N)}\label{eq:l_parameter}
\end{align}
with $N=n+m$ and $\psi=k+N$. 
Then the corresponding $\W_\infty$-algebra $\W_\infty[c_{N,m}(k),\lambda_{N,m}(k)]$ has a maximal ideal $\mathscr{I}_{N,m}(k)$ whose lowest weight component occurs in weight $(m+1)(N+1)$ and has the form
\begin{align}\label{truncation}
    W_{(m+1)(N+1)}=P(L,W_3,\dots,W_{(m+1)(N+1)-1})
\end{align}
for some differential polynomial $P(L,W_3,\dots,W_{(m+1)(N+1)-1})$.
The regular $\W$-algebra $\W^k(\sll_N,f_N)$ is of strong generating type
$$\W^k(\sll_N,f_N)=\W(2,3,\dots,N)$$
and is obtained as the quotient 
\begin{align*}
    \W^k(\sll_N,f_N)\simeq \W_\infty[c_{N,0}(k),\lambda_{N,0}(k)]/\mathscr{I}_{N,0}(k)
\end{align*}
as a deformable family of vertex algebras depending on the level $k$.

More generally, let us consider the hook-type $\W$-algebra
$\W^k(\sll_{N},f_{1^m,n})$. Recall from \S\ref{sec:WslN} that it contains the affine vertex subalgebra 
\begin{align*}
   V^{k^\sharp}(\gl_m):=\pi^{J_{N,m}}\otimes V^{k^\sharp}(\sll_m)\hookrightarrow \W^k(\sll_{N},f_{1^m,n})
\end{align*}
with $k^\sharp=k+n-1$ and the Heisenberg field $J_{N,m}$ satisfies the relation
\begin{align}\label{Heisenberg norm}
    \lbr{J_{N,m}}{J_{N,m}}=\left(-m+\frac{n m}{N}(k+N)\right)\Lambda.
\end{align}
Then the (universal) affine coset subalgebra 
\begin{align*}
C^k(\sll_{N},f_{1^m,n}):=\Com{V^{k^\sharp}(\gl_m), \W^k(\sll_{N},f_{1^m,n})}
\end{align*} 
is of strong generating type 
\begin{align}\label{strong generating type of affine coset}
    C^k(\sll_{N},f_{1^m,n})=\W(2,3,\dots,(m+1)(N+1)-1)
\end{align}
and indeed obtained as the quotient 
\begin{align*}
    C^k(\sll_{N},f_{1^m,n})\simeq \W_\infty[c_{N,m}(k),\lambda_{N,m}(k)]/\mathscr{I}_{N,m}(k)
\end{align*}
as a deformable family of vertex algebras depending on the level $k$. 
We note that the strong generators of $C^k(\sll_{N},f_{1^m,n})$, which have conformal weights $2,3,\dots,(m+1)(N+1)-1$ in \eqref{strong generating type of affine coset}, do not form a PBW base in general contrary to the case $\W^k(\sll_N,f_N)$ mentioned previously. 
The existence of the strong generators comes from the classical orbifold theory appearing in the large level limit $k\rightarrow\infty$, thanks to \cite{CL5}.
As we will use it below, we explain it briefly following \cite{CL1}.

In general the $\W$-algebra $\W^k(\g,f)$ admits an integral form $\W^{\mathbf{k}}(\g,f)$, which is defined and free over $R=\C[\mathbf{k}]$ so that the specializations $\mathbf{k}=k$ with $k\in \C$ recover the original one:
\begin{align*}
    \W^{\mathbf{k}}(\g,f)/(\mathbf{k}-k)\simeq \W^k(\g,f).
\end{align*}
In the case of the affine vertex algebra $V^\mathbf{k}(\g)$, one of the meaningful large level limit $k\rightarrow\infty$ within the theory of vertex algebras is introduced in the following manner. 
We replace $R$ with $\widehat{R}=\C[\mathbf{k}^{\pm1/2}]$ and consider $V^\mathbf{k}_{\widehat{R}}(\g):=V^\mathbf{k}(\g)\otimes_{R}\widehat{R}$. 
Then we introduce the integral form $V^\mathbf{k}_{\widehat{R}_\infty}(\g)$ over $\widehat{R}_\infty=\C[\mathbf{k}^{-1/2}]$, which is strongly generated by $a_{\widehat{R}_\infty}:=\mathbf{k}^{-1/2}a$ for $a\in \g$. 
The large level limit $V^\infty(\g)=\lim_{k\rightarrow\infty}V^k(\g)$ is by definition
$$V^\infty(\g):=V^\mathbf{k}_{\widehat{R}_\infty}(\g)/(\mathbf{k}^{-1/2}).$$
Since 
\begin{align*}
    \lbr{a_{\widehat{R}_\infty}}{b_{\widehat{R}_\infty}}=(a,b)\Lambda+\mathbf{k}^{-1/2}[a,b]_{\widehat{R}_\infty}\qquad (a,b\in \g),
\end{align*}
$V^\infty(\g)$ is strongly generated by the images $a_{\infty}$ of $a_{\widehat{R}_\infty}$ ($a\in \g$) satisfying the OPEs
\begin{align*}
    \lbr{a_{\infty}}{b_{\infty}}=(a,b)\Lambda\qquad (a,b\in \g).
\end{align*}
Hence, $V^\infty(\g)$ is nothing but the Heisenberg vertex algebra associated with $\g$ and the normalized invariant bilinear form $(\text{-},\text{-})$. 
As for the $\W$-algebras $\W^k(\g,f)$ in general, we replace the BRST complex $C_f(V^k(\g))$ with an integral form $C_f(V^\mathbf{k}_{\widehat{R}}(\g))$ defined over $\widehat{R}$ equipped with a slightly rescaled differential, see \cite[\S 3.1]{CL1} for details.
We introduce $\W^\mathbf{k}_{\widehat{R}_\infty}(\g,f)$ as a vertex subalgebra over $\widehat{R}_\infty$ of the zero-th cohomology $H^0(C_f(V^\mathbf{k}_{\widehat{R}}(\g)))$ (denoted as $\W^k_{\sigma}(\g,f)$ with $\sigma^2=\epsilon$ in \emph{loc.cit}), which satisfies 
\begin{align*}
    \W^\mathbf{k}_{\widehat{R}_\infty}(\g,f)/(\mathbf{k}^{-1/2}-k^{-1/2})\simeq \W^k(\g,f),
\end{align*}
and then define the large level limit $k\rightarrow \infty$ to be 
\begin{align}\label{large level limit}
    \W^\infty(\g,f):=\W^\mathbf{k}_{\widehat{R}_\infty}(\g,f)/(\mathbf{k}^{-1/2}).
\end{align}
Then $\W^\infty(\g,f)$ is strongly generated by elements in $\g^f$ satisfying the OPEs
\begin{align*}
    \lbr{a}{b}=\delta_{p,q} B_p(a,b)\lm{2p+1}\qquad (a\in \g_{-p}^f,\ b\in \g^f_{-q})
\end{align*}
where 
\begin{align*}
    B_p(\text{-},\text{-})\colon \g^f_{-p}\times \g^f_{-p}\rightarrow \C,\quad (a,b)\mapsto (-1)^{2p}(\ad_f^{2p}a,b)
\end{align*}
by \cite[Thm. 3.5]{CL1}.

Now, consider the $\W$-algebra $\W^k(\sll_N,f_{1^m,n})$ which is of strong generating type 
\begin{align*}
    \W^k(\sll_N,f_{1^m,n})=\W(1^{m^2},[2]_{n-1},(\tfrac{n+1}{2})^{2m})
\end{align*}
by Proposition \ref{strong generating type}.
The weight $1$ fields generate the affine vertex algebra $V^{k^\sharp}(\gl_m)$, the fields of weights $2,\dots,n$ give the $(n-1)$ first lowest conformal weight strong generators, say $\{W_2,\dots W_n\}$, of $C^k(\sll_{N},f_{1^m,n})$, and the $2m$ fields of weights $(n+1)/2$ form bases of the Weyl modules induced from the natural representation $\C^m$ and its dual $\overline{\C}^m$
\begin{align}\label{weyl modules in module extesion}
    \C^m\subset \weyl_{\gl_m}^{k^\sharp}(\varpi_1),\quad \overline{\C}^m\subset \weyl_{\gl_m}^{k^\sharp}(\varpi_m),
\end{align}
as $V^{k^\sharp}(\gl_m)$-modules. By \eqref{large level limit}, we have 
\begin{equation*}
\W^\infty(\sll_N,f_{1^m,n})\simeq \left\{
\begin{array}{ll}
\mathcal{O}_{\text{ev}}(m^2,2) \otimes \big( \bigotimes_{i=2}^n \mathcal{O}_{\text{ev}}(1,2i)\big) \otimes \mathcal{S}_{\text{ev}}(m, n+1), & n \ \text{even},\\ 
\mathcal{O}_{\text{ev}}(m^2,2) \otimes \big( \bigotimes_{i=2}^n \mathcal{O}_{\text{ev}}(1,2i)\big) \otimes \mathcal{O}_{\text{ev}}(2m, n+1), & n \ \text{odd}.
\end{array} 
\right.
\end{equation*}
Here we have used the generalized free field algebras
\begin{itemize}
    \item $\mathcal{S}_{\mathrm{ev}}(n,\ell)$: strongly generated by the fields $a^i$, $b^i$ ($i=1,\dots,n$) of weight $\ell/2$ satisfying the OPEs
    \begin{align*}
        \lbr{a^i}{b^j}=\delta_{i,j}\lm{\ell-1},\quad \lbr{a^i}{a^j}=0,\quad \lbr{b^i}{b^j}=0,
    \end{align*}
    \item $\mathcal{O}_{\mathrm{ev}}(n,2\ell)$: strongly generated by the fields $a^i$ ($i=1,\dots,n$) of weight $\ell$ satisfying the OPEs
    \begin{align*}
        \lbr{a^i}{a^j}=\delta_{i,j}\lm{\ell-1}.
    \end{align*}
\end{itemize}
In the limit $k\rightarrow \infty$, the factor $\mathcal{O}_{\mathrm{ev}}(m^2,2)$ corresponds to $V^{k^\sharp}(\gl_m)$, the factor $\mathcal{S}_{\mathrm{ev}}(m,n+1)$ (resp. $\mathcal{O}_{\mathrm{ev}}(2m,n+1)$) to the differential polynomials generated by $\C^m$ and $\overline{\C}^m$, and the factor $\mathcal{O}_{\mathrm{ev}}(1,2i)$ to those generated by $W_i$ for $i=2,\dots,n$.
Then, by \cite[Lemma 4.2]{CL1}, the large level limit $k\rightarrow \infty$ of $C^k(\sll_{N},f_{1^m,n})$ coincides with the orbifold
\begin{align*}
    C^\infty(\sll_{N},f_{1^m,n})\simeq  \left\{
\begin{array}{ll}
    \mathcal{S}_{\mathrm{ev}}(m,n+1)^{\mathrm{GL}_m}\otimes \bigotimes_{i=2}^n \mathcal{O}_{\mathrm{ev}}(1,2i),& n \ \text{even},\\
    \mathcal{O}_{\mathrm{ev}}(2m,n+1)^{\mathrm{GL}_m}\otimes \bigotimes_{i=2}^n \mathcal{O}_{\mathrm{ev}}(1,2i), & n \ \text{odd}. 
    \end{array} 
\right.
\end{align*}

The strong generating type of $\mathcal{S}_{\mathrm{ev}}(m,n+1)^{\mathrm{GL}_m}$ is obtained by using methods from the invariant theory first developed in \cite{L2,L3}.
The generators of $\mathcal{S}_{\mathrm{ev}}(m,n+1)^{\mathrm{GL}_m}$ (resp. $\mathcal{O}_{\mathrm{ev}}(2m,n+1)^{\mathrm{GL}_m}$) transform as $\C^m \oplus \overline{\C}^m$ under $\mathrm{GL}_m$ and so do the derivatives of generators.
By Weyl's first fundamental theorem of invariant theory, the $\mathrm{GL}_m$-invariants are generated by quadratics obtained by the pairing 
\begin{align*}
    \omega_{p,q}=\sum_{i=1}^m \partial^p a^i\ \partial^q b^i,\qquad (p,q \geq 0)
\end{align*}
which one may restrict to $p=0$ by using the derivative $\partial$.
This implies that $\mathcal{S}_{\mathrm{ev}}(m,n+1)^{\mathrm{GL}_m}$ is at most of strong generating type $\W([n+1]_{\infty})$. 
Note that this is \emph{not} a minimal set of strong generators as there are relations between them coming from Weyl's second fundamental theorem. 
These relations are generated by $(n+1) \times (n+1)$ determinants in the generators. 
Thus, it has a strong generating type \cite[Thm.\ 4.6, 4.9]{CL1}
$$\mathcal{S}_{\mathrm{ev}}(m,n+1)^{\mathrm{GL}_m}=\W(n+1,\dots,(m+1)(N+1)-1).$$
Similarly, $\mathcal{O}_{\mathrm{ev}}(2m,n+1)^{\mathrm{GL}_m}$ has the same strong generating type by \cite[Thm.\ 4.9]{CL1}. Hence, we have 
\begin{align*}
    C^\infty(\sll_{N},f_{1^m,n})=\W(2,3,\dots,(m+1)(N+1)-1)
\end{align*}
and so is $ C^k(\sll_{N},f_{1^m,n})$ for generic $k\in \C$ as desired.

 \subsection{Height-two partitions} 
Let us consider the family of $\W$-algebras $\W^k(\sll_{N},f_{n+r,m+r})$ for $r\geq 0$ with $m> n\geq 2$ fixed. 
By Proposition \ref{strong generating type}, we have
\begin{align*}
    \W^k(\sll_{N},f_{n+r,m+r})=\W(1,[2]_{n+r-1},[2]_{m+r-1},[\sigma]_{n+r},[\sigma]_{n+r})
\end{align*}
with $\sigma=\tfrac{m-n}{2}+1$. 
Assuming Conjecture \ref{conj:decomposition}, the strong generators of conformal weights $1,[2]_{n+r-1},[2]_{m+r-1}$ give an embedding of mutually-commuting vertex subalgebras
\begin{align}\label{subalgebra of length two W}
 \pi^{J_{N,n+r}}\otimes \mathscr{W}_1^{[r]} \otimes \mathscr{W}_2^{[r]} \hookrightarrow \W^k(\sll_{N},f_{n+r,m+r})
\end{align}
 where the Heisenberg vertex algebra $\pi^{J_{N,n+r}}$ is defined in \eqref{Heisenberg norm} and
\begin{align*}
    \mathscr{W}_1^{[r]}=C^{k+m+r-1}(\sll_{n+r},f_{n+r}),\quad \mathscr{W}_2^{[r]}=C^k(\sll_{n+m+2r},f_{1^{n+r},m+r}).
\end{align*}
The remaining fields of conformal weights $[\sigma]_{n+r},[\sigma]_{n+r}$ (including the multiplicities) describe modules extending the embedding \eqref{subalgebra of length two W}, which can be taken as the bases of $\C^m$ and $\overline{\C}^m$ in \eqref{weyl modules in module extesion} under the BRST reduction $H_{f_m}$ after some quantum corrections, thanks to Conjecture \ref{conj:successive_reductions}.
Therefore, we expect that this family of $\W$-algebras is obtained as quotients of a universal vertex algebra $\W_{\infty}^{[2;\sigma]}[c,\lambda]$ such that 
\begin{itemize}
    \item it is a vertex algebra freely defined over the polynomial ring $\C[c,\lambda]$ and of strong generating type 
    \begin{align*}
    \W_{\infty}^{[2;\sigma]}[c,\lambda]=\W(1,[2]_{\infty},[2]_{\infty},[\sigma]_{\infty},[\sigma]_{\infty}),
\end{align*}
   \item it is an extension of two copies of $\W_\infty$-algebras \cite{L}
   \begin{align}\label{relation of parameters}
    \pi\otimes \W_\infty[c_1,\lambda_1]\otimes \W_\infty[c_2,\lambda_2]\hookrightarrow \W_{\infty}^{[2;\sigma]}[c,\lambda].
\end{align}
\end{itemize} 
\begin{remark}
The enveloping algebra of $\W_\infty[c,\lambda]$ is isomorphic to the quantum group called \emph{the affine Yangian associated with $\gl_1$} after appropriate completions for both of them. Vertex algebra extensions \eqref{relation of parameters} is discussed in physics in this language \cite{Li20, LP19}.
The relation of the parameters in \eqref{relation of parameters} is concisely expressed by using the Yangian parameters $(\mu_1,\mu_2,\mu_3)$ uniquely determined by 
\begin{align*}
    \frac{1}{\lambda_1}+\frac{1}{\lambda_2}+\frac{1}{\lambda_3}=0,\quad \frac{N}{\mu_2}+\frac{m}{\mu_3}=1,\quad -\frac{\mu_1}{\mu_2}=\Psi
\end{align*}
in terms of $(c,\lambda)$ by solving $c=c_{N,m}(k)$, $\lambda=\lambda_{N,m}(k)$, and $\Psi=k+N$ see \eqref{eq:c_parameter}-\eqref{eq:l_parameter}.
Let $(\mu_{i;1},\mu_{i;2},\mu_{i;3})$ be the Yangian parameter for $\W_\infty[c_i,\lambda_i]$ in \eqref{relation of parameters} for $i=1,2$. 
These six parameters satisfy 
\renewcommand{\arraystretch}{1.3}
\begin{align*}
\begin{array}{ll}
   \displaystyle{\frac{1}{\mu_{1;1}}+\frac{1}{\mu_{1;2}}+\frac{1}{\mu_{1;3}}=0},  & \displaystyle{\frac{1}{\mu_{2;1}}+\frac{1}{\mu_{2;2}}+\frac{1}{\mu_{2;3}}=0}, \\
   &\\
    \displaystyle{\frac{\mu_{1;1}}{\mu_{1;2}}-\frac{\mu_{2;1}}{\mu_{2;2}}=1}, & \displaystyle{\mu_{1;1}-\frac{\mu_{1;1}}{\mu_{1;2}}\hat{\sigma}=-\mu_{2;2}},
\end{array}
\end{align*}
\renewcommand{\arraystretch}{1}
with $\hat{\sigma}=2(\sigma-1)$.
Therefore, one has two remaining parameters, corresponding to $(c,\lambda)$ for $\W_{\infty}^{[2;\sigma]}[c,\lambda]$.  
\end{remark}

As in the case of $\W_\infty[c,\lambda]$ in \S \ref{sec: W_infty-algebra}, we expect that after specializing $(c,\lambda)$, quotients of $\W_{\infty}^{[2;\sigma]}[c,\lambda]$ recover the following affine cosets of $\W$-superalgebras,
\begin{align*}
    \Com{V^{k^\sharp}(\gl_a),\W^k(\sll_{N}, f_{1^a,n,m})},\quad \Com{V^{k^\sharp}(\gl_a),\W^k(\sll_{n+m|a}, f_{n,m|1^a})}
\end{align*}
with $a>n+m$.
As a non-trivial check, we derive the (minimal) strong generating type of the first affine coset. 
We can derive the second in a similar manner, and thus, we omit it.
By Proposition \ref{strong generating type}, we have
$$\W^k(\sll_{n+m+a}, f_{1^a,n,m})=\W(1^{a^2+1},[2]_{n-1},[2]_{m-1},[\sigma]_{n},[\sigma]_{n},(\tfrac{n+1}{2})^{2a},(\tfrac{m+1}{2})^{2a})$$
with $\sigma=\tfrac{m-n}{2}+1$. The strong generators of conformal weights $1^{a^2+1},[2]_{n-1},[2]_{m-1}$ (including the multiplicity) give a conformal embedding
\begin{align*}
    \pi^{J_N,a+n}\otimes V^{\ell}(\gl_a)\otimes C^{k^\sharp}(\sll_{a+n},f_{1^a,n})\otimes C^{k}(\sll_{N},f_{1^{a+n},m})\hookrightarrow \W^k(\sll_{n+m+a}, f_{1^a,n,m})
\end{align*}
with $\ell=k+m+n-2$ and $k^\sharp=k+m-1$. 
On the other hand, the remaining strong generators give module extension of this embedding:
those of conformal weights $[\sigma]_{n},[\sigma]_{n}$ commute with the affine subalgebra, and those of conformal weights $(\tfrac{n+1}{2})^{2a}$ and $(\tfrac{m+1}{2})^{2a}$ transform as the bases of 
$\C_1\otimes \C^m\oplus \C_{-1}\otimes \overline{\C}^m$ with respect to the $(\gl_1, \gl_m)$-action.
Then it follows from \cite[Thm. 3.5]{CL1} that 
\begin{align*}
\W^\infty(\sll_{n+m+a},f_{1^a,n,m})
\simeq \cO_{\text{ev}}(a^2,2) \otimes \mathcal{A}_1\otimes \mathcal{A}_2 \otimes  \mathcal{B}_1\otimes \mathcal{B}_2
\end{align*}
with 
\begin{align*}
 & \mathcal{A}_1=\bigotimes\limits_{i=1}^{n} \cO_{\text{ev}}(1,2i)\otimes \bigotimes_{i=2}^{m} \cO_{\text{ev}}(1,2i),\quad 
 \mathcal{A}_2=\begin{cases}
      \bigg(\bigotimes\limits_{0\leq i<m} \cO_{\text{ev}}(2,2(\sigma+i))\bigg), &(n\equiv m),\\
      \bigg(\bigotimes\limits_{0\leq i<m} \mathcal{S}_{\text{ev}}(1,2(\sigma+i))\bigg), &(n\not\equiv m),
  \end{cases}
\end{align*}
and 
\begin{align*}
    \mathcal{B}_1=\begin{cases}
   \mathcal{S}_{\text{ev}}(a, m+1), & (m\equiv0),\\
   \cO_{\text{ev}}(2a, m+1), & (m\equiv1),
  \end{cases}\quad 
  \mathcal{B}_2=\begin{cases}
   \mathcal{S}_{\text{ev}}(a, n+1), & (n\equiv0),\\
   \cO_{\text{ev}}(2a, n+1), & (n\equiv1),
  \end{cases}
\end{align*}
depending on $n,m$ modulo $2$.
Therefore, the affine coset
$$C^k(\sll_N, f_{1^a,n,m}):=\Com{V^{\ell}(\gl_a),\W^k(\sll_{N}, f_{1^a,n,m})}$$
admits the large level limit 
$$C^\infty(\sll_N, f_{1^a,n,m})\simeq \mathcal{A}_1\otimes \mathcal{A}_2 \otimes  (\mathcal{B}_1 \otimes \mathcal{B}_2)^{\mathrm{GL}_a}.$$
The strong generating type of $(\mathcal{B}_1 \otimes \mathcal{B}_2)^{\mathrm{GL}_a}$ is described using the invariant theory and an argument similar to the one in the previous section.
It is at most
$$\W([m+1]_{\infty}, [n+1]_{\infty}, [\sigma+n]_{\infty}, [\sigma+n]_{\infty}).$$ 
Here the generators of type $[m+1]_{\infty}$ and $[n+1]_{\infty}$ come from $\mathcal{B}_1^{\mathrm{GL}_a}$ and $\mathcal{B}_2^{\mathrm{GL}_a}$ respectively. 
The remaining generators (of type $[\sigma+n]_{\infty}$) come from the crossing pairing between the generators of $\mathcal{B}_1$ and $\mathcal{B}_2$.

By Weyl's second fundamental theorem, $[m+1]_{\infty} $ and $[n+1]_{\infty}$ truncate to $[m+1]_{R_m} $ and $[n+1]_{R_n}$ with $R_t=(a+t+1)(a+1)-(t+1)$ as before.
Moreover by using the relations appearing in this truncation, one can then construct relations for the remaining generators of $(\mathcal{B}_1 \otimes \mathcal{B}_2)^{\mathrm{GL}_a}$ in the blocks of type $[\sigma+n]_{\infty}$ starting at weight $(a+m+1)(a+1)$. 
It gives the following description, which is not necessarily the minimal strong generating set,
\begin{align*}
    C^\infty(\sll_N, f_{1^a,n,m})=\W(1,[2]_{r-1},[2]_{s-1},[\sigma]_{t-1},[\sigma]_{t-1}),
\end{align*}
for some $r,s,t \leq (a+m+1)(a+1)$ which are bounded below by a quadratic function of $a$. Therefore $C^k(\sll_N, f_{1^a,n,m})$ has the same strong generating type at generic levels. This can be seen as a truncation of $\W_{\infty}^{[2;\sigma]}[c,\lambda]$ as desired.

\subsection{More partitions}
We generalize the previous considerations to additional families of $\W$-algebras.
\begin{example}The $\W$-algebra $\W^k(\gl_{nr},f_{n,\dots,n})$ is of strong generating type
$$\W^k(\gl_{nr},f_{n,\dots,n})=\W(1^{n^2},2^{n^2},\dots r^{n^2}),$$
the conformal weight-one fields generate the affine vertex subalgebra $V^{k^\sharp}(\gl_n)$, and the higher conformal weight fields form the adjoint representations. 
By taking the limit $r\rightarrow \infty $, we expect that there exists a two-parameter vertex algebra of strong generating type 
$$\W(1^{n^2},2^{n^2}, 3^{n^2},\dots)$$
whose quotients recover the affine cosets 
\begin{align*}
\W^k(\gl_{nr},f_{n^r}),\quad \Com{V^{k^\sharp}(\gl_a), \W^k(\sll_{nr+a}, f_{1^a,n^r})},\quad \Com{V^{k^\sharp}(\gl_a), \W^k(\sll_{nr|a}, f_{n^r|1^a})}.
\end{align*}
The existence of such a vertex algebra is partially checked in \cite{EP} as a matrix extended $\W_{1+\infty}$-algebra, which conjecturally corresponds to the affine Yangian for $\gl_n$, see also \cite{KU}. 
\end{example}

\begin{example} The $\W$-algebra $\W^k(\gl_{N+nr},f_{\lambda_1+r,\dots,\lambda_n+r})$ with $\lambda_1>1$ is of strong generating type
$$\W^k(\gl_{N+nr},f_{\lambda_1+r,\dots,\lambda_n+r})=\W\left(\bigsqcup_{i,j=1}^n\left[\tfrac{|\lambda_j-\lambda_i|}{2}+1\right]_{\mathrm{min}\{\lambda_i,\lambda_j\}+r}\right).$$
Here again, by taking the limit $r\rightarrow \infty $, we expect that there exists a two-parameter vertex algebra of strong generating type 
$$\W\left(\bigsqcup_{i,j=1}^n\left[\tfrac{|\lambda_j-\lambda_i|}{2}+1\right]_{\infty}\right)$$
whose quotients recover the affine cosets 
\begin{align*}
\W^k(\gl_{N+nr},f_{\lambda_1+r,\dots,\lambda_n+r}),\quad \Com{V^{k^\sharp}(\gl_a), \W^k(\gl_{N+a}, f_{1^a,\lambda})},\quad  \Com{V^{k^\sharp}(\gl_a), \W^k(\gl_{N|a},  f_{\lambda|1^a})}
\end{align*}
where $\lambda=(\lambda_1,\dots,\lambda_n)$.
Such a vertex algebra is expected to correspond to a shifted generalization of the affine Yangian for $\gl_n$.
\end{example}

One can also adapt the previous discussion for the $\W$-superalgebras $\W^k(\gl_{N|M,f_{\lambda|\mu}})$ in type $A$. 
In the super-setting, the nilpotent orbits are taken to be even by definition and parameterized by two partitions
$$f=f_{\lambda|\mu}=f_{\lambda|0}+f_{0|\mu}\in \gl_N\oplus \gl_M\subset \gl_{N|M},\quad \lambda\vdash N,\ \mu\vdash M$$
where $f_{\lambda|0}, f_{0|\mu}$ are taken as in \eqref{nilpotent element}.
Then a similar pyramid as in Figure \ref{fig:sym_grading_gl9} defines a good grading and thus $\W^k(\gl_{N|M,f_{\lambda|\mu}})$ has a strong generating type as in Proposition \ref{strong generating type} decorated by the parity information, namely
\begin{align}\label{conformal wt for super case}
    \W^k(\gl_{N|M},f_{\lambda|\mu})=\W\left(I^+_{\lambda,\mu};I^-_{\lambda,\mu} \right)
\end{align}
where the conformal weights of even generators are
\begin{align*}
   I^+_{\lambda,\mu}= \bigsqcup_{i,j}\left[\tfrac{|\lambda_j-\lambda_i|}{2}+1\right]_{\mathrm{min}\{\lambda_i,\lambda_j\}} \sqcup \bigsqcup_{i,j}\left[\tfrac{|\mu_j-\mu_i|}{2}+1\right]_{\mathrm{min}\{\mu_i,\mu_j\}}
\end{align*}
and the conformal weights of odd generators are
\begin{align*}
 I^-_{\lambda,\mu}= \bigsqcup_{i,j}\left[\tfrac{|\mu_j-\lambda_i|}{2}+1\right]_{\mathrm{min}\{\lambda_i,\mu_j\}}.
\end{align*}

\begin{example}We consider the $\W$-superalgebras $\W^k(\gl_{n+r|m+r},f_{n+r|m+r})$ with $n\geq m$. They are of strong generating type
$$\W^k(\gl_{n+r|m+r},f_{n+r|m+r})=\W([1]_{n+r}\sqcup [1]_{m+r};\ [\tfrac{n+m}{2}+1]_{m+r}\sqcup [\tfrac{n+m}{2}+1]_{m+r})$$
By taking the limit $r\rightarrow \infty $, we expect that there exists a two-parameter vertex algebra of strong generating type 
$$\W([1]_{\infty}\sqcup [1]_{\infty};\ [\tfrac{n+m}{2}+1]_{\infty}\sqcup [\tfrac{n+m}{2}+1]_{\infty})$$
whose quotients recover the affine cosets 
 \begin{equation*}
 \begin{gathered}
    \W^k(\gl_{n+r|m+r},f_{n+r|m+r}),\\
    \Com{V^{k^\sharp}(\gl_a), \W^k(\gl_{n+a| m}, f_{1^a,n| m})},\quad
    \Com{V^{k^\sharp}(\gl_a), \W^k(\gl_{n| m + a}, f_{n|1^a, m})}.
 \end{gathered}
 \end{equation*}
\end{example}

\begin{example} The $\W$-superalgebra $\W^k(\gl_{N+nr|M+mr},f_{\lambda_1+r,\dots,\lambda_n+r|\mu_1+r,\dots,\mu_m+r })$, with $\lambda_1,\mu_1>1$, is of strong generating type 
$\W(I^{+[r]}_{\lambda,\mu};I^{-[r]}_{\lambda,\mu})$ with $I^{\pm[r]}_{\lambda,\mu}$ obtained from $I^{\pm}_{\lambda,\mu}$ replacing
$$\mathrm{min}\{\lambda_i,\lambda_j\}\mapsto \mathrm{min}\{\lambda_i,\lambda_j\}+r,\quad
\mathrm{min}\{\mu_i,\mu_j\}\mapsto \mathrm{min}\{\mu_i,\mu_j\}+r,\quad
\mathrm{min}\{\lambda_i,\mu_j\}\mapsto \mathrm{min}\{\lambda_i,\mu_j\}+r$$
in \eqref{conformal wt for super case}.
By taking the limit $r\rightarrow \infty $, we expect that there exists a two-parameter vertex algebra of strong generating type $\W(I^{+[\infty]}_{\lambda,\mu};I^{-[\infty]}_{\lambda,\mu} )$
whose quotients recover the affine cosets 
 \begin{equation*}
 \begin{gathered}
     \Com{V^{k^\sharp}(\gl_a), \W^k(\gl_{N+s|M}, f_{1^a,\lambda|\mu})},\quad 
     \Com{V^{k^\sharp}(\gl_s), \W^k(\gl_{N|M+a}, f_{\lambda|1^a,\mu})},\\
     \Com{V^{k^\sharp}(\gl_{r|a}), \W^k(\gl_{N+r|M+a}, f_{1^r,\lambda|1^a,\mu})}.
 \end{gathered}
 \end{equation*}
\end{example}

\subsection{Conjectures on duality}\label{Conjectures on duality}
The obvious generalization of Conjecture~\ref{conj:successive_reductions} in the super-setting identifies the webs of $\W$-algebras in \cite{PR} with boundary conditions associated with resolved conifold diagrams and gives a conjecture on Feigin--Frenkel type duality for $\W$-superalgebras. 
The set-up of such boundary conditions under flip relations is described in Figure \ref{fig: resolved conifold diagrams}. 
The corresponding webs of $\W$-algebras are 
\begin{align*}
    H_{f_{0|n}}H_{f_{n+r|0}}(V^k(\gl_{n+r|n})),\quad \Com{V^{\ell+r}(\gl_r),H_{f_{1^n,r}}(V^k(\gl_{n+r}))\otimes \mathcal{F}},
\end{align*}
respectively. Here $\mathcal{F}$ is the free field algebras
\begin{align}\label{eq:freefield}
    \mathcal{F}=\mathrm{SB}^{0|n}\ (r> 0),\quad \mathcal{F}=\mathrm{SB}^{n|n}\ (r=0)
\end{align}
where $\mathrm{SB}^{p|q}=\beta\gamma^{\otimes p}\otimes bc^{\otimes q}$ is the symplectic bosons which correspond to the three-dimensional hypermultiplets on the horizontal lines.
\setlength{\unitlength}{1mm}
\begin{figure}[htbp]
\centering
\begin{picture}(80,22)(0,0)
\put(20,13){\line(1,0){8}}
\put(20,13){\line(0,1){8}}
\put(20,13){\line(-1,-1){5}}
\put(15,8){\line(-1,0){8}}
\put(15,8){\line(0,-1){8}}
\put(23,16){\footnotesize$n+r$}
\put(21,7){\footnotesize$0$}
\put(13,13){\footnotesize$n$}
\put(10,4){\footnotesize$0$}
\put(60,13){\line(1,0){8}}
\put(60,13){\line(0,1){8}}
\put(60,13){\line(-1,-1){5}}
\put(55,8){\line(-1,0){8}}
\put(55,8){\line(0,-1){8}}
\put(63,16){\footnotesize$n+r$}
\put(61,7){\footnotesize$n$}
\put(53,13){\footnotesize$0$}
\put(50,4){\footnotesize$0$}
\end{picture}
\caption{Resolved conifold diagrams under horizontal flip relation}
\label{fig: resolved conifold diagrams}
\end{figure}
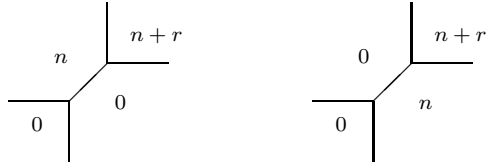
Hence, the reduction by stages for $\W$-superalgebras implies the following conjecture of $\W$-superalgebras.

\begin{conjecture}\label{new type of duality} Let $k,\ell$ be generic levels satisfying $(k+r)(\ell + n+r) = 1$.
  \begin{enumerate}[wide, labelindent=0pt]
    \item For $r\geq 1$, there exists an isomorphism of vertex superalgebras
    $$\W^k(\sll_{n+r|n}, f_{n+r|n}) \simeq \Com{V^{\ell + r}(\gl_n), \W^{\ell}(\sll_{n+r}, f_{1^n,r}) \otimes \mathrm{SB}^{0|n}}.$$
    \item There exists an isomorphism of vertex superalgebras
  $$\W^k(\sll_{n|n}, f_{n|n}) \simeq \Com{V^{\ell }(\gl_n), V^{\ell}(\sll_{n})\otimes  \mathrm{SB}^{n|n}}.$$ 
\end{enumerate} 
\end{conjecture}
\begin{remark} The conjecture implies the strong rationality of the simple quotient 
$\W_k(\sll_{n+r|n},f_{n+r|n})$ at levels 
\begin{align*}
    k=-r+\frac{r}{p},\quad (p\geq n+r\ \mathrm{and}\ r>n).
\end{align*}
\end{remark}

The second part of Conjecture \ref{new type of duality} has been studied in \cite{CL4} and, in particular, has been proven for $n=2$.
The first part unifies the following three cases known in the literature:
\begin{itemize}
\item the case $n=0$, $r\geq2$ can be read as the Feigin--Frenkel duality \cite{FF91}
$$\W^k(\sll_r, f_r) \simeq \W^{\ell}(\sll_{r},f_r).$$
\item the case $n=1$, $r\geq1$ is the Kazama--Suzuki duality \cite{CGN}
$$\W^k(\sll_{r+1|1}, f_{r+1|1}) \simeq \Com{\pi, \W^{\ell}(\sll_{r+1}, f_{1,r}) \otimes V_\Z}$$
\item the case $n\geq1$, $r=1$ is a well-known conjecture by Ito \cite{Ito} 
$$\W^{k}(\sll_{n+1|n}) \cong \Com{V^{\ell+1}(\gl_n), V^{\ell}(\sll_{n+1}) \otimes V_{\Z^n}},$$
which is established for $n=1$ and $n=2$ \cite{GL}.
\end{itemize}

The coincidence of vertex subalgebras corresponding to the vertex algebras at the (two) corners in Figure \ref{fig: resolved conifold diagrams} supports Conjecture \ref{new type of duality}:
\begin{itemize}
\item For $r\geq1$: 
\begin{align*}
&\W^k(\sll_{n+r|n},f_{n+r|n})\supset \Com{V^{k^\sharp}(\gl_n),\W^k(\sll_{n+r|n},f_{n+r|1^n})}\otimes \W^{k^\sharp}(\gl_n,f_n),\\
&\Com{V^{\ell + r}(\gl_n), \W^{\ell}(\sll_{n+r}, f_{1^n,r}) \otimes \mathrm{SB}^{0|n}} \\
&\hspace{1cm}\supset \Com{V^{\ell^\sharp}(\gl_n),\W^\ell(\sll_{n+r},f_{1^n,r})}\otimes \Com{V^{\ell + r}(\gl_n),V^{\ell^\sharp}(\gl_n)\otimes \mathrm{SB}^{0|n}}\\
&\hspace{1cm}\supset \Com{V^{\ell^\sharp}(\gl_n),\W^\ell(\sll_{n+r},f_{1^n,r})}\otimes \W^{s}(\gl_n,f_{n}),
\end{align*}
\item For $r=0$:
\begin{align*}
&\W^k(\sll_{n|n},f_{n|n})\supset \Com{V^{k^\sharp}(\gl_n),\W^k(\sll_{n|n},f_{n+r|1^n})}\otimes \W^{k^\sharp}(\gl_n,f_n),\\
&\Com{V^{\ell}(\gl_n), V^{\ell}(\sll_{n}) \otimes \mathrm{SB}^{n|n}} \\
&\hspace{1cm}\supset \Com{V^{\ell^\sharp}(\gl_n),V^{\ell}(\sll_{n})\otimes \mathrm{SB}^{n|0}}\otimes \Com{V^{\ell}(\gl_n),V^{\ell^\sharp}(\gl_n)\otimes \mathrm{SB}^{0|n}}\\
&\hspace{1cm}\supset \Com{V^{\ell^\sharp}(\gl_n),V^{\ell}(\sll_{n})\otimes \mathrm{SB}^{n|0}} \otimes   \W^{s}(\gl_n,f_{n}),
\end{align*}
\end{itemize}
with $k^\sharp=-(k+n+r)+1$, $\ell^\sharp=\ell+r-1$, $s=-n+(\ell^\sharp+n)/(\ell^\sharp+n+1)$.
The affine coset subalgebras on the right-hand side are isomorphic \cite{CL1} whereas the regular $\W$-algebras are identical since $k^\sharp=s$ as desired.

Another set-up of resolved conifold diagrams under flip relations is given in Figure \ref{fig: resolved conifold diagrams 2}. 
\setlength{\unitlength}{1mm}
\begin{figure}[htbp]
\centering
\begin{picture}(80,22)(0,0)
\put(20,13){\line(1,0){8}}
\put(20,13){\line(0,1){8}}
\put(20,13){\line(-1,-1){5}}
\put(15,8){\line(-1,0){8}}
\put(15,8){\line(0,-1){8}}
\put(23,16){\footnotesize$n$}
\put(21,7){\footnotesize$0$}
\put(12,13){\footnotesize$N$}
\put(10,4){\footnotesize$m$}
\put(60,13){\line(1,0){8}}
\put(60,13){\line(0,1){8}}
\put(60,13){\line(-1,-1){5}}
\put(55,8){\line(-1,0){8}}
\put(55,8){\line(0,-1){8}}
\put(63,16){\footnotesize$m$}
\put(61,7){\footnotesize$0$}
\put(52,13){\footnotesize$N$}
\put(50,4){\footnotesize$n$}
\end{picture}
\caption{Resolved conifold diagrams under vertical flip relation}
\label{fig: resolved conifold diagrams 2}
\end{figure}
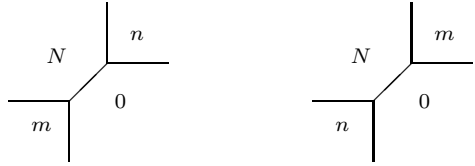
By assuming the condition $N>n,m\geq2$, the corresponding webs of $\W$-algebras are 
\begin{align*}
    \Com{V^{k^\sharp}(\gl_{m}),H_{f_{0|1^m,N-m}}H_{f_{n|0}}(V^k(\sll_{n|N})},\quad \Com{V^{\ell^\sharp}(\gl_{n}),H_{f_{0|1^n,N-n}}H_{f_{m|0}}(V^k(\sll_{m|N})}
\end{align*}
with 
$$k^\sharp=-k-N+(n+m),\quad \ell^\sharp=-\ell-N+(n+m).$$
Therefore, it gives rise to the following conjecture of $\W$-superalgebras.
\begin{conjecture}\label{new type of duality 2} For $N>n,m\geq2$ and the generic levels $k,\ell$ satisfying $(k+n-N)(\ell+m-N) = 1$, there exists an isomorphism of vertex superalgebras
 $$\Com{V^{k^\sharp}(\gl_{m}), \W^k(\sll_{n|N}, f_{n|1^m,N-m})} \simeq \Com{V^{\ell^\sharp}(\gl_{n}), \W^\ell(\sll_{m|N}, f_{m|1^n,N-n})}.$$
\end{conjecture}
   As a non-trivial check of Conjecture \ref{new type of duality 2}, one can check again the coincidence of vertex subalgebras corresponding to the corners:
\begin{align*}
&\Com{V^{k^\sharp}(\gl_{m}), \W^k(\sll_{n|N}, f_{n|1^m,N-m})}\\
&\hspace{1cm}\supset \Com{V^{\overline{k}}(\gl_{N}), \W^k(\sll_{n|N}, f_{n|1^N})}\otimes \Com{V^{k^\sharp}(\gl_{m}), \W^{\overline{k}}(\sll_{N}, f_{1^m,N-m})},\\
&\Com{V^{\ell^\sharp}(\gl_{n}), \W^k(\sll_{m|N}, f_{m|1^n,N-n})}\\
&\hspace{1cm}\supset \Com{V^{\overline{\ell}}(\gl_{N}), \W^\ell(\sll_{m|N}, f_{m|1^N})}\otimes \Com{V^{\ell^\sharp}(\gl_{n}), \W^{\overline{\ell}}(\sll_{N}, f_{1^n,N-n})}
\end{align*}
with $\overline{k}=-(k+n)+1$, $\overline{\ell}=-(\ell+m)+1$.
Then, the affine coset subalgebras on the right-hand side are isomorphic as justified by the relations on their levels \cite{CL1}
\begin{align*}
    (\overline{\ell}+N)+\frac{1}{k+n-N}=1,\quad (\overline{k}+N)+\frac{1}{\ell+m-N}=1.
\end{align*}

\section{Non-trivial examples in low ranks}\label{sec: Non-trivial examples in low ranks}
This section is devoted to examples verifying Conjectures~\ref{conj:successive_reductions}, \ref{conj:decomposition} and \ref{conj:embedding}. The main goal is to relate hook-type and non-hook-type $\W$-algebras. Hence we consider $\W$-algebras corresponding to the partition $(2,2)$ in $\sll_4$ -- often referred has the rectangular $\W$-algebra in $\sll_4$ -- and partitions $(2,3)$ and $(1,2^2)$ in $\sll_5$. The proofs are based on the descriptions of the $\W$-algebras that can be found in \S\ref{appx1}.
\subsection{The case \tWtype{4}{2,2} }
In this section, we relate the rectangular $\W$-algebra $\W^k(\sll_4,f_{2,2})$ and the minimal $\W$-algebra $\W^k(\sll_4,f_{1^2,2})$. 
We prove the following theorem.
\begin{theorem}\label{thm: the case 2,2} \hspace{0mm}
\begin{enumerate}[wide, labelindent=0pt]
    \item We have an isomorphism of vertex algebras
    \begin{align*}
        H_{f_{2}}(\W^k(\sll_4,f_{1^2,2}))\simeq \W^k(\sll_4,f_{2,2}).
    \end{align*}
    \item We have an isomorphism of vertex algebras
        \begin{align*}
        \W^k(\sll_4,f_{1^2,2})\simeq (\W^k(\sll_4,f_{2,2})\otimes \Pi[0])^{\SL_2},
    \end{align*}
    where $M^{\SL_2}$ denotes the maximal $\SL_2$-submodule of the $\sll_2$-module $M$.
    \item Moreover, if $k$ is generic (e.g. $k\notin \Q$) then we have a conformal embedding
    \begin{align*}
       C^k(\sll_4,f_{1^2,2})\otimes \W^{k+1}(\sll_2)\otimes \pi^J\hookrightarrow \W^k(\sll_4,f_{2,2}).
    \end{align*}
\end{enumerate} 
\end{theorem}
\noindent

Before we prove this theorem, we recall the structure of the main actor, that is $\W^k(\sll_4,f_{2,2})$.
By \S \ref{2,2}, the vertex algebra $\W^k(\sll_4,f_{2,2})$ at levels $k\neq-2,-3,-4$ has three mutually-commuting Virasoro vectors, namely $L_1, L_2$ and $L_J:=\frac{1}{2(k+2)}J^2$, whose sum $L_{\tot}=L_1+L_2+L_J$ is the conformal vector defined for non-critical levels, see \cite{KW04,KW22}. The central charges of $L_1, L_2$ agree with $c_{4,2}(k)$ and $c_{2,0}(k+1)$ given by the formula \eqref{eq:c_parameter}.
For $k\neq-1,-2,-3,-4$, let us introduce the following element of conformal weight-three in $\W^k(\sll_4,f_{2,2})$:
\begin{align*}
    \Omega_3=&-\frac{1}{3(2+k)}J^3 - \frac{(2+k)}{(1+k)}L_1J + \frac{(4+k)}{(1+k)}L_2J + \frac{(3+k)}{(1+k)}(JJ'-Jv_1^+v_1^-)\\
    &\quad+ \frac{(2+k)}{2(1+k)}(v_1^+v_2^- +v_1^-v_2^+) + \frac{(16+11k+2k^2)}{4(1+k)}({v_1^+}'v_1^- -{v_1^-}'v_1^+) - \frac{2(2+k)(5+2k)}{3(1+k)}J''.
\end{align*}
By direct computation, we check that the element $\Omega_3$ satisfies the following properties.
\begin{proposition}\label{cl-parameter} 
    For $k\neq-1,-2,-3,-4$, the element $\Omega_3$ in $\W^k(\sll_4,f_{2,2})$ satisfies
    \begin{align*}
        \lbr{J}{\Omega_3}&=\lbr{L_1}{\Omega_3}=0,\qquad\lbr{L_2}{\Omega_3}=3\Omega_3\lm{1}+\Omega_3',\\
        \lbr{\Omega_3}{\Omega_3}&=n_k\left(\frac{c_{4,2}(k)}{3}\lm{5}+2L_2\lm{3}+L_2'\lm{2}+W_4\lm{1}+\frac{1}{2}W_4'-\frac{1}{12}L_2^{(3)}\right)
    \end{align*}
    with $n_k=\frac{(2+k)(3+k)(4+k)^2}{(1+k)}$ and $W_4=\tfrac{1}{n_k}\Omega_3{}_{(1)}\Omega_3$. 
    Moreover, the element $W_4$ is of conformal weight four and satisfies 
    \begin{equation*}
        \Omega_{3}{}_{(3)}W_4=-(-31+16(2+c_{4,2}(k))\lambda_{4,2}(k))\Omega_3.
    \end{equation*}
\end{proposition}
\noindent
Therefore, the element $\Omega_3$ satisfies the same OPEs up to normalization as the weak generator $W_3$ in the universal $\W_\infty$-algebra $\W_\infty[c,\lambda]$ with parameters $(c,\lambda)=(c_{4,2}(k),\lambda_{4,2}(k))$.
This suggests that $L_1$ and $\Omega_3$ generate a quotient of $\W_\infty[c,\lambda]$, namely the affine coset 
$$C^k(\sll_4,f_{1^2,2})=\Com{V^{k+2}(\gl_2),\W^k(\sll_{4}, f_{1^2,2})}.$$ 
Since $L_1$ and $\Omega_3$ are only the weak generators, we can not conclude at this point as we need to check all the relations of the actual strong generators, especially the truncation \eqref{truncation} at conformal weight 15. To avoid doing this directly, we relate  $\W^k(\sll_4,f_{2,2})$ to $\W^k(\sll_4,f_{1^2,2})$, whose OPEs are calculated in \S~\ref{2,2}-\ref{2,1,1}. 

Recall that $\W^\ell(\sll_2)$ at generic level $\ell$ is isomorphic to the simple Virasoro vertex algebra $\vir{c}$ of central charge $c=c_{2,0}(\ell)$ and that the category $\W^\ell(\sll_2)\mod$ of $C_1$-cofinite modules is semisimple with simple objects $L_{r,s}^{\ell}:=\vir{c,h_{r+1,s+1}}$ ($r,s\geq 0$) of highest weight 
\begin{align}\label{virasoro highest weihgt}
    h_{r+1,s+1}=\frac{((r+1)(\ell+2)-(s+1))^2-(\ell+1)^2}{4(\ell+2)},
\end{align}
see \cite{CJHRY}.
In particular, the modules $L_{0,s}^{\ell}$ are realized as the BRST reduction of the Weyl modules $\weyl_{s\varpi}^\ell$ over $V^\ell(\sll_2)$, that is $L_{0,s}^{\ell}\simeq H_{f_{2}}(\weyl_{s\varpi}^\ell)$, and obtained as the simple quotient of the Virasoro Verma module $M^\ell_{1,s+1}$ of highest weight $h_{1,s+1}$:
\begin{align}\label{Virasoro modules}
0\rightarrow L_{0,-s-2}^{\ell} \rightarrow M^\ell_{1,s+1} \rightarrow L_{0,s}^{\ell} \rightarrow 0.
\end{align}
For $v\in \W^k(\sll_4,f_{2,2})$, we denote by $\mathcal{U}(\mathrm{Vir}^2)v\subset \W^k(\sll_4,f_{2,2})$ the submodule generated by $v$ over the Virasoro vertex subalgebra $\langle L_2\rangle \simeq \W^{k+1}(\sll_2)$.
\begin{lemma}\label{strong generators are C1}
For the strong generators $v_{i}^\pm$ $(i=1,2)$, 
\begin{align*}
    \mathcal{U}(\mathrm{Vir}^2)v_{i}^\pm\simeq L_{0,1}^{\ell}
\end{align*}
as $\W^{k+1}(\sll_2)$-modules.
\end{lemma}
\proof
By the OPE formulas in \S \ref{2,2}, $v_{i}^\pm$ are highest weight vectors of highest weight 
$h_{1,2}=-\frac{2k+3}{4(k+3)}$ $(\ell=k+1)$. Hence $\mathcal{U}(\mathrm{Vir}^2)v_{i}^\pm$ is a quotient of the Virasoro Verma module 
$M^\ell_{1,2}$. One can check directly that $\mathcal{U}(\mathrm{Vir}^2)v_{i}^\pm$ cannot be isomorphic to $M^\ell_{1,2}$ since $v_{i}^\pm$ satisfies 
\begin{align*}
    (L_{2,-2}-(k+3)L_{2,-1}^2)v_{i}^\pm=0.
\end{align*} 
The assertion follows from \eqref{Virasoro modules}.
\endproof
The category $\W^\ell(\sll_2)\mod$ is a braided tensor category \cite{CJHRY} and in particular, we have
\begin{align*}
    L_{0,s_1}^{\ell}\boxtimes L_{0,s_2}^{\ell} \simeq \bigoplus_{\begin{subarray}\ s_3=|s_1-s_2|\\ s_3\equiv s_1+s_2 \text{ }(2)\end{subarray}}^{s_1+s_2} L_{0,s_3}^\ell,
\end{align*}
see \cite{CJHRY,FZ}. 
Recal that $\W^\ell(\sll_2)\mod$ is semisimple and that $\W^k(\sll_4,f_{2,2})$ is strongly generated by $L_1,L_2,J$ and $v_{i}^\pm$. It follows that $\W^k(\sll_4,f_{2,2})$ is an injective limit of $C_1$-cofinite $\W^{k+1}(\sll_2)$-modules.
In particular, $\W^k(\sll_4,f_{2,2})$ is semisimple as $\W^{k+1}(\sll_2)$-modules and only $L_{0,s}^{\ell}$ ($s\geq0$) appear in the direct summands. 
On the other hand, $v_{i}^\pm$ generates the Fock module $\pi^J_{\pm1}$, which has the highest weights $J_0=\pm1$.
Let us introduce the coset algebra 
\begin{align*}
\mathscr{C}^k:=\Com{\W^{k+1}(\sll_2)\otimes \pi^J,\W^k(\sll_4,f_{2,2})}
\end{align*}
Then we have obtained the following decomposition.

\begin{corollary}\label{first decomp of 2,2} 
For $k\notin \Q$, the $\W$-algebra $\W^k(\sll_4,f_{2,2})$ decomposes into
    \begin{align*}
        \W^k(\sll_4,f_{2,2})\simeq \bigoplus_{\begin{subarray}c
            a\geq 0,b\in \Z\\ a\equiv b\ (2)
        \end{subarray}} \mathscr{C}^k[a,b]\otimes L_{0,a}^{k+1}\otimes \pi^J_b
    \end{align*}
    as a $\mathscr{C}^k\otimes \W^{k+1}(\sll_2)\otimes \pi^J$-module, with some $\mathscr{C}^k$-modules $\mathscr{C}^k[a,b]$.
\end{corollary}
The inverse Hamiltonian reduction for the affine vertex algebra $V^\ell(\sll_2)$ is given by the existence of a conformal embedding \cite{Ad,Sem94}:
\begin{align}\label{eq:IHRsl2}
    \mu_{\mathrm{iHR}}\colon V^\ell(\sll_2)\hookrightarrow \W^\ell(\sll_2)\otimes \Pi[0],
\end{align}
where $\Pi[0]$ is the chiral differential operators for $\mathrm{GL}_1$ usually called the half-lattice vertex algebra. It is a lattice vertex algebra extension of a rank-$2$ Heisenberg algebra strongly and freely generated by the elements $x$, $y$ of respective norms $1$ and $-1$:
\begin{align*}
    \Pi[0]=\bigoplus_{n\in \Z} \pi^{x,y}_{n(x+y)}\subset V_{\Z x\oplus \Z y} 
\end{align*}
The map $\mu_{\mathrm{iHR}}$ is determined by 
\begin{align}\label{sl2 iHR}
    e\mapsto \hwt{x+y},\quad h\mapsto \tfrac{\ell}{2}(x+y)+(x-y),\quad f\mapsto \left((\ell+2)L-(\mathbf{b}_\ell^2+(\ell+1)\mathbf{b}_\ell' \right)\hwt{-(x+y)}
\end{align}
with $\mathbf{b}_\ell=-\frac{\ell}{4}(x+y)+\frac{1}{2}(x-y)$.
Under this embedding, we have an isomorphism of $V^\ell(\sll_2)$-modules
\begin{align}\label{identification of highest weight vectors}
\weyl^\ell_{a\varpi} \simeq \begin{cases}
    (L_{0,a}^\ell\otimes \Pi[0])^{\SL_2} & (a\text{: even}),\\
    (L_{0,a}^\ell\otimes \Pi[\tfrac{1}{2}])^{\SL_2} & (a\text{: odd})
\end{cases}
\end{align}
where $\Pi[\tfrac{1}{2}]$ is the (simple) $\Pi[0]$-module defined by $\Pi[\tfrac{1}{2}]=\bigoplus_{n\in \frac{1}{2}+\Z} \pi^{x,y}_{n(x+y)}$ and $M^{\SL_2}$ denotes the maximal $\SL_2$-submodule of the $\sll_2$-module $M$.
Using the larger vertex algebra $\Pi^{\frac{1}{2}}[0]:=\Pi[0]\oplus \Pi[\tfrac{1}{2}]$, we introduce the following vertex subalgebra 
    \begin{align*}
        \mathscr{V}^k:=(\W^k(\sll_4,f_{2,2})\otimes \Pi^{\frac{1}{2}}[0])^{\SL_2}, 
    \end{align*}
    which decomposes, as a $\mathscr{C}^k\otimes V^{k+1}(\sll_2)\otimes \pi^J$-module, into 
    \begin{align}\label{decomposition of iHR for 2,2}
        \mathscr{V}^k\simeq \bigoplus_{a,b} \mathscr{C}^k[a,b]\otimes \weyl^{k+1}_{a\varpi}\otimes \pi^J_b.
    \end{align}
On the other hand, the hook-type $\W$-algebra $\W^k(\sll_4,f_{1^2,2})$ decomposes into 
    \begin{align}\label{decomposition of 2,1,1}
        \W^k(\sll_4,f_{1^2,2}) \simeq \bigoplus_{a,b} C^k[a,b] \otimes \weyl^{k+1}_{a\varpi}\otimes \pi^J_b
    \end{align}
    as $C^k(\sll_4,f_{1^2,2})\otimes V^{k+1}(\sll_2)\otimes \pi^J$-modules with some $C^k(\sll_4,f_{1^2,2})$-modules $C^k[a,b]$.
    Here we keep the Heisenberg field symbol $J$ as they satisfy the same OPEs. 
One can relate these two vertex algebras by extending the map $\mu_{\mathrm{iHR}}$. 
\begin{proposition}\label{inverse Hamiltonian reduction for minimal to rectangular}
    There exists an embedding of vertex algebras 
    \begin{align*}
        \widehat{\mu}_{\mathrm{iHR}}\colon \W^k(\sll_4,f_{1^2,2})\hookrightarrow \mathscr{V}^k\subset \W^k(\sll_4,f_{2,2})\otimes \Pi^{\frac{1}{2}}[0],
    \end{align*}
    which satisfies
    \begin{align*}
        &L\mapsto L_1,\quad J\mapsto J,\quad h \mapsto \tfrac{k+1}{2}(x+y)+(x-y),\\
        &e\mapsto \hwt{x+y},\quad f\mapsto \left((k+3)L_2-(\mathbf{b}_{k+1}^2+(k+2)\mathbf{b}_{k+1}' \right)\hwt{-(x+y)}\\
        &v_{1}^\pm\mapsto v_1^{\pm}\hwt{\frac{1}{2}(x+y)},\quad 
        v_2^\pm\mapsto (\pm v_2^\pm-(\mathbf{b}_{k+1}\pm \tfrac{k+3}{k+2}J)v_1^\pm+(k+3)v_1^\pm{}')\hwt{-\frac{1}{2}(x+y)}.
    \end{align*}
    for $k\neq -2,-3,-4$.
\end{proposition}
\noindent
The proof is based on the OPE formulas at $k\neq-2,-3,-4$ which extends to all levels.

We are now in a position to prove the main result of this subsection.
\proof[Proof of Theorem \ref{thm: the case 2,2}]
By applying the BRST reduction to the embedding $ \widehat{\mu}_{\mathrm{iHR}}$ in Proposition \ref{inverse Hamiltonian reduction for minimal to rectangular}, we obtain a vertex algebra homomorphism 
\begin{align*}
    [\widehat{\mu}_{\mathrm{iHR}}]\colon H_{f_2}(\W^k(\sll_4,f_{1^2,2}))\rightarrow H_{f_2}(\W^k(\sll_4,f_{2,2})\otimes \Pi^{\frac{1}{2}}[0]).
\end{align*}
As $\widehat{\mu}_{\mathrm{iHR}}(e)=\hwt{x+y}\in \Pi[0]$, we have 
$$H_{f_2}(\W^k(\sll_4,f_{2,2})\otimes \Pi^{\frac{1}{2}}[0])=\W^k(\sll_4,f_{2,2})\otimes H_{f_2}(\Pi^{\frac{1}{2}}[0]).$$
By \cite[Prop. 7]{ACGY}, the cohomology $H_{f_2}(\Pi^{\frac{1}{2}}[0])$ is isomorphic to the group algebra $\C[\Z_2]$, that is 
\begin{align}\label{group algebra appears}
    \C[\Z_2]\xrightarrow{\simeq} H_{f_2}(\Pi^{\frac{1}{2}}[0]),\quad \bar{0}\mapsto [\mathbf{1}],\ \bar{1}\mapsto [\hwt{\tfrac{1}{2}(x+y)}]
\end{align}
where $\bar{0}$ is the unit and $\bar{1}$ the generator satisfying $\bar{1}^2=\bar{0}$.
Actually, we have the following equivalences of cohomology classes
$$[\hwt{n(x+y)}]=[\mathbf{1}],\quad [\hwt{(\frac{1}{2}+n)(x+y)}]=[\hwt{\tfrac{1}{2}(x+y)}]\quad (n\in\Z).$$
Notice that $\W^k(\sll_4,f_{1^2,2})$ is $\Z_2$-graded by $\frac{1}{2}h_0$-eigenvalues modulo $\Z$.
Comparing this with the embedding in Proposition \ref{inverse Hamiltonian reduction for minimal to rectangular}, one finds that $[\widehat{\mu}_{\mathrm{iHR}}]$ preserves the $\Z_2$-gradings. Thus, $[\widehat{\mu}_{\mathrm{iHR}}]$ induces a vertex algebra homomorphism 
\begin{align*}
    [\widehat{\mu}_{\mathrm{iHR}}]\colon H_{f_2}(\W^k(\sll_4,f_{1^2,2}))\rightarrow \W^k(\sll_4,f_{2,2}).
\end{align*}

(1) We show that $[\widehat{\mu}_{\mathrm{iHR}}]$ is an isomorphism. 
Since we have the equality of the $q$-characters by Theorem \ref{Conj at the level of q}, it suffices to show the surjectivity of $[\widehat{\mu}_{\mathrm{iHR}}]$. 
By the OPE formulas in \S \ref{2,2}, one finds that $\W^k(\sll_4,f_{2,2})$ is weakly generated by the fields $J,L_1,L_2, v_1^\pm$ as $v_2^\pm=(L_1)_{-1}v_1^\pm$.
By \eqref{group algebra appears}, it is clear that the cohomology classes $[J]$ and $[L]$ map to $J$ and $L_1$ respectively. 
By construction, $[\widehat{\mu}_{\mathrm{iHR}}]$ restricts to $\W^{k+1}(\sll_2)\rightarrow \W^k(\sll_4,f_{2,2})$ whose image contains $L_2$. 
Since $v_1^\pm$ are the highest weight vectors of the Weyl modules $\weyl_{\varpi}^{k+1}\subset \W^k(\sll_4,f_{1^2,2})$, they define cohomology classes $[v_1^\pm]$ which map to the cohomology classes $[v_1^\pm\hwt{\frac{1}{2}(x+y)}]$ through $[\widehat{\mu}_{\mathrm{iHR}}]$, that is, $v_1^\pm$ in $\W^k(\sll_4,f_{2,2})$. This completes the proof.

(3) By (1) and Corollary \ref{first decomp of 2,2}, we have $\mathscr{C}^k[a,b]\simeq C^k[a,b]$ for all $a,b$ at generic levels. Hence the assertion follows from $\mathscr{C}^k\simeq C^k(\sll_4,f_{1^2,2})$ corresponding to $(a,b)=(0,0)$.

(2) By Proposition \ref{inverse Hamiltonian reduction for minimal to rectangular}, we have the embedding 
$\W^k(\sll_4,f_{1^2,2})\hookrightarrow (\W^k(\sll_4,f_{2,2})\otimes \Pi[0])^{\SL_2}$ for all levels. 
Recall that the $q$-character of $\W^k(\sll_4,f_{1^2,2})$ is independent of the level $k$ (see e.g. Proposition \ref{strong generating type}), so is $\W^k(\sll_4,f_{2,2})\otimes \Pi[0]$ and the subalgebra $(\W^k(\sll_4,f_{2,2})\otimes \Pi[0])^{\SL_2}$ as one can see from the formula \eqref{sl2 iHR}. 
Hence, it suffices to show $\W^k(\sll_4,f_{1^2,2})\simeq (\W^k(\sll_4,f_{2,2})\otimes \Pi[0])^{\SL_2}$ for generic levels. In the latter case, it follows from $\mathscr{C}^k[a,b]\simeq C^k[a,b]$ for all $a,b$. This completes the proof.
\endproof

\subsection{The case \tWtype{5}{2,3}}
Here, we relate the $\W$-algebra $\W^k(\sll_5,f_{2,3})$ to the hook-type $\W$-algebra $\W^k(\sll_5,f_{1^2,3})$ and show the following results.
\begin{theorem}\label{thm: the case 3,2}\hspace{0mm}
\begin{enumerate}[wide, labelindent=0pt]
    \item We have an isomorphism of vertex algebras
    \begin{equation*}
        H_{f_{2}}(\W^k(\sll_5,f_{1^2,3}))\simeq \W^k(\sll_5,f_{2,3}).
    \end{equation*}
       \item We have an isomorphism of vertex algebras
    \begin{equation*}
        \W^k(\sll_5,f_{1^2,3})\simeq (\W^k(\sll_5,f_{2,3})\otimes \Pi[0])^{\SL_2}.
    \end{equation*}
    \item Moreover, if $k$ is generic (e.g. $k\notin\Q$), then we have a conformal embedding
    \begin{equation*}
        C^k(\sll_5,f_{1^2,3})\otimes \W^{k+2}(\sll_2)\otimes \pi^J\hookrightarrow \W^k(\sll_5,f_{2,3}).
    \end{equation*}
\end{enumerate} 
\end{theorem}
The proof is very similar to the one for Theorem \ref{thm: the case 2,2} so we only give the main ideas. 
By \S \ref{3,2}, $\W^k(\sll_5,f_{2,3})$ is strongly generated by the following elements:
\renewcommand{\arraystretch}{1.3}
\begin{align*}
    \begin{array}{|c|c|c|c|c|c|} \hline
      \text{Conformal weight} &1  & 3/2 & 2 & 5/2 & 3\\ \hline
      \text{Strong generator}  &J  & v_{1}^\pm  & L_1, L_2 &v_{2}^\pm  &\Omega_{1,3}\\ \hline
    \end{array}
\end{align*}
\renewcommand{\arraystretch}{1}
Here $L_1,L_2$ are two mutually-commuting Virasoro elements of central charges $c_{3,2}(k)$ and $c_{2,0}(k+2)$ respectively and they commute with the Heisenberg element $J$ with $\lbr{J}{J}=\frac{2}{5}(3k+10)\Lambda$.
The element $\Omega_{1,3}$ commutes with $L_2$, $J$ and satisfies 
\begin{align*}
    \lbr{L_1}{\Omega_{1,3}}=3\Omega_{1,3}\Lambda+\Omega_{1,3}',
\end{align*}
The elements $v^\pm_{1}$ are the highest weight vectors for $L_2$ and $J$ and generate
$$ L_{0,1}^{k+2}\otimes \pi^J_1,\quad  L_{0,1}^{k+2}\otimes \pi^J_{-1},$$
respectively. Then the elements $v_2^\pm$ are obtained as
\begin{align*}
    v_2^\pm=\pm (k+4)L_{1,-1}v_1^\pm.
\end{align*}
Similarly, by \S \ref{3,1,1}, $\W^k(\sll_5,f_{1^2,3})$ is strongly generated by the following elements:
\renewcommand{\arraystretch}{1.3}
\begin{align*}
    \begin{array}{|c|c|c|c|} \hline
      \text{Conformal weight} &1  & 2 &3 \\ \hline
      \text{Strong generator}  &e,h,f,J  & L_1, v_{1}^\pm, v_{2}^\pm  & \Omega_{1,3}\\ \hline
    \end{array}
\end{align*}
\renewcommand{\arraystretch}{1}
The fields $e,h,f$ generate the affine vertex algebra $V^{k+2}(\sll_2)$, $J$ defines the Heisenberg vertex algebra $\pi^J$ with $\lbr{J}{J}=\frac{2}{5}(3k+10)\Lambda$ as before, and $L_1$ the Virasoro vertex algebra of central charge $c_{3,2}(k)$, which give the conformal embedding
\begin{align*}
    \vir{c_{3,2}(k)}\otimes V^{k+2}(\sll_2)\otimes \pi\hookrightarrow \W^k(\sll_5,f_{1^2,3}).
\end{align*}
The field $\Omega_{1,3}$ commutes with $V^{k+2}(\sll_2)\otimes \pi$ and satisfies 
\begin{align*}
    \lbr{L_1}{\Omega_{1,3}}=3\Omega_{1,3}\Lambda+\Omega_{1,3}',
\end{align*}
and thus extends only the factor $\vir{c_{3,2}(k)}$.
The fields $v_{1}^\pm, v_{2}^\pm$ give bases of the Weyl modules for $V^{k+2}(\sll_2)\otimes \pi$ associated with the natural representation $\C^2$ and its dual $\overline{\C}^2$:
\begin{align*}
    \Span\{v_{1}^+, v_{2}^+\}\simeq \C^2\subset \weyl_\omega^{k+2}\otimes \pi^J_1,\quad \Span\{v_{1}^-, v_{2}^-\}\simeq \overline{\C}^2\subset \weyl_\omega^{k+2}\otimes \pi^J_{-1}.
\end{align*}

Again, we extend the conformal embedding \eqref{eq:IHRsl2} to define the inverse Hamiltonian reduction which relates $\W^k(\sll_5,f_{2,3})$ with $\W^k(\sll_5,f_{1^2,3})$.
\begin{proposition}\label{inverse Hamiltonian reduction for 1,1,3 to 2,3}
    There exists an embedding of vertex algebras 
    \begin{align*}
        \widehat{\mu}_{\mathrm{iHR}}\colon \W^k(\sll_5,f_{1^2,3})\hookrightarrow  \W^k(\sll_5,f_{2,3})\otimes \Pi^{\frac{1}{2}}[0],
    \end{align*}
    which satisfies 
    \begin{align*}
        &L_1\mapsto L_1,\quad \Omega_{1,3}\mapsto \Omega_{1,3},\quad J\mapsto J,\quad h \mapsto \tfrac{k+2}{2}(x+y)+(x-y),\\
        &e\mapsto \hwt{x+y},\quad f\mapsto \left((k+4)L_2-(\mathbf{b}_{k+2}^2+(k+3)\mathbf{b}_{k+2}' \right)\hwt{-(x+y)},\\
        &v_{1}^\pm\mapsto v_1^{\pm}\hwt{\frac{1}{2}(x+y)},\quad 
        v_2^\pm\mapsto (\pm v_2^\pm-(\mathbf{b}_{k+2}\pm \tfrac{5(k+4)}{2(3k+10)}J)v_1^\pm+(k+4)v_1^\pm{}')\hwt{-\frac{1}{2}(x+y)}.
    \end{align*}
     for $k\neq -{10}/{3},-4,-5$.
\end{proposition}
\noindent
The proof is again by direct computation based on the OPE formulas at $k\neq -{10}/{3},-4,-5$ in \S \ref{3,2}-\ref{3,1,1}, which extends to all levels by using the global forms. 
Now, Theorem \ref{thm: the case 3,2} is proven using the same argument as in the proof of Theorem \ref{thm: the case 2,2}.

\subsection{The case \tWtype{5}{1,2^2}}
Our last example is the $\W$-algebra $\W^k(\sll_5,f_{1,2^2})$, that we relate to the minimal $\W$-algebra $\W^k(\sll_5,f_{1^3,2})$. We show the following.
\begin{theorem}\label{thm: the case 2,2,1} \hspace{0mm}
\begin{enumerate}[wide, labelindent=0pt]
    \item We have an isomorphism of vertex algebras
    \begin{equation*}
        H_{f_{1,2}}(\W^k(\sll_5,f_{1^3,2}))\simeq \W^k(\sll_5,f_{1,2^2}).
    \end{equation*}
    \item We have an isomorphism of vertex algebras
    \begin{equation*}
        \W^k(\sll_5,f_{1^3,2})\simeq (\W^k(\sll_5,f_{1,2^2})\otimes\Pi^{\frac{1}{3}}
        [0]\otimes \beta\gamma)^{\SL_3}
    \end{equation*}
    \item Moreover, if $k$ is a generic level (e.g. $k\notin \Q$), then we have a conformal embedding of vertex algebra
    \begin{equation*}
       C^k(\sll_5,f_{1^3,2})\otimes\W^{k+1}(\sll_3,f_{1,2})\otimes\pi^J\hookrightarrow \W^k(\sll_5,f_{1,2^2}).
    \end{equation*}
\end{enumerate} 
\end{theorem}
In the theorem, $\Pi^{\frac{1}{3}}[0]$ denotes the following simple current extension of $\Pi[0]$:
$$\Pi^{\frac{1}{3}}[0]=\bigoplus_{n\in\frac{1}{3}\Z}\pi^{x,y}_{n(x+y)}.$$
\begin{remark}\label{CDO of GL1}
    For an integer $m\neq 0$, there is an isomorphism of vertex algebras 
    \begin{align*}
        \Pi[0]\xrightarrow{\simeq} \Pi^{\frac{1}{m}}[0]:=\bigoplus_{n\in\frac{1}{m}\Z}\pi^{x,y}_{n(x+y)},
    \end{align*}
    satisfying 
    \begin{align*}
        \hwt{n(x+y)}\mapsto \hwt{\frac{n}{m}(x+y)},\ x+ y\mapsto \frac{1}{m}(x+y),\quad x- y\mapsto m(x-y).
    \end{align*}
\end{remark}
By \S\ref{2,2,1}, the $\W$-algebra $\W^k(\sll_5,f_{1,2^2})$ is strongly generated by the fields
\renewcommand{\arraystretch}{1.3}
\begin{align*}
    \begin{array}{|c|c|c|c|} \hline
      \text{Conformal weight} &1  & 3/2 & 2\\ \hline
      \text{Strong generator} &e_*,h,f_*,J & v^\pm_1, v^\pm_2  & L_1,L_2, E,F\\ \hline
    \end{array}
\end{align*}
\renewcommand{\arraystretch}{1}
where $L_1,L_2$ are two commuting Virasoro elements of central charge $c_{2,3}(k)$ and $c_{2,1}(k)$ respectively. 
Moreover, there is another pair of mutually-commuting Virasoro elements 
\begin{align*}
 &\mathbb{L}_1=L_2-\tfrac{1}{4+k}(L_1+L_2)-\tfrac{1}{2(4+k)}(\tfrac{1}{(5+2 k)}h^2-2e_*f_*+h'),\\
 &\mathbb{L}_2=L_1+\tfrac{1}{4+k}(L_1+L_2)+\tfrac{1}{2(4+k)}(\tfrac{1}{(5+2 k)}h^2-2e_*f_*+h'),
\end{align*}
of the same central charges $c_{2,3}(k)$ and $c_{2,1}(k)$, satisfying 
\begin{align*}
    L_1+L_2=\mathbb{L}_1+\mathbb{L}_2.
\end{align*}
These four Virasoro elements commute with the Heisenberg elements $J$ and $h$ satisfying 
$$\lbr{J}{J}=\tfrac{2}{5}(2k+5)\lm{1},\quad\lbr{h}{h}=2(2k+5)\lm{1},\quad\lbr{J}{h}=0.$$

\begin{proposition}\label{prop:structure_2,2,1} 
For $k\neq -4,-5,-5/2$, there is an embedding of vertex algebras 
    \begin{align*}
        \W^{k+1}(\sll_3,f_{1,2})\hookrightarrow \W^k(\sll_5,f_{1,2^2})
    \end{align*}
    explicitely given by 
    $$J_{\BP}\mapsto \tfrac{1}{6}(h+5J),\quad L_{1,\BP}\mapsto\mathbb{L}_2,\quad v_{\BP}^+\mapsto v^+_1,\quad v_{\BP}^-\mapsto-v_2^-.$$
\end{proposition}
The Virasoro field $\mathbb{L}_1$ and the Heisenberg field $J_\perp=-\tfrac{1}{2}(h - J)$ commute with the image of $\W^{k+1}(\sll_3,f_{1,2})$ inside $\W^k(\sll_5,f_{1,2^2})$.
Note that the conformal vector $L$ of $\W^k(\sll_5,f_{1,2^2})$ decomposes into
    \begin{align*}
        L=&L_1+L_2+\tfrac{1}{4(5+2k)}(h^2+5J^2)
         =\mathbb{L}_1+(L_{\BP}+\tfrac{3}{2(5+2k)}J_{\BP}^2)+\tfrac{1}{6(5+2k)}J_\perp^2.
    \end{align*}
Therefore, Proposition~\ref{prop:structure_2,2,1} implies the conformal embedding
\begin{align*}
    \vir{c_{2,3}(k)}\otimes \W^{k+1}(\sll_3,f_{1,2})\otimes \pi^{J_\perp}\hookrightarrow \W^k(\sll_5,f_{1,2^2}).
\end{align*}
In terms of this embedding, the remaining strong generators $\{f_*,v^+_2,F\}$ and $\{e_*,v^-_1,E\}$ give the extension part. 
Indeed, $f_*$ and $e_*$ are highest weight vectors satisfying
\begin{align*}
  \mathbb{L}_{1,n}f_*&=\delta_{n,0}\Delta f_*, &L_{\BP,n}f_*&=\delta_{n,0}\Delta_{\BP} f_*, &J_{\BP,n}f_*&=-\tfrac{1}{3}\delta_{n,0}f_*, &J_{\perp,n}f_*&=\delta_{n,0}f_*,\\
  \mathbb{L}_{1,n}e_*&=\delta_{n,0}\Delta f_*,  &L_{\BP,n}e_*&=\delta_{n,0}\Delta_{\BP} e_*, &J_{\BP,n}e_*&=\tfrac{1}{3}\delta_{n,0}e_*, &J_{\perp,n}e_*&=-\delta_{n,0}e_*,
\end{align*}
for $n\geq 0$ where $\Delta=\frac{(2+k)(10+3k)}{(4+k)(5+2k)}$ and $\Delta_{\BP}=-\frac{(4+3k)}{6(4+k)}$ and the other generators are obtained as
\begin{align}
\label{remaining fields1} v^+_{\BP,0}f_*&=-v^+_2,&v^-_{\BP,0}v^+_2&=-F-\tfrac{5(2+k)}{2(5+2k)}J_{\perp} f_*+\tfrac{(2+k)}{2}f_*',\\
\label{remaining fields2} v^-_{\BP,0}e_*&=v^-_1,&v^-_{\BP,0}v^-_1&=E+\tfrac{5(2+k)}{2(5+2k)}J_{\perp} e_*+\tfrac{(2+k)}{2}e_*'.
\end{align}

On the other hand, according to \S \ref{2,1,1,1} the $\W$-algebra $\W^k(\sll_5,f_{1^3,2})$ is strongly generated by the following elements.
\renewcommand{\arraystretch}{1.3}
\begin{align*}
    \begin{array}{|c|c|c|c|} \hline
      \text{Conformal weight} &1  & 3/2 & 2\\ \hline
      \text{Strong generator} &\begin{array}{c}
           e_{i,j} (1\leq i\neq j \leq 3)  \\
           h_1, h_2, J
      \end{array} & \begin{array}{c}
           v_1^+, v_2^+, v_3^+   \\
           v_1^-, v_2^-, v_3^-
      \end{array} & L_1\\ \hline
    \end{array}
\end{align*}
\renewcommand{\arraystretch}{1}
The $\W$-algebra $\W^k(\sll_5,f_{1^3,2})$ contains a vertex subalgebra isomorphic to $V^{k+1}(\sll_3)$ generated by the fields $h_1,h_2,e_{i,j}$ ($1\leq i\neq j \leq 3$) which commutes with the Heisenberg vertex algebra $\pi^J$ with $\lbr{J}{J}=\tfrac{3}{5}(5+2k)\lm{1}$ and the Virasoro vertex algebra with central charge $c_{2,3}(k)$ corresponding to $L_1$.
Hence, we have a conformal embedding
\begin{align*}
    \vir{c_{2,3}(k)}\otimes V^{k+1}(\sll_3)\otimes \pi\hookrightarrow \W^k(\sll_5,f_{1^3,2}).
\end{align*}
Finally, the fields $v_{i}^\pm$ ($1\leq i\leq 3$) give bases of the Weyl modules for $V^{k+1}(\sll_3)\otimes \pi^J$ associated with the natural representation $\C^3$ and its dual $\overline{\C}^3$:
\begin{align*}
\Span\{v_1^+,v_2^+,v^+_3\}\subset \weyl_{\varpi_1}^{k+1}\otimes \pi^J_1,\quad \Span\{v_1^-,v_2^-,v_3^-\}\subset \weyl_{\varpi_2}^{k+1}\otimes \pi^J_{-1}.
\end{align*}

To prove Theorem \ref{thm: the case 2,2,1}, we use the inverse Hamiltonian reduction which relates $V^\ell(\sll_3)$ with $\W^{\ell}(\sll_3,f_{1,2})$ established first in \cite{ACG}.
We use a slightly different embedding of vertex algebras given explicitly by:
\begin{align}
      \nonumber &\mu_{\mathrm{iHR}}\colon V^\ell(\sll_3)\hookrightarrow \W^\ell(\sll_3,f_{1,2})\otimes \Pi[0] \otimes \beta\gamma\\
      \nonumber &h_1\mapsto -J_{\BP}-\beta\gamma+\tfrac{2l}{3}(x+y)+\tfrac{1}{2}(x - y),\quad 
        h_2 \mapsto 2J_{\BP}+2 \beta\gamma-\tfrac{\ell}{3}(x+y)+\tfrac{1}{2}(x-y),\\
      \nonumber &e_{1,2}\mapsto \gamma\hwt{x+y},\quad 
        e_{2,3} \mapsto \beta,\quad 
        e_{1,3} \mapsto \hwt{x+y},\\
      \label{iHP for sl3} &e_{2,1} \mapsto \bigg( v^+_{\mathrm{BP}}- J_{\BP}\beta -\tfrac{\ell}{3} \beta (x+y) +\tfrac{1}{2} \beta (x-y)+(\ell+1) \beta' \bigg)\hwt{-(x+y)},\\
       \nonumber&e_{3,2}\mapsto v^-_{\mathrm{BP}} -2 J_{\BP}\gamma +\tfrac{\ell}{3} \gamma(x+y)- \tfrac{1}{2}\gamma(x-y) - \gamma^2\beta- \ell \gamma',\\
       \nonumber&e_{3,1} \mapsto \bigg( (3+\ell)L_{\BP}+\tfrac{(3-\ell)}{2(3+2\ell)}J_{\BP}^2+J_{\BP}(\x+\y+\beta \gamma)-(\x+\y)^2+(\x+\y) \beta\gamma\\
       \nonumber&\hspace{2cm}-(v^-_{\BP}\beta+v^+_{\BP}\gamma)-(1+\ell)\beta'\gamma-+(1+\ell)(\x+\y-\tfrac{1}{2}J_{\BP})'\bigg)\hwt{-(x+y)},
\end{align}
where we set 
$$\x=\tfrac{(-3+2\ell)}{6}x,\quad \y=\tfrac{(3+2\ell)}{6}y.$$
We extend the embedding to the $\W$-algebras $\W^k(\sll_5,f_{1^3,2})$ and $\W^k(\sll_5,f_{1,2^2})$.
\begin{proposition} \label{inverse Hamiltonian reduction for 2,1,1,1 to 2,2,1}
    There exists an embedding of vertex algebras 
    \begin{align*}
        \widehat{\mu}_{\mathrm{iHR}}\colon \W^k(\sll_5,f_{1^3,2})\hookrightarrow  \W^k(\sll_5,f_{1,2^2})\otimes \Pi^{\frac{1}{3}}[0]\otimes \beta\gamma,
    \end{align*}
    which extends $\mu_{\mathrm{iHR}}$ in \eqref{iHP for sl3} for the subalgebra $V^{k+1}(\sll_3)$ and satisfies 
    \begin{align*}
        &J\mapsto J_{\perp},\quad L_1\mapsto \mathbb{L}_1\\
        &v^+_1\mapsto f_*\hwt{\frac{2}{3}(x+y)},\quad v^+_2\mapsto (f_*\beta-v^+_2)\hwt{-\tfrac{1}{3}(x+y)},\\
        &v^+_3\mapsto \bigg(F+(\x+\y-\beta\gamma-\tfrac{(10+3k)}{2(5+2k)}J_{\BP}+\tfrac{5(2+k)}{6(5+2k)}J_{\perp})f_*+v^+_2\gamma -\tfrac{(2+k)}{2}f_*'\bigg)\hwt{-\frac{1}{3}(x+y)},\\
        &v^-_1\mapsto -e_*\hwt{\frac{1}{3}(x+y)},\quad v^-_2\mapsto (v^-_1-e_*\gamma)\hwt{\frac{1}{3}(x+y)},\\
        &v^-_3\mapsto \bigg(-E+(\x+\y-\tfrac{k}{2(5+2k)}J_{\BP}-\tfrac{5(2+k)}{2(5+2k)}J_{\perp})e_*-v^-_1\beta-\tfrac{(2+k)}{2} e_*'\bigg) \hwt{-\frac{2}{3}(x+y)}
    \end{align*}
    for $k\neq -5/2,-4,-5$. 
\end{proposition}
\noindent Again, the proof is by direct computation based on the OPE formulas at $k\neq -5/2,-4,-5$ and extends to all levels by using the global forms. 
To finish the proof of Theorem \ref{thm: the case 2,2,1}, we apply the BRST reduction $H_{f_{1,2}}$ to the embedding $\widehat{\mu}_{\mathrm{iHR}}$ in Proposition \ref{inverse Hamiltonian reduction for 2,1,1,1 to 2,2,1}. Then we have a vertex algebra homomorphism
\begin{align}\label{eq: 2,2,1 case}
    [\widehat{\mu}_{\mathrm{iHR}}]\colon  H_{f_{1,2}}(\W^k(\sll_5,f_{1^3,2}))\rightarrow  H_{f_{1,2}}(\W^k(\sll_5,f_{1,2^2})\otimes \Pi^{\frac{1}{3}}[0]\otimes \beta\gamma).
\end{align}
By \eqref{iHP for sl3}, we have 
\begin{align*}
   H_{f_{1,2}}(\W^k(\sll_5,f_{1,2^2})\otimes \Pi^{\frac{1}{3}}[0]\otimes \beta\gamma) \simeq \W^k(\sll_5,f_{1,2^2})\otimes H_{f_{1,2}}( \Pi^{\frac{1}{3}}[0]\otimes \beta\gamma).
\end{align*}
\begin{lemma}\label{lemma for 2,2,1}
There is an isomorphism of commutative vertex algebras 
\begin{align*}
   \C[\Z_3]\xrightarrow{\simeq} H_{f_{1,2}}( \Pi^{\frac{1}{3}}[0]\otimes \beta\gamma).
\end{align*}
\end{lemma}
\proof
To apply $H_{f_{1,2}}$, we take the nilpotent element $f=e_{3,1}$ and the good grading $\Gamma$ such that $\Gamma(e_{1,2})=0$, $\Gamma(e_{2,3})=1$ for the simple root vectors. Then the BRST complex $C_{f_{1,2}}(\Pi^{\frac{1}{3}}[0]\otimes \beta\gamma)$ is 
\begin{align*}
 C_{f_{1,2}}(\Pi^{\frac{1}{3}}[0]\otimes \beta\gamma)= \Pi^{\frac{1}{3}}[0]\otimes \beta\gamma \otimes \mathcal{F}(\sll_3, f_{1,2})
\end{align*}
where $\mathcal{F}(\sll_3, f_{1,2})\simeq bc^{\otimes2}$ is the two copies of the bc-system vertex superalgebra generated by the odd fields
\begin{align*}
    b_i(z)= \sum_{n\in \Z} b_{i,n}z^{-n-1},\quad c_i(z)= \sum_{n\in \Z} c_{i,n}z^{-n}\quad (i=1,2)
\end{align*}
satisfying $[b_{i,n},c_{j,m}]=\delta_{i,j}\delta_{n+m}$ and $[b_{i,n},b_{j,m}]=[c_{i,n},c_{j,m}]=0$, see \S\ref{sec:BRST}.
The complex $C_{f_{1,2}}(\Pi^{\frac{1}{3}}[0]\otimes \beta\gamma)$ is equipped with differential
\begin{align*}
    d&=\int Y((\mu_{\mathrm{iHR}}(e_{1,3})+1)c_1,z)dz+ \int Y(\mu_{\mathrm{iHR}}(e_{2,3})c_2,z)dz \\
     &=\underbrace{\int Y(\hwt{(x+y)}+1)c_1,z)dz}_{d_1}+ \underbrace{\int Y(\beta c_2,z)dz}_{d_2}.
\end{align*}
Then it is straightforward to see that the complex decomposes as 
\begin{align*}
    (C_{f_{1,2}}(\Pi^{\frac{1}{3}}[0]\otimes \beta\gamma \otimes bc^{\otimes2}),d)\simeq (\Pi^{\frac{1}{3}}[0]\otimes bc_1,d_1)\otimes (\beta\gamma \otimes bc_2,d_2)
\end{align*}
and that we have isomorphisms of commutative vertex algebras
\begin{align*}
  H^\bullet (\Pi^{\frac{1}{3}}[0]\otimes bc_1,d_1)\simeq \delta_{\bullet,0}\bigoplus_{a\in \frac{1}{3}\Z/\Z}\C [\hwt{a(x+y)}]\simeq \C[\Z_3],\quad   H^\bullet (\beta\gamma \otimes bc_2,d_2)\simeq \delta_{\bullet,0} \C[\mathbf{1}]\simeq \C.
\end{align*}
The assertion follows from the K\"{u}nneth formula. 
\endproof

Now, the proof of Theorem \ref{thm: the case 2,2,1} can be completed.
\proof[Proof of Theorem \ref{thm: the case 2,2,1}]
We apply Lemma \ref{lemma for 2,2,1} to \eqref{eq: 2,2,1 case} and obtain a homomorphism 
\begin{align*}
        [\widehat{\mu}_{\mathrm{iHR}}]\colon  H_{f_{1,2}}(\W^k(\sll_5,f_{1^3,2}))\rightarrow  \W^k(\sll_5,f_{1,2^2})\otimes \C[\Z_3].
\end{align*}
The $\W$-algebra $\W^k(\sll_5,f_{1^3,2})$ is equipped with a $\Z_3$-grading by setting 
\begin{align*}
    e_{i,j}=0,\quad h_i=0,\quad J=0,\quad L_1=0,\quad v^\pm_i=\pm\tfrac{1}{3}
\end{align*}
modulo $\Z$, which comes from the quotient group of the weight lattice by the root lattice for $\sll_3$. 
It induces the corresponding $\Z_3$-grading on $H_{f_{1,2}}(\W^k(\sll_5,f_{1^3,2}))$.
Since $\widehat{\mu}_{\mathrm{iHR}}$ in Proposition \ref{inverse Hamiltonian reduction for 2,1,1,1 to 2,2,1} preserves the $\Z_3$-gradings, it induces a vertex algebra homomorphism 
\begin{align*}
    [\widehat{\mu}_{\mathrm{iHR}}]\colon  H_{f_{1,2}}(\W^k(\sll_5,f_{1^3,2}))\rightarrow  \W^k(\sll_5,f_{1,2^2}).
\end{align*}
mapping
\begin{align*}
    [J] \mapsto J_{\perp},\ [L_1] \mapsto \mathbb{L}_1,\ [f_*\hwt{\frac{2}{3}(x+y)}] \mapsto v^+_1,\ [e_*\hwt{\frac{1}{3}(x+y)}]\mapsto -v^-_1.
\end{align*}
Hence the image contains all the weak generators of $\W^k(\sll_5,f_{1,2^2})$ by \eqref{remaining fields1}-\eqref{remaining fields2}. Therefore, $[\widehat{\mu}_{\mathrm{iHR}}]$ is surjective and thus is an isomorphism by \cite{CDM}. This completes the proof of (1). 
(2)-(3) are proved using the same argument as in the proof of Theorem \ref{thm: the case 2,2}.
\endproof

\section{Collapsing levels}\label{sec: collapsing levels}
In this section, we consider the simple quotient $\W_{k}(\sll_N,f)$ of the $\W$-algebras of type $A$ at levels where some strong generators belong to the maximal ideal of $\W^{k}(\sll_N,f)$. Their images drop out in the quotient.
The levels considered include the \emph{collapsing levels} \cite{AKMPP1, AFP, AvEM} which correspond to the case where all the strong generators of conformal weight greater than one belong to the maximal ideal of $\W^k(\sll_N,f)$. Then $\W_{k}(\sll_N,f)$ is rather trivial or isomorphic to the simple affine vertex algebra $L_{k^\sharp}(\g^\sharp)$ where $\g^\sharp$ is the centralizer of $f$ and $x$.

As we consider examples in low-rank cases, such levels can be detected by computing the Shapovalov form \cite{Li94}. The latter is obtained from the explicit OPEs and it degenerates if at least one of the strong generators drops out. We remove the critical level $k=-h^\vee$ from consideration, as the corresponding $\W$-algebras do not admit conformal vectors.

From the examples that we compute explicitly (see Tables \ref{Tab:Field contents for 2,2}--\ref{table:collaps_sl5_min}), we are able to deduce new families of collapsing levels (Theorems \ref{thm: collapsing for higherrank hooks} and \ref{thm:collapse_general}). We also discuss in Conjecture \ref{conj:reduction_simple_quotient} the compatibility between reduction by stages and simple quotients of $\W$-algebras.

\subsection{The case \tWtype{4}{2,2}} This case has already been studied in the literature \cite{AFP,AMP,CKLR19}.
We include the results for completeness of the paper. 
The fields $J, L_1,L_2,v^\pm_1,v^\pm_2$ appearing in \S~\ref{2,2} are not well-defined for all levels, but only for $k\neq -2,-3,-4$.
We twist them slightly to remove the singularities and the final set of global strong generators -- $\widehat{J}, \widehat{L}_1,\widehat{L}_2,\widehat{v}^\pm_1,\widehat{v}^\pm_2$ -- is well-defined for any non-critical level $k$. 
For instance, the explicit form of $L_1$ uses $-\tfrac{1}{(2+k)(k+4)}J^2$ in the Kac--Wakimoto construction \cite{KW04} so that $L_1$ commutes with $J$. The global form $\widehat{L}_1$ is obtained from $L_1$ by removing this term.
We obtain the following set of {quasi-primary}\footnote{The field $v$ is said \emph{quasi-primary} if it has a homogenous conformal weight and the total conformal vector acts by $L_{\mathrm{tot},1}v=0$.} fields
\begin{align*}
    &\widehat{J}=J,\quad \widehat{L}_1=(k+4)L_1+L_2+\tfrac{1}{(2+k)}J^2,\quad \widehat{L}_2=(k+3)L_2, \\
    &\widehat{v}^\pm_1=v^\pm_1,\quad \widehat{v}^\pm_2=v^\pm_2+\tfrac{(3+2k)(8+3k)}{2(2+k)}Jv^\pm_1.
\end{align*}

For $n\in\frac{1}{2}\Z$, denote by $\det(n)$ the determinant of the Shapovalov form restricted to the subspace of conformal weight $n$. For $n=1,2$, they are respectively given by
\begin{align*}
    &\det(1)=4(2+k)^3,\\
    &\det(2)=-1728(2+k)^{12}(4+k)(3+2k)^9(5+2k)(8+3k)^4. 
\end{align*}
Therefore, we are interested in the levels 
\begin{align*}
    k=-2,-\tfrac{8}{3},-\tfrac{5}{2},-\tfrac{3}{2}.
\end{align*}
In Table \ref{Tab:Field contents for 2,2}, we summarize the strong generators which survive (indicated by $\checkmark$) in the simple quotient. The dropped-out strong generators lie in the maximal ideal through certain corrections by differential polynomials of strong generators of strictly lower conformal weights. 
The symbols $*_1,*_2$ indicate linear relations in the simple quotient, namely:
\begin{align*}
    &*_1: \widehat{L}_1-\widehat{L}_2=0,\quad 
    *_2: \widehat{L}_1+\widehat{L}_2 - \widehat{J}^2+\tfrac{3}{2}(\widehat{J}'- \widehat{v}_1^+\widehat{v}_1^-)=0.
\end{align*}
\renewcommand{\arraystretch}{1.3}
\begin{table}[htbp]
\begin{tabular}{l|clllllll}
\hline
Level $k$ & $\widehat{J}, \widehat{v}_1^\pm$ & $\widehat{L}_1$ & $\widehat{L}_2$ & $\widehat{v}_2^\pm$  \\
\hline
{$-2$} &  & $*_1$ & $*_1$ &     \\[0.1em]
{$-8/3$} & \ok &  &  &   \\
{$-5/2$} & \ok  & $*_2$ & $*_2$   &  \ok   \\
{$-3/2$} & \ok &  &  &  &  \\
\hline
\end{tabular}
\caption{{Strong generators of $\W_k(\sll_4,f_{2,2})$.}}\label{Tab:Field contents for 2,2}
\end{table}
\renewcommand{\arraystretch}{1}

These levels provide interesting isomorphisms \cite{AMP}. Moreover, in  \cite{AMP}, the level $k=-\frac{7}{3}$ is also considered. It is not a collapsing level, but it is still special because the central charge of $L_2$ is zero.
\begin{proposition}[\cite{AMP}]\label{prop:collaps_sl4}
We have the isomorphisms of vertex algebras
\begin{align*}
    & \W_{-{3}/{2}}(\sll_4,f_{2,2})\simeq L_{1}(\sll_2), && \W_{-{8}/{3}}(\sll_4,f_{2,2})\simeq L_{-{4}/{3}}(\sll_2),\\
    & \W_{-{5}/{2}}(\sll_4,f_{2,2})\simeq \mathrm{FT}_1^0(\sll_2), && \W_{-{7}/{3}}(\sll_4,f_{2,2})\simeq (\W_{-{1}/{3}}(\osp_{3|2},f_{\mathrm{min}})\otimes \mathcal{F})^{\Z_2},\\
    & \W_{-2}(\sll_4,f_{2,2})\simeq \vir{1},
\end{align*}
where $\vir{c}$ denotes the simple Virasoro vertex algebra of central charge $c$ and $\mathcal{F}$ is the rank-one free fermion.
\end{proposition}

\begin{remark}
    In Proposition~\ref{prop:collaps_sl4}, the vertex algebra $\mathrm{FT}_1^0(\sll_2)$ is the singlet-type subalgebra of the affine Feigin--Tipunin algebra $\mathrm{FT}_p(\sll_2)$ with $p=1$ studied in \cite{ACGY,CNS}. The isomorphism $\W_{-{5}/{2}}(\sll_4,f_{2,2})\simeq \mathrm{FT}_1^0(\sll_2)$ was first noted in \cite[Cor.5.3]{C}.
    The simple $\W$-algebra $\W_{-5/2}(\sll_4,f_{2,2})$ decomposes into
    \begin{align*}
    \W_{-5/2}(\sll_4,f_{2,2})\simeq \Com{\pi,\beta\gamma^2} \simeq \bigoplus_{n\geq0}L_{-1}(\sll_2,n\alpha).
    \end{align*}
\end{remark}
    
Note that the simple $\W$-algebra $\W_{-7/3}(\sll_4,f_{2,2})$ also admits a decomposition in terms of the simple affine vertex algebra modules $L_{k}(\sll_2,\lambda)$ of highest weight $\lambda$ and the simple Virasoro modules $\vir{c,h}$ of highest weight $h$:
\begin{align*}
    \W_{-7/3}(\sll_4,f_{2,2})\simeq L_{-{2}/{3}}(\sll_2)\otimes \vir{\tfrac{1}{2},0}\oplus L_{-{2}/{3}}(\sll_2,\alpha)\otimes \vir{\tfrac{1}{2},\tfrac{1}{2}}.
\end{align*}

\subsection{The case \tWtype{5}{2,3}}
Applying the same reasoning as in the previous section, one can take the global quasi-primary strong generators to be 
\begin{align*}
    &\widehat{J}=J,\quad \widehat{L}_1=(k+5)L_1+L_2+\tfrac{25}{12(10+3k)}J^2,\quad \widehat{L}_2=(k+4)L_2,\quad \\
    &\widehat{\Omega}_{1,3}=\Omega_{1,3}-\tfrac{5(2+k)(5+k)}{2(10+3k)}JL_1-\tfrac{5}{2}JL_2-\tfrac{125(64+30k+3k^2)}{144(10+3k)^2}J^3,\\
    &\widehat{v}^\pm_1=v^\pm_1,\quad \widehat{v}^\pm_2=v^\pm_2+\tfrac{3(3+k)(15+4k)}{(10+3k)}Jv^\pm_1.
\end{align*}
The determinants of the Shapovalov form on the generators of conformal weights $n=1,2,3,3/2,5/2$ are given by 
\begin{align*}
    &\det(1)=-\tfrac{2}{5}(10+3k),\\
    &\det(2)=-\tfrac{144}{125}(3+k)^2(5+k)(5+2k)(10+3k)^4(15+4k),\\
    &\det(3)=\tfrac{2654208}{125}(2+k)(3+k)^6(5+k)^3(5+2k)^2(8+3k)(10+3k)^9(15+4k)^3,\\
    &\det(3/2)=-4(3+k)^2(10+3k)^2,\\
    &\det(5/2)=-\tfrac{9216}{25}(3+k)^8(5+2k)^2(10+3 k)^6(15+4 k)^2.
\end{align*}
Therefore, we are interested in the levels 
\begin{align*}
    k=-\tfrac{10}{3},-2,-3,-\tfrac{5}{2},-\tfrac{15}{4},-\tfrac{8}{3}.
\end{align*}
In Table \ref{Tab:Field contents for 3,2}, we summarize the strong generators that survive in the simple quotient. Again, the dropped-out strong generators lie in the maximal ideal through certain corrections by differential polynomials of strong generators of strictly lower conformal weights.
{ 
\renewcommand{\arraystretch}{1.3}
\begin{table}[htbp]
\begin{tabular}{l|llllllll}
\hline
Level $k$ & $\widehat{J}$ & $\widehat{L}_1$ & $\widehat{L}_2$ & $\widehat{\Omega}_{1,3}$ &$\widehat{v}^\pm_1$ & $\widehat{v}^\pm_2$  \\
\hline
{$-10/3$} &  &  &  &  &  &   \\[0.1em]
{$-3$} & \ok  & & \ok &  &  &   \\
{$-5/2$} & \ok & \ok &  &\ok  & \ok &   \\
{$-15/4$} & \ok &  & \ok & & \ok &   \\
{$-2$} &\ok  & \ok & \ok &  &\ok  & \ok  \\
{$-8/3$} &\ok & \ok &\ok &  & \ok & \ok  \\
\hline
\end{tabular}
\caption{{Field contents of $\W_k(\sll_5,f_{2,3})$.}}
\label{Tab:Field contents for 3,2}
\end{table}
\renewcommand{\arraystretch}{1}
}

In the below, the simple Virasoro vertex algebra $\vir{c}$ of central charge $c=c_{p,q}=1-\frac{6(p-q)^2}{pq}$ will be denoted by $\sL^{p,q}$ and its simple module $\vir{c,h}$ of highest weight $h=h_{r+1,s+1}$ (see \eqref{virasoro highest weihgt} by $\sL^{p,q}_{r,s}$ for simplicity.
We describe the $\W$-algebra $\W_k(\sll_5,f_{2,3})$ at collapsing levels.
\begin{proposition}\label{prop:collaps_sl5_32}
We have the isomorphisms of vertex algebras
\begin{align*}
    &\W_{-10/3}(\sll_5,f_{2,3})\simeq \C, &&\W_{-3}(\sll_5,f_{2,3})\simeq \sL^{1,1}\otimes \pi,\\
    &\W_{-5/2}(\sll_5,f_{2,3})\simeq \bigoplus_{n\in \Z}S_{2;n}\otimes \pi_{n},
    && \W_{-15/4}(\sll_5,f_{2,3}) \simeq \W_{-9/4}(\sll_3,f_{1,2}),\\
    &\W_{-8/3}(\sll_5,f_{2,3})\simeq\bigoplus_{a\in \Z_4} \mathscr{V}_a\otimes V_{2\ssqrt{5}\Z(1+a/4)}, && 
\end{align*}
where
\begin{align*}
    &\mathscr{V}_0=\sL^{4,7}\otimes \sL^{4,3}\oplus \sL^{4,7}_{0,2}\otimes  \sL^{4,3}_{0,2},\quad 
    \mathscr{V}_2=\sL^{4,7}_{0,0}\otimes \sL^{4,3}_{0,2}\oplus \sL^{4,7}_{0,2}\otimes \sL^{4,3}_{0,0},\quad 
    \mathscr{V}_{\pm1}=\sL^{4,7}_{0,1}\otimes \sL^{4,3}_{0,1}.
\end{align*}
\end{proposition}
Here we used the singlet algebra $\mathcal{M}(p)$ ($p\geq 2$) extending $\sL^{p,1}$ and its simple current modules $S_{p;n}:=\mathcal{M}_{n+1,1}(p)$, which are of conformal dimension $h_{n+1,1}$ and satisfy 
\begin{align}\label{simple currents for singlet}
    S_{p;n}\boxtimes S_{p;m}\simeq S_{p;n+m}, 
\end{align}
see \cite{CMY}.
It is interesting to note that $\W_{-15/4}(\sll_5,f_{2,3})$ is also an extension of the singlet algebra \cite{ACGY}
\begin{align*}
    \W_{-9/4}(\sll_3,f_{1,2})\simeq \bigoplus_{n\in \Z}S_{4;n}\otimes \pi_{n\ssqrt{-2}}.
\end{align*}

The isomorphisms at levels $k=-10/3,-3$ in Proposition~\ref{prop:collaps_sl5_32} can be deduced directly from the OPEs in \S\ref{3,2}. We prove the remaining cases.
\proof
\underline{Case $k=-15/4$}: The $\W$-algebra $\W_k(\sll_5,f_{2,3})$ is strongly generated by 
$J,\mathbb{L}=L_2-J^2, v_1^\pm$, which satisfy the OPEs
\begin{align*}
    &\lbr{J}{J}=-\tfrac{1}{2}\lm{1},\quad \lbr{\mathbb{L}}{\mathbb{L}}=-\tfrac{23}{4}\lm{3}+2\mathbb{L}\lm{1}+\mathbb{L}',\quad \lbr{J}{v_1^\pm}=\pm v_1^\pm
\end{align*} 
and 
\begin{align*}
    \lbr{\mathbb{L}}{v_1^\pm}&=\tfrac{3}{2}v_1^\pm\lm{1}+v_1^\pm{}'\pm4 v_2^\pm \\
                             &\equiv \tfrac{3}{2}v_1^\pm\lm{1}+v_1^\pm{}',\\
    \lbr{v_1^+}{v_1^-}&=\tfrac{15}{8}\lm{2}-\tfrac{15}{4}J\lm{1}+\tfrac{15}{4}J^2-\tfrac{5}{4}L_1-\tfrac{3}{4}L_2-\tfrac{15}{8}J'\\ 
                      & \equiv \tfrac{15}{8}\lm{2}-\tfrac{15}{4}J\lm{1}+3J^2+\tfrac{15}{8}J'-\tfrac{3}{4}\mathbb{L}.
\end{align*}
Here $\equiv$ is the equality in the simple quotient. They agree with the OPEs for $\W^\ell(\sll_3,f_{1,2})$ at level $\ell=-9/4$ (see \S\ref{2,1}). 
Thus we have a surjection $\W^\ell(\sll_3,f_{1,2})\twoheadrightarrow \W_k(\sll_5,f_{2,3})$ by \cite{DSK05}, which gives the isomorphism $\W_\ell(\sll_3,f_{1,2})\simeq \W_k(\sll_5,f_{2,3})$.

\noindent
\underline{Case $k=-5/2$}:
By direct computation, we have a homomorphism of vertex algebras 
\begin{align*}
    \W^{-5/2}(\sll_5,f_{2,3})\rightarrow \beta\gamma \otimes L_1(\sll_2) 
\end{align*}
which satisfies
\begin{align*}
    &L_1\mapsto L_{\beta\gamma},\quad L_2\mapsto 0,\quad J\mapsto h-\beta\gamma,\\
    &v^+_1\mapsto \gamma e,\quad v^-_1\mapsto -\tfrac{5}{2}\beta f,\\
    &v^+_2\mapsto -\tfrac{3}{2}(2\gamma'+\beta\gamma^2)e,\quad v^-_2\mapsto -\tfrac{15}{4}(2\beta'-\beta^2\gamma)f,\\
    &\Omega_{1,3}\mapsto \tfrac{35}{6}h^3+\tfrac{5}{12}(\beta\gamma)^3+\tfrac{15}{8}\beta^2 \gamma '\gamma+35 \beta\gamma ef -\tfrac{35}{2} \beta\gamma h^2-\tfrac{35}{2} \beta  \gamma h'+\tfrac{5}{8} \beta \gamma ''\\
    &\hspace{1cm}-35 e'f+\tfrac{35}{2} h'h-\tfrac{15}{8} \beta '\beta \gamma^2-\tfrac{5}{2}\beta '\gamma+\tfrac{5}{8} \beta ''\gamma+\tfrac{35}{6}h''
\end{align*}
where 
$$L_{\beta\gamma}=\beta\gamma'-\beta'\gamma+\frac{1}{2}(\beta\gamma)^2.$$
Since the image of $\W^{-5/2}(\sll_5,f_{2,3})$ commutes with $J_{\perp}:=h - 2 \beta\gamma$, we obtain the map
\begin{align*}
    \rho\colon \W^{-5/2}(\sll_5,f_{2,3})\rightarrow \Com{\pi^{J_\perp},\beta\gamma \otimes L_1(\sll_2)}.
\end{align*}

We use the following decomposition (see for instance \cite[Sect.\ 8]{CNS})
\begin{align}\label{betagamma}
    \beta\gamma \simeq \bigoplus_{n\in\Z}S_{2;n}\otimes \pi^y_{-n}, 
\end{align}
where the Heisenberg field $y$ has norm $-1$. 
Then, the fields
$L_{\beta\gamma}$ and $H=\beta\gamma$
identify respectively with the conformal vector of the singlet algebra $\mathcal{M}(2)$ and the Heisenberg vector $y$. 
Moreover,
\begin{align*}
    \beta\in S_{2;1}\otimes \pi^v_{-1},\quad \gamma\in S_{2;-1}\otimes \pi^v_{1},
\end{align*}
are the highest weight vectors. 
Hence
\begin{align*}
    \Com{\pi^{J_\perp},\beta\gamma \otimes L_1(\sll_2)}\simeq \bigoplus_{n\in \Z}S_{2;n}\otimes \pi_n.
\end{align*}

The coset $\com(\pi^{J_\perp},\beta\gamma \otimes L_1(\sll_2))$ is a simple current extension of $\mathcal{M}(2)\otimes \pi$, and thus simple as a vertex algebra. Since $\rho(L_1), \rho(J), \rho(v_1^\pm)$ weakly generate $\com(\pi^{J_\perp},\beta\gamma \otimes L_1(\sll_2))$, $\rho$ is surjective and
\begin{align*}
    \W_{-5/2}(\sll_5,f_{2,3})\simeq \Com{\pi^{J_\perp},\beta\gamma \otimes L_1(\sll_2)}\simeq\bigoplus_{n\in \Z}S_{2;n}\otimes \pi_{n}.
\end{align*}

\noindent\underline{Case $k=-8/3$}: 
Set $\mathscr{W}=\W_{-8/3}(\sll_5,f_{2,3})$. It is an exceptional $\W$-algebra and the branching rule can be obtained using the asymptotic data \cite{AvEM}. 
Since the central charges of $L_1, L_2$ are those of minimal models $c_{4,7}=13/16$ and $c_{4,3}=1/2$, the equality of the asymptotic growth
\begin{align*}
    g(\mathscr{W})=\tfrac{16}{7}=\tfrac{11}{14}+\tfrac{1}{2}+1=g(\sL^{4,7})+g(\sL^{4,3})+g(\pi^J)
\end{align*}
implies that the Virasoro vertex algebras generated by $L_1, L_2$ are indeed simple by \cite[Lem 2.8]{AvEM}.
Then the remaining strong generators $v_{1}^\pm$, $v_{2}^\pm$ generate the module $\sL^{4,7}_{0,1}\otimes \sL^{4,3}_{0,1}\otimes \pi_{\pm \ssqrt{5}/2}$. 
Considering the fusion rules of the Virasoro modules, we conclude that the possible $(\sL^{4,7}\otimes \sL^{4,3})$-submodules appearing in $\mathscr{W}$ are 
$\mathscr{L}_{a,b}:=\sL^{4,7}_{0,a}\otimes \sL^{4,3}_{0,b}$ ($a,b\in\{0,2\}$) and $\sL^{4,7}_{0,1}\otimes \sL^{4,3}_{0,1}$.
By using the strong generating type and the conformal weights appearing in $\mathscr{W}$, we obtain the decomposition 
\begin{align*}
    \mathscr{W}\simeq &(\mathscr{L}_{0,0}\oplus \mathscr{L}_{2,2}^{\oplus a_1})\otimes V_{\ssqrt{5}2N\Z}\oplus \mathscr{L}_{1,1}\otimes V_{\ssqrt{5}(2N\Z+1/2)}\\
                 &\oplus (\mathscr{L}_{0,2}^{\oplus a_2}\oplus \mathscr{L}_{2,0}^{\oplus a_3})\otimes V_{\ssqrt{5}(2N\Z+1)}\oplus \dots \oplus \mathscr{L}_{1,1}\otimes V_{\ssqrt{5}(2N\Z-1/2)}
\end{align*}
with $a_1\geq0, N>1$ and $a_2,a_3\in\{0,1\}$. The formula of the asymptotic growth in \cite[Prop 4.10]{AvEM} gives only one possible decomposition, which is the one in the assertion. This completes the proof.

\endproof

As the representation theory of the singlet algebras is now well-developed, one may study the representation theory of $\W_{-5/2}(\sll_4,f_{2,2})$ (see Proposition \ref{prop:collaps_sl4}) by using the vertex algebra extension. We find another such decomposition outside the collapsing levels.

\begin{proposition}\label{prop:singlet}
    We have the isomorphisms of vertex algebras
    \begin{align*}
        \W_{-11/3}(\sll_5,f_{2,3})\simeq \bigoplus_{n\in \Z}S_{2;n}\otimes S_{3;n}\otimes \pi_{\ssqrt{-5/2}n}.
    \end{align*}
\end{proposition}
\proof
By direct computation, we get a homomorphism of vertex algebras 
\begin{align*}
    \rho \colon\W^{-11/3}(\sll_5,f_{2,3})\rightarrow \Com{\pi^{H_\perp}, \beta\gamma\otimes L_{-4/3}(\sll_2)}
\end{align*}
which satisfies
\begin{align*}
    &L_1\mapsto L_{\beta\gamma},\quad L_2\mapsto L_{\sll_2},\quad J\mapsto H,\\
    &v^+_1\mapsto \beta e,\quad v^-_1\mapsto -\gamma f,\\
    &v^+_2\mapsto \tfrac{1}{3}( e\beta^2\gamma -2\beta ')e,\quad 
    v^-_2\mapsto -\tfrac{1}{3}(\beta \gamma^2+2 \gamma ')f,\\
    &\Omega_{1,3}\mapsto -\tfrac{7}{12} ef'+\tfrac{7}{8} hef+\tfrac{7}{32} h^3 +\tfrac{7}{12} e'f-\tfrac{7}{16}
h'h+\tfrac{7}{72}h''\\
    &\hspace{1cm}-\tfrac{20}{27}(\beta\gamma)^3-\tfrac{10}{3} \beta^2\gamma '\gamma -\tfrac{10}{9} \beta\gamma ''+\tfrac{10}{3} \beta '\beta\gamma^2+\tfrac{40}{9} \beta'\gamma '-\tfrac{10}{9} \beta ''\gamma
\end{align*}
where 
\begin{align*}
    L_{\sll_2}=\tfrac{3}{4}(\tfrac{1}{2}h^2+ef+fe)-\tfrac{3}{16}h^2
\end{align*}
is the conformal vector of the Heisenberg coset $\com(\pi^h,L_{-4/3}(\sll_2))\simeq \mathcal{M}(3)$ and $$H=\tfrac{3}{10}h - \tfrac{2}{5}\beta\gamma,\quad
H_\perp=h+2\beta\gamma\quad  \in \beta\gamma \otimes L_{-4/3}(\sll_2).$$

We use \eqref{betagamma} and the following decomposition \cite[Sect.~4.1]{CRW}
\begin{align}\label{B(P=3)}
    L_{-4/3}(\sll_2)\simeq \bigoplus_{n\in\Z}S_{3;n}\otimes \pi^h_{2n}
\end{align}
so that
\begin{align*}
\Com{\pi^{H_\perp}, \beta\gamma\otimes L_{-4/3}(\sll_2)}\simeq \bigoplus_{n\in \Z} S_{2;n}\otimes S_{3;n}\otimes \pi^H_n.
\end{align*}
The coset is a simple vertex algebra as it is a simple current extension of $\mathcal{M}(2)\otimes \mathcal{M}(3)\otimes \pi^H$.

The image of $\rho$ contains a set of weak generators, namely $\rho(L_1), \rho(L_2), \rho(J), \rho (v^\pm_1)$. Hence $\rho$ induces an isomorphism 
\begin{align*}
    \W_{-11/3}(\sll_5,f_{2,3})\simeq \Com{\pi^{H_\perp}, \beta\gamma\otimes L_{-4/3}(\sll_2)}\simeq \bigoplus_{n\in \Z} S_{2;n}\otimes S_{3;n}\otimes \pi^H_n.
\end{align*}
\endproof

The isomorphism in Proposition \ref{prop:singlet} seems to extend to the following family.
\begin{conjecture} 
We have the following decomposition
\begin{align*}
    \W_{-(p+2)+\frac{p+1}{p}}(\sll_{p+2},f_{2,p})\simeq \bigoplus_{n\in \Z}S_{2;n}\otimes S_{p;n}\otimes \pi_{\ssqrt{-\frac{p+2}{2}}n}.
\end{align*}
\end{conjecture}

\subsection{The case \tWtype{5}{1,2^2}}
One can take the global quasi-primary strong generators to be
\begin{align*}
    &\widehat{J}=J,\quad \widehat{h}=h,\quad\widehat{e}_*=e_*,\quad\widehat{f}_*=f_*,\\
    &\widehat{L}_1=(k+5)L_1+L_2+\tfrac{9}{16(5+2k)}h^2-\tfrac{5(9+4k)}{8(5+2k)}hJ+\tfrac{25}{16(5+2k)}J^2,\\
    &\widehat{L}_2=(k+4)L_2-\tfrac{(19+8k)}{16(5+2k)}h^2-\tfrac{5}{8(5+2k)}hJ+\tfrac{25}{16(5+2k)}
    J^2,\quad \\
    &\widehat{E}=E-\tfrac{5(2+k)(3+2k)}{4(5+2k)}Je_*,\quad \widehat{F}=Fw-\tfrac{5(2+k)(3+2k)}{4(5+2k)}Jf_*,\\
    &\widehat{v}^\pm_1=v^\pm_1,\quad \widehat{v}^\pm_2=v^\pm_2.
\end{align*}
The determinants of the Shapovalov form are then given by
\begin{align*}
    &\det(1)=-\tfrac{4}{5}(5+2k)^4,\\
    &\det(2)=-\tfrac{5435817984}{15625}(2+k)^{13}(5+k)(5+2k)^{18}(10+3k)^4,\\
    &\det(3/2)=(2+k)^4(5+2k)^4.
\end{align*}
Therefore, we are interested in the levels 
\begin{align*}
    k=-\tfrac{5}{2},-2,-\tfrac{10}{3}.
\end{align*}
In Table \ref{Tab:Field contents for 2,2,1}, we summarize the strong generators which survive in the simple quotient.
\renewcommand{\arraystretch}{1.3}
\begin{table}[htbp]
\begin{tabular}{l|ccccccccc}
\hline
Level $k$ & $e_*,h,f_*$ & $\widehat{J}$ & $\widehat{L}_1$ & $\widehat{L}_2$ & $\widehat{E}$ &$\widehat{F}$& $v_{1,2}^+$ &$v_{1,2}^-$ \\
\hline
{$-5/2$} &  &  &  &  &  &   \\[0.1em]
{$-2$} &\ok&\ok&  &  &  &   \\
{$-10/3$} & \ok&\ok &  &  & & & \ok&\ok  \\
\hline
\end{tabular}
\caption{{Field contents of $\W_k(\sll_5,f_{1,2^2})$.}}
\label{Tab:Field contents for 2,2,1}
\end{table}
\renewcommand{\arraystretch}{1}
\begin{proposition}\label{prop:collaps_sl5_221}
We have the isomorphisms
\begin{align*}
    &\W_{-{5}/{2}}(\sll_5,f_{1,2^2})\simeq \C, && \W_{-2}(\sll_5,f_{1,2^2})\simeq L_1(\sll_2)\otimes \pi,\\
    &\W_{-{10}/{3}}(\sll_5,f_{1,2^2})\simeq \bigoplus_{n\in \Z} \mathrm{FT}_3^n(\sll_2) \otimes \pi_{\ssqrt{-{3}/{2}}n}. &&
\end{align*}
\end{proposition}
\proof The cases $k=-5/2, -2$ are immediate by looking at the strong generators.
We show the case $k=-10/3$. We use the inverse Hamiltonian reduction \eqref{eq:IHRsl2}.
By \cite{CNS}, it extends to the singlet subalgebra of the affine Feigin--Tipunin algebra 
$$\mathrm{FT}_{p=3}^0(\sll_2)\hookrightarrow \mathcal{M}(3)\otimes \Pi[0].$$
By \eqref{B(P=3)}, one can replace $\mathcal{M}(3)$ with $L_{-4/3}(\sll_2)$ and obtain the decomposition 
\begin{align*}
    (L_{-4/3}(\sll_2)\otimes \Pi^{\tfrac{1}{2}}[0])^{\SL_2}\simeq \bigoplus_{n\in \Z} \mathrm{FT}_{p=3}^n(\sll_2)\otimes \pi^{H}_{2n}
\end{align*}
where $H$ corresponds to the semisimple element in $L_{-4/3}(\sll_2)$ satisfying $\lbr{H}{H}=-2/3\lm{1}$. 
Since the OPEs $\lbr{H}{H}$ and $\lbr{J}{J}$ agree, it suffices to construct a surjective homomorphism of vertex algebras  
\begin{align}\label{Collapsing map for 221 at -10/3}
\rho\colon \W_{-{10}/{3}}(\sll_5,f_{1,2^2})\twoheadrightarrow (L_{-4/3}(\sll_2)\otimes \Pi^{\frac{1}{2}}[0])^{\SL_2}.
\end{align}
By direct computation, one finds that $\rho$ is realized by setting
\begin{align*}
    &h\mapsto x-y-\tfrac{5}{6} (x+y),\quad J\mapsto H,\quad  e_*\mapsto \hwt{x+y},\\
    &f_*\mapsto (-\tfrac{1}{144}(11x-y)^2+\tfrac{1}{18}(11x-y)'+\tfrac{1}{2} ef+\tfrac{3}{16}h^2-\tfrac{1}{4}h')\hwt{-(x+y)},\\
    &v^+_1\mapsto e\hwt{\frac{1}{2}(x+y)},\quad v^-_1\mapsto f\hwt{\frac{1}{2}(x+y)},\\
    &v^+_2\mapsto \tfrac{1}{12} (-e(11x-y)+3he+4e')\hwt{\frac{1}{2}(x+y)},\\
    &v^-_2\mapsto \tfrac{1}{12} (-f(11x-y)-3hf+4f')\hwt{-\frac{1}{2}(x+y)},\\
    &L_1\mapsto 0,\quad E\mapsto 0,\quad F\mapsto 0,\\
    &L_2\mapsto \tfrac{1}{240}(11x-y)^2-\tfrac{1}{12}(11x-y)'+\tfrac{3}{2}ef+\tfrac{9}{16}h^2-\tfrac{3}{4}h'.
\end{align*}
This completes the proof.
\endproof

\subsection{Hook-type \tW-algebras}
Finally, we consider the collapsing levels for the hook-type $\W$-algebras $\W_k(\sll_{N},f_{1^m,n})$ for $N=2,\dots,5$.
We follow the same strategy as before and thus omit the details. 
Collapsing levels are listed in Tables \ref{table:collaps_sl2_pr}-\ref{table:collaps_sl5_min}.
Some of the collapsing levels already appear in the literature. They are gathered in the following theorem.
\begin{theorem}\label{thm:known_hook_type}
We have the following isomorphisms of vertex algebras.
\begin{enumerate}[wide, labelindent=0pt]
\item \cite[Thm. 8.6.]{AvEM} For $n\geq 2$,
\begin{align*}
    \W_{-n+\frac{n+1}{n}}(\sll_n,f_n)\simeq \C,\quad \W_{-n+\frac{n}{n-1}}(\sll_n,f_{1,n-1})\simeq \C.
\end{align*}
\item \cite[Cor. 6.5.]{CL1} For $n,m\geq 3$ with $(n,m)=1$, 
\begin{align*}
    \W_{-n+\frac{n+m}{n}}(\sll_n,f_n)\simeq \W_{-m+\frac{n+m}{m}}(\sll_m,f_m).
\end{align*}
\item \cite[Thm. 10.1]{L} For $n\geq 3$,
\begin{align*}
    &\W_{-n+\frac{n-1}{n}}(\sll_{n})\simeq \W_{-3+\frac{2}{3}}(\sll_3)\simeq \mathcal{M}(2).
\end{align*}
\item \cite[Thm 1.3]{AKMPP1} For $n\geq 3$, 
\begin{align*}
    \W_{-n+\frac{n-1}{n-2}}(\sll_n,f_{1,n-1})\simeq \pi.
\end{align*}
\item \cite{AFP,CGNS} For $n\geq 2$, 
\begin{align*}
    &\W_{-n+\frac{n+1}{n-1}}(\sll_n,f_{1,n-1})\simeq \begin{cases} \pi & (n\colon \text{odd}),\\ V_{\sqrt{n}\Z} & (n\colon \text{even}).\end{cases}
\end{align*}
\item \cite{CGNS} For $n\geq 2$ satisfying $(n+2,n-1)=1$,
\begin{align*}
    &\W_{-n+\frac{n+2}{n-1}}(\sll_n,f_{1,n-1})\simeq \sL^{3,n+2}\otimes V_{\sqrt{2n}\Z}\oplus \sL^{3,n+2}_{0,n}\otimes V_{\sqrt{n/2}+\sqrt{2n}\Z}.
\end{align*}
\item \cite{ACGY} The $\mathcal{B}^{(n+1)}$-algebra $\W_{-n+\frac{n}{n+1}}(\sll_n,f_{1,n-1})$ and the $\mathcal{R}^{(n-1)}$-algebra $\W_{-n+\frac{n}{n-1}}(\sll_n,f_{1^2,n-2})$ satisfy
\begin{align*}
    &\W_{-n+\frac{n}{n+1}}(\sll_n,f_{1,n-1})\simeq \bigoplus_{a\in \Z}S_{n+1;a}\otimes \pi_{\ssqrt{-\frac{n+1}{2}}a},\quad (n\geq 2),\\
    &\W_{-n+\frac{n}{n-1}}(\sll_n,f_{1^2,n-2})\simeq \bigoplus_{a\in \Z}S_{n-1;a}\otimes \pi_{\ssqrt{-\frac{n-1}{2}}a},\quad (n\geq 3).
\end{align*}
\item \cite{AFP} For $n,m\geq2$,
\begin{align*}
    &\W_{-(n+m)+\frac{n+m-1}{n-1}}(\sll_{n+m},f_{1^m,n})\simeq \pi\otimes L_{-m+\frac{m}{n-1}}(\sll_m),\\
    &\W_{-(n+m)+\frac{n+m}{n}}(\sll_{n+m},f_{1^m,n})\simeq L_{-m+\frac{m}{n}}(\sll_m).
\end{align*}
\item \cite[Cor 8.3]{AKFPP} The minimal $\W$-algebra $\W_{-2}(\sll_{5},f_{1^3,2})$ decomposes into
\begin{align*}
    &\W_{-2}(\sll_{5},f_{1^3,2})\simeq \bigoplus_{n\in \Z} L_{-1}(\sll_3,\lambda_n)\otimes \pi_{\ssqrt{5/3}n},
\end{align*} 
as $L_{-1}(\sll_3)\otimes \pi$-modules with $\lambda_n=n\varpi_1\ (n\geq 0)$, $-n\varpi_2\ (n\leq0)$.
\item \cite[Thm 9.5]{CL1} For $n\geq2$,
\begin{align*}
    &\W_{-n+\frac{n+1}{n}}(\sll_{n},f_{1,n-1})\simeq \beta\gamma^{\Z_n}\simeq \bigoplus_{p\in n\Z}S_{2;p}\otimes \pi^y_{-p}.
\end{align*}
\end{enumerate} 
\end{theorem}

The following isomorphisms are new as far as we know.
\begin{theorem}\label{thm: collapsing for hook type in cases}
For $(n,m)=(3,1), (4,1), (3,2)$,
\begin{align*}
    \W_{-n+1}(\sll_{n+m},f_{1^m,n})\simeq \pi\otimes \W_{-(n-1)+(m+1)}(\sll_{n-1}).
\end{align*}
Moreover, we have
\begin{align*}
    &\W_{-4+3/2}(\sll_4)\simeq \W_{-2+1/3}(\sll_2),\qquad
    &&\W_{-4+4/5}(\sll_4,f_{1,3})\simeq \W_{-3+3/5}(\mathfrak{sp}_4,f_{\mathrm{sub}}),\\
    &\W_{-4+5/4}(\sll_4,f_{1^2,2})\simeq (\beta\gamma\otimes V_{\sqrt{2}\Z})^{\Z_2},\qquad
    &&\W_{-4+4/3}(\sll_4,f_{1^2,2})\simeq\W_{-5+5/3}(\sll_5,f_{1,2^2}),\\
    &\W_{-5+7/3}(\sll_5,f_{1^2,3})\simeq\bigoplus_{a\in \Z_4} \widetilde{\mathscr{V}}_a\otimes V_{2\ssqrt{5}\Z(1+a/4)},
    &&\W_{-5+5/3}(\sll_5,f_{1^3,2})\simeq  \bigoplus_{n\in \Z}  L_{-3+2/3}(\sll_3,\lambda_n)\otimes \pi_{\ssqrt{-1}n},&&\\
    &\W_{-5+5/3}(\sll_5,f_{1,4})\simeq \bigoplus_{n\in\Z} S_{2;n}\otimes \pi_{\sqrt{3}n},&&
\end{align*}
where $\lambda_n=n\varpi_1\ (n\geq 0)$, $-n\varpi_2\ (n\leq0)$ and   
\begin{align*}
    &\widetilde{\mathscr{V}}_0=\sL^{4,7}\otimes L_{-2+4/3}(\sll_2,0)\oplus \sL^{4,7}_{0,2}\otimes  L_{-2+4/3}(\sll_2,2\varpi_1),\quad 
    \widetilde{\mathscr{V}}_{\pm1}=\sL^{4,7}_{0,1}\otimes L_{-2+4/3}(\sll_2,\varpi_1),\\
    &\widetilde{\mathscr{V}}_2=\sL^{4,7}_{0,0}\otimes L_{-2+4/3}(\sll_2,2\varpi_1)\oplus \sL^{4,7}_{0,2}\otimes L_{-2+4/3}(\sll_2,0).
\end{align*}
\end{theorem}

\proof
We show the isomorphisms case by case.

\noindent\underline{Case $\W_{-4+5/4}(\sll_4,f_{1^2,2})$}:
Using the asymptotic growth \cite[Cor 3.8, Prop 4.2]{AvEM}, we have an equality 
\begin{align*}
    g(\W_{-4+5/4}(\sll_4,f_{1^2,2}))=3=2+1= g(L_{-2+3/2}(\sll_2))+g(\pi),
\end{align*}
see the proof of Theorem \ref{thm:branching_rules_hook_type} below for the detailed computation.
Then it follows that $L_{-2+3/2}(\sll_2)\otimes \pi\hookrightarrow \W_{-4+5/4}(\sll_4,f_{1^2,2})$ by \cite[Thm 3.9]{AvEM}. Since the Kazhdan--Lusztig modules over $L_{-2+3/2}(\sll_2)$ are completely reducible with simple modules $L_{-2+3/2}(\sll_2,0)$, $L_{-2+3/2}(\sll_2,\varpi_1)$, the set of strong generators imply the decomposition 
\begin{align}\label{Proof 1:decomposition}
    \W_{-4+5/4}(\sll_4,f_{1^2,2})\simeq L_{-2+3/2}(\sll_2,0)\otimes V_{2\sqrt{2}\Z}\oplus L_{-2+3/2}(\sll_2,\varpi_1)\otimes V_{\sqrt{2}+2\sqrt{2}\Z},
\end{align}
and thus $\W_{-4+5/4}(\sll_4,f_{1^2,2})$ is a simple current extension of $L_{-2+3/2}(\sll_2,0)\otimes V_{2\sqrt{2}\Z}$. Recall that we have a conformal embedding 
\begin{align*}
L_{-2+3/2}(\sll_2)\hookrightarrow \beta\gamma,\quad e\mapsto \tfrac{1}{2}\beta^2,\ h\mapsto \beta\gamma,\ f\mapsto -\tfrac{1}{2}\gamma^2,
\end{align*}
and the decomposition 
\begin{align*}
    \beta\gamma\simeq L_{-2+3/2}(\sll_2,0)\oplus L_{-2+3/2}(\sll_2,\varpi_1),
\end{align*}
which is, again, a simple current extension. Hence, \eqref{Proof 1:decomposition} implies the desired isomorphism. 

\noindent
\underline{Case $\W_{-4+4/3}(\sll_4,f_{1^2,2})$}: One can use the realization \eqref{Collapsing map for 221 at -10/3} to show the isomorphism directly:
\begin{align*}
    \W_{-4+4/3}(\sll_4,f_{1^2,2})&\twoheadrightarrow (L_{-4/3}(\sll_2)\otimes \Pi^{\frac{1}{2}}[0])^{\SL_2} \xleftarrow[\rho]{\simeq} \W_{-{10}/{3}}(\sll_5,f_{1,2^2})\\
    e,h,f, J &\mapsto \rho(e_*), \rho(h),\rho(f_*), \rho(J),\\
    v_{1}^\pm, v_2^\pm &\mapsto \rho(v_{1}^\pm), \rho(v_{2}^\pm),\quad  L_1 \mapsto \rho(L_2).
\end{align*}

\noindent
\underline{Case $\W_{-8/3}(\sll_5,f_{1^2,3})$}: The decomposition follows from Proposition \ref{prop:collaps_sl5_32} and Theorem \ref{nice embedding} below.

\noindent
\underline{Case $\W_{-5+5/3}(\sll_5,f_{1^3,2})$}: By \cite[Thm 6.4]{AKFPP}, we have a conformal embedding 
$L_{-2+2/3}(\sll_3)\otimes \pi^{J}\hookrightarrow \W_{-5+5/3}(\sll_5,f_{1^3,2})$
and the complete reducibility 
\begin{align}
\W_{-5+5/3}(\sll_5,f_{1^3,2})\simeq \bigoplus_{n\in \Z}L_{-2+2/3}(\sll_3,\lambda_n)\otimes \pi^J_n
\end{align}
as a module over $L_{-2+2/3}(\sll_3)\otimes \pi^{J}$ for some highest weights $\lambda_n$ ($n\in \Z$).
We apply the BRST reduction $H_{f_{1,2}}$ to this decomposition and obtain
\begin{align}\label{second reduction for 2111}
H_{f_{1,2}}(\W_{-5+5/3}(\sll_5,f_{1^3,2}))\simeq \bigoplus_{n\in \Z} H_{f_{1,2}}(L_{-2+2/3}(\sll_3,\lambda_n))\otimes \pi^J_n.
\end{align}
By Theorem \ref{thm: the case 2,2,1}, $H_{f_{1,2}}(\W_{-5+5/3}(\sll_5,f_{1^3,2}))$ is a quotient of $\W^{-5+5/3}(\sll_5,f_{1,2^2})$.
On the other hand, we have $H_{f_{1,2}}(L_{-2+2/3}(\sll_3))\simeq \W_{-2+2/3}(\sll_3,f_{1,2})$ and $H_{f_{1,2}}(L_{-2+2/3}(\sll_3,\lambda_n))$ are simple or zero by \cite{Ar3}. 
It follows from \eqref{Collapsing map for 221 at -10/3} that the elements $e_*^n$ ($n\geq 0$) are non-zero in the simple quotient.
By using the OPEs in \S\ref{2,2,1}, we find that they are highest weight vectors for $\W_{-2+2/3}(\sll_3,f_{1,2})\otimes \pi^{J}$ such that 
\begin{align*}
    L_{\BP,0}e_*^n=\tfrac{(2n-3k-6)}{6(k+4)}e_*^n,\quad  J_{\BP,0}e_*^n=\tfrac{n}{3}e_*^n,\quad J_{0}e_*^n=ne_*^n.
\end{align*}
It follows that $\lambda_n=n\varpi_1$ $(n\geq0)$ and then $\lambda_n=-n\varpi_2$ $(n\leq0)$
by the self-duality of $\W_{-5+5/3}(\sll_5,f_{1,2^2})$. 
Now, the assertion follows from \eqref{second reduction for 2111} and $[J_\Lambda J]=-\Lambda$.

\noindent
\underline{Case $\W_{-5+5/3}(\sll_5,f_{1,4})$}: This case is a particular case of Theorem~\ref{thm:collapse_general} below.

\noindent The remaining cases can be proven directly by comparing the OPEs in Appendix \ref{appx1} and \cite{Fa} for $\W^k(\mathfrak{sp}_4,f_{\mathrm{sub}})$.
\endproof 

Two cases in Theorem \ref{thm: collapsing for hook type in cases} are generalized in higher rank cases.
\begin{theorem}\label{thm: collapsing for higherrank hooks}
For $m\geq 1$ and $n\geq 3$, there is an isomorphism of vertex algebras
\begin{align*}
    \W_{-(n+m)+(m+1)}(\sll_{n+m},f_{1^m,n})\simeq \pi\otimes \W_{-(n-1)+(m+1)}(\sll_{n-1}).
\end{align*}
\end{theorem}

\proof We divide the proof into the cases $m=1$ and $m\geq2$.\\
\underline{Case $m=1$}: Recall that the $\W$-algebra $\W^{k}(\sll_{n+1},f_{1,n})$ with $k=-n+1$ is generated by a Heisenberg element $J$, which gives the simple subalgebra $\pi^J$, the highest weight vectors $G^\pm$ of highest weights $\pm1$ and the elements of conformal weights $2,\dots,n$ lying in the coset $\com(\pi^J,\W^k(\sll_{n+1},f_{1,n}))$.
Combining \cite[Thm. 3.1, Eq. (A.1)]{FS} together with \cite{CGN,G}, the leading term of the OPE between $G^\pm$ is given by
\begin{align}\label{leading term of OPE for subreg}
    G^+(z)G^-(w)\sim \frac{\Lambda_{n-1}(n,k)}{(z-w)^{n+1}}+\dots,
\end{align}
with $\Lambda_{n}(n+1,k)=\prod_{i=1}^{n-1}(i(k+n)-1)=0$. 
Then the Shapovalov form implies that $G^\pm$ belong in the maximal ideal and thus their image is zero in the simple quotient $\W_{k}(\sll_{n+1},f_{1,n})$. 
Therefore, $\W_{k}(\sll_{n+1},f_{1,n})$ is the simple quotient of $\pi^J\otimes \mathrm{Com}(\pi^J,\W^k(\sll_{n+1},f_{1,n}))$. By \cite[Cor 6.5]{CL1},
$$\Com{\pi^J,\W_k(\sll_{n+1},f_{1,n})}\simeq \W_{-(n-1)+(m+1)}(\sll_{n-1}),$$
and the desired isomorphism follows.

\noindent
\underline{Case $m\geq2$}:
Recall that the universal $\W$-algebra $\W^{k}(\sll_{n+m},f_{1^m,n})$ (with $k=-n+1$) contains the affine vertex subalgebra $V^{k_\sharp}(\gl_m):=\pi^J\otimes V^{k+n-1}(\sll_m)$, the Weyl modules corresponding to the natural representation $\C^m$ and its dual $\overline{\C}^m$, and the elements of conformal weights $2,\dots,n$ lying in the affine coset $\mathrm{Com}(V^{k_\sharp}(\gl_m),\W^k(\sll_{n+1},f_{1,n}))$.
Since $k+n-1=0$, the Shapovalov form implies that the image of $V^{k+n-1}(\sll_m)$ in the simple quotient $\W_{k}(\sll_{n+m},f_{1^m,n})$ is trivial and so are the Weyl modules. Therefore, $\W_{k}(\sll_{n+m},f_{1^m,n})$ is the simple quotient of $\mathrm{Com}(V^{k_\sharp}(\gl_m),\W^k(\sll_{n+1},f_{1,n}))$, and thus is isomorphic to the desired vertex algebra, see \cite[Cor 6.5]{CL1} again.
\endproof

\begin{theorem}\label{thm:collapse_general}
For $n=3,5,\dots$, we have the following decomposition:
\begin{align*}
    \W_{-n+\frac{n}{n-2}}(\sll_{n},f_{1,n-1})\simeq \bigoplus_{p\in \Z} S_{2;p}\otimes \pi_{\sqrt{n-2}p}.
\end{align*}
\end{theorem}

\proof
The central charge of $\com(\pi^J,\W_k(\sll_{n},f_{1,n-1}))$ is $-2$. Since the universal Heisenberg coset $\mathrm{Com}(\pi^J,\W^k(\sll_{n},f_{1,n-1}))$ is of type $\W(2,\dots,n-1)$ \cite{CGN}, the simple quotient 
$$\mathcal{C}_k=\Com{\pi^J,\W_k(\sll_{n},f_{1,n-1})}$$
is either the Virasoro vertex algebra $\sL^{1,2}$ or the extension $\W_{-3+2/3}(\sll_3)\simeq\mathcal{M}(2)$ by \cite[Thm. 10.1]{L}. We show that $\mathcal{C}_k\simeq \mathcal{M}(2)$. 
Consider the decomposition
\begin{align*}
\W_k(\sll_{n},f_{q,n-1})\simeq \bigoplus_{p\in \Z} \mathcal{C}_{k,p}\otimes \pi_{\sqrt{n-2}p}
\end{align*}
as $\mathcal{C}_k\otimes \pi^J$-modules so that $\mathcal{C}_{k,p}$ are all simple $\mathcal{C}_k$-modules by \cite[Cor 5.6]{CGN}.
The vectors $G^\pm$ are highest weight vectors for $\mathcal{C}_k\otimes \pi^J$. 
It is straightforward to check that they generate the simple modules $\sL^{2,1}_{2,1}\otimes \pi_{\pm\sqrt{n-2}}$ or $G^\pm\in S_{2;\pm1}\otimes \pi_{\pm\sqrt{n-2}}$. 
Since the assignment $\pi_{\sqrt{n-2}p}\mapsto \mathcal{C}_{k,p}$ induces a braided-reverse equivalence between the categories of modules over $\pi^J$ and $\mathcal{C}_k$ \cite{CKM}, $\mathcal{C}_{k,p}$ must be all simple currents, which is the case only for $\mathcal{C}_k\simeq \mathcal{M}(2)$ by \cite{MY} and \eqref{simple currents for singlet}. Now, the desired decomposition follows.
\endproof

The decomposition for $n=3$, i.e. $\W_{0}(\sll_{3},f_{1,2})$, can also be recover using the explicit embedding
\begin{align*}
    &\W_0(\sll_3,f_{1,2})\hookrightarrow \beta\gamma\otimes L_1(\sll_2)\\
    &J\mapsto \beta\gamma+h,\quad L\mapsto L_{\beta\gamma}+\frac{1}{2}( \beta\gamma+h)^2,\\
    &G^+\mapsto \beta e,\quad G^-\mapsto \gamma f,
\end{align*}
where $L_{\beta\gamma}$ is defined in \eqref{betagamma}. By Proposition \ref{prop:collaps_sl5_32}, we obtain the following.
\begin{corollary}\label{cor: 3.2 v.s. 2.1}
    There is an isomorphism of vertex algebras
    \begin{align*}
        \W_{-5+5/2}(\sll_5,f_{2,3})\simeq \W_0(\sll_3,f_{1,2}).
    \end{align*}
\end{corollary}

The classification of simple $\W_0(\sll_3,f_{1,2})$-modules is established in \cite{AK} and is based on another realization which uses the rank-one symplectic fermion algebra $\mathcal{SF}$:
\begin{align*}
    \mathcal{SF}:=\langle \hwt{x},\partial\hspace{0.3mm} \hwt{-x}\rangle \simeq \mathrm{Ker}_{V_\Z} \int Y(\hwt{x},z)dz.
\end{align*}
These two realizations are related through the Kazama--Suzuki duality (see for instance \cite{CGNS}) 
\begin{align*}
    \mathcal{SF}\simeq \mathrm{H}_{\mathrm{rel}}(\widehat{\gl}_1; \beta\gamma\otimes V_\Z), 
\end{align*}
which induces the embedding in \cite{AK}: 
\begin{align*}
    \W_0(\sll_3,f_{1,2})&\simeq \mathrm{H}_{\mathrm{rel}}(\widehat{\gl}_1; \W_0(\sll_3,f_{1,2}))\\
                        &\hookrightarrow \mathrm{H}_{\mathrm{rel}}(\widehat{\gl}_1; (\beta\gamma\otimes V_\Z)\otimes V_\Z)\simeq  \mathcal{SF}\otimes V_\Z
\end{align*}
since $L_1(\sll_2)\simeq V_{\sqrt{2}\Z}\hookrightarrow V_{\Z^2}$.

Related to the remaining case of interest in Table \ref{table_collapsing_sub5}, i.e. $\W_{-5+5/7}(\sll_5,f_{1,4})$, there is a conformal embedding 
\begin{align*}
    \W_{-3+2/7}(\sll_3)\otimes \pi^J\hookrightarrow \W_{-5+5/7}(\sll_5,f_{1,4})
\end{align*}
by \cite[Cor 6.5]{CL1}. We conjecture the following decomposition.
\begin{conjecture}
    There is an isomorphism of $\W_{-3+2/7}(\sll_3)\otimes \pi$-modules
    \begin{align*}
        \W_{-5+5/7}(\sll_5,f_{1,4})\simeq \bigoplus_{n\in \Z}H_{f_3}(L_{-3+2/7}(\sll_3,\lambda_n))\otimes \pi_{\ssqrt{-7/3}n}
    \end{align*}
    where $\lambda_n=n\varpi_1\ (n\geq 0)$, $-n\varpi_2\ (n\leq0)$. 
\end{conjecture}
\renewcommand{\arraystretch}{1.1}
\begin{table}[htbp]
\begin{minipage}[l]{0.5\linewidth}
\centering
\scalebox{0.9}{
\begin{tabular}{c|c|c}
\hline
Level $k$ & $L$ & $\W_k(\sll_2,f_{2})$ \\
\hline
{$-1/2,-4/3$} & & $\C$\\
\hline
\end{tabular}}
\caption{{$\W_k(\sll_2,f_{2})$.}}
\label{table:collaps_sl2_pr}
\end{minipage}%
\begin{minipage}[l]{0.5\linewidth}
\centering
\scalebox{0.9}{
\begin{tabular}{c|cc|c}
\hline
Level $k$  & $L$ & $\Omega_3$ & $\W_k(\sll_3,f_{3})$ \\
\hline
{$-5/3,-9/4$}  & && $\C$\\
{$-4/3,-15/5$} & \ok && $\sL^{2,5}$\\
\hline
\end{tabular}}
\caption{{$\W_k(\sll_3,f_{3})$.}}
\end{minipage}
\end{table}

\begin{table}[htbp]
\begin{minipage}[l]{0.5\linewidth}
\centering
\scalebox{0.9}{
\begin{tabular}{c|ccc|c}
\hline
Level $k$ & $J$ & $v^\pm$ & $L$ & $\W_k(\sll_3,f_{1,2})$ \\
\hline
{$-3/2$}  &&&& $\C$\\
{$-1$} &  \ok &&& $\pi$\\
\hline
\end{tabular}}
\caption{{$\W_k(\sll_3,f_{1,2})$.}}
\end{minipage}%
\begin{minipage}[l]{0.5\linewidth}
\centering
\scalebox{0.9}{
\begin{tabular}{c|ccc|c}
\hline
Level $k$  & $L$ & $\Omega_3$ & $\Omega_4$ & $\W_k(\sll_4,f_{4})$ \\
\hline
${-11/4,-16/5}$  &&&& $\C$ \\
${-5/2,-10/3}$  &\ok&&& $\sL^{1,3}$\\
${-9/4,-24/7}$  &\ok&\ok&& $\W_{-2/3}(\sll_3,f_3)$\\
$-8/3,-13/4$  &\ok&\ok&& $\W_{-7/3}(\sll_3)$\\
\hline
\end{tabular}}
\caption{{$\W_k(\sll_4,f_{4})$.}}
\end{minipage}
\end{table}

\begin{table}[htbp]
\begin{minipage}[l]{0.5\linewidth}
\centering
\scalebox{0.9}{
\begin{tabular}{c|cccc|c}
\hline
Level $k$  & $J$ & $v^\pm$ & $L_1$ & $\Omega_{1,3}$ & $\W_k(\sll_4,f_{1,3})$ \\
\hline
$-8/3$ &&&& & $\C$ \\
$-5/2$ &\ok&&& & $\pi$ \\
{$-7/3$} &\ok&\ok&&& {$V_{2\Z}$} \\
{$-2$} &\ok&&\ok&& $\pi\otimes \sL^{1,2}$\\
$-16/5$ &\ok&\ok&\ok&& $\W_{-12/5}(\mathfrak{sp}_4,f_{\mathrm{sub}})$\\
\hline
\end{tabular}}
\caption{{$\W_k(\sll_4,f_{1,3})$.}}
\end{minipage}%
\begin{minipage}[l]{0.5\linewidth}
\centering
\scalebox{0.9}{
\begin{tabular}{c|cccc|c}
\hline
Level $k$  & $e,h,f$ & $J$ & $v^\pm_{1,2}$ & $L_1$ & $\W_k(\sll_4,f_{1,1,2})$ \\
\hline
$-1$  & &\ok&&& $\pi$\\
$-2$  & \ok&&&& $L_{-1}(\sll_2)$\\
{$-3/2$}  & \ok&\ok&\ok&& $(\beta\gamma\otimes V_{\sqrt{2}\Z})^{\Z_2}$ \\
{$-8/3$}  & \ok&\ok&\ok&& $\W_{-10/3}(\sll_5,f_{1,2,2})$\\
\hline
\end{tabular}}
\caption{{$\W_k(\sll_4,f_{1,1,2})$.}}\label{table:collaps_sl4_min}
\end{minipage}
\end{table}

\begin{table}[htbp]
\begin{minipage}[l]{0.5\linewidth}
\centering
\scalebox{0.9}{
\begin{tabular}{c|cccc|c}
\hline
Level $k$ & $L$ & $\Omega_3$ & $\Omega_4$ & $\Omega_5$ & $\W_k(\sll_5,f_5)$\\
\hline
${-19/5,-25/6}$ & &&&& $\C$ \\
${-18/5,-30/7}$ & \ok&&&& $\sL^{2,7}$\\
$-17/5,-35/8$ & \ok&\ok&&& $\W_{-1/3}(\sll_3,f_3)$\\
$-16/5,-40/9$ & \ok&\ok&\ok&& $\W_{-7/4}(\sll_4,f_4)$\\
$-15/4,-21/5$ & \ok&\ok&&& $\W_{-7/3}(\sll_3)$\\ \hline
\end{tabular}}
\caption{{$\W_k(\sll_5,f_5)$.}}
\end{minipage}%
\begin{minipage}[l]{0.5\linewidth}
\centering
\scalebox{0.9}{
\begin{tabular}{c|ccccc|c}
\hline
Level $k$  & $J$ & $L_1$ & $v^\pm$ & $\Omega_{1,3}$ & $\Omega_{1,4}$ & $\W_k(\sll_5,f_{1,4})$\\
\hline
$-15/4$ & &&&&& $\C$\\
$-7/2$ & \ok&&&&& $\pi$\\
$-11/3$ & \ok&&&&& $\pi$\\
$-3$ & \ok&\ok&&\ok&& $\pi\otimes\W_{-1}(\sll_3,f_3)$\\
$-13/4$ & \ok&\ok&\ok&&& \text{ext. of }$\sL^{3,7}\otimes V_{\sqrt{10}\Z}$\\
$-25/6$ & \ok&\ok&\ok&&& \text{ext. of }$\mathcal{M}(6)\otimes \pi$\\
$-10/3$& \ok&\ok&\ok&\ok&& \text{ext. of }$\mathcal{M}(2)\otimes \pi$\\
$-19/5$& \ok&\ok&\ok&\ok&& $\beta\gamma^{\Z_5}$\\
$-30/7$& \ok&\ok&\ok&\ok&& \\
\hline
\end{tabular}}
\caption{{$\W_k(\sll_5,f_{1,4})$.}}\label{table_collapsing_sub5}
\end{minipage}
\end{table}

\begin{table}[htbp]
\centering
\scalebox{0.9}{
\begin{tabular}{c|ccccc|c}
\hline
Level $k$  & $e,h,f$ & $J$ & $L_1$ & $v^\pm_{1,2}$ & $\Omega_{1,3}$ & $\W_k(\sll_5,f_{1,1,3})$\\
\hline
$-2$&  &\ok&\ok&&& $\pi\otimes \sL^{1,3}$\\
$-10/3$& \ok&&&&& $L_{-4/3}(\sll_2)$\\
$-3$&  \ok&\ok&&&& $\pi\otimes L_{-1}(\sll_2)$ \\
$-15/4$&  \ok&\ok&&\ok&& \text{ext. of }$\mathrm{FT}_4(\sll_2)\otimes \pi$ \\
$-8/3$&  \ok&\ok&\ok&\ok&&  \text{ext. of }$\sL^{4,7}\otimes L_{-2/3}(\sll_2)\otimes V_{2\sqrt{5}\Z}$ \\
\hline 
\end{tabular}}
\caption{{$\W_k(\sll_5,f_{1,1,3})$.}}\label{table:collaps_sl5_311}
\end{table}
\begin{table}[htbp]
\centering
\scalebox{0.9}{
\begin{tabular}{c|cccc|c}
\hline
Level $k$  & $\sll_3$ & $J$ & $v^\pm_{1,2,3}$ & $L_1$ & $\W_k(\sll_5,f_{1,1,1,2})$\\
\hline
$-1$&  &\ok&&& $\pi$\\
$-5/2$&  \ok&&&& $L_{-3/2}(\sll_3)$\\
$-2$&  \ok&\ok&\ok&& ext. of $L_{-1}(\sll_3)\otimes \pi$\\
$-10/3$& \ok&\ok&\ok&& ext. of $L_{-7/3}(\sll_3)\otimes \pi$ \\
\hline
\end{tabular}}
\caption{{$\W_k(\sll_5,f_{1,1,1,2})$.}}
\label{table:collaps_sl5_min}
\end{table}
\renewcommand{\arraystretch}{1}


Finally, note that most of the collapsing levels appearing in Propositions \ref{prop:collaps_sl4}, \ref{prop:collaps_sl5_32} and \ref{prop:collaps_sl5_221} are also collapsing levels for certain hook-type $\W$-algebras (see Tables \ref{table:collaps_sl2_pr}-\ref{table:collaps_sl5_min}).
This suggests that the reduction by stages descends to the simple $\W$-algebras. More generally, we conjecture the following phenomenon which is consistent with the long-standing conjecture by Kac--Roan--Wakimoto \cite{KRW03, KW08} asserting that $H_{f_\lambda}(L_k(\sll_n))$ is simple provided that it does not vanish.
\begin{conjecture}\label{conj:reduction_simple_quotient}
    Let $\lambda=(\lambda_1, \dots ,\lambda_n)\vdash N$ and set $\widehat{\lambda}_i=(1,\dots,1,\lambda_i)\vdash N_i$ with $N_i=N-\sum_{j>i}\lambda_{j}$.
    If for all $i=1,\dots,n$ the affine part of $H_{f_{\widehat{\lambda}_{i}}}H_{f_{\widehat{\lambda}_{i+1}}}\dots H_{f_{\widehat{\lambda}_n}}(V^k(\sll_N))$ does not belong to the maximal ideal then $H_{f_{\widehat{\lambda}_1}}H_{f_{\widehat{\lambda}_{2}}}\dots H_{f_{\widehat{\lambda}_n}}(L_k(\sll_N))$ is simple unless it vanishes:
    \begin{equation*}
        H_{f_{\widehat{\lambda}_1}}H_{f_{\widehat{\lambda}_{2}}}\dots H_{f_{\widehat{\lambda}_n}}(L_k(\sll_N))\simeq\left\{\begin{aligned}
            &\W_k(\sll_N,f_\lambda),\\&0.
        \end{aligned}\right.
    \end{equation*}
\end{conjecture}

\noindent In particular, Theorem~\ref{thm:known_hook_type} (9) gives a non-trivial example for the conjecture out of admissible levels:
$$H_{f_{1,2}}(\W_{-2}(\sll_5,f_{1^3,2}))\simeq \W_{-2}(\sll_5,f_{1,2^2})$$
since Proposition~\ref{prop:collaps_sl5_221} and Theorem~\ref{thm:known_hook_type} (4) give
\begin{align*}
   H_{f_{1,2}}(\W_{-2}(\sll_5,f_{1^3,2}))\simeq \bigoplus_{n\in \Z}\pi_{\ssqrt{1/3}n}\otimes \pi_{\ssqrt{5/3}n}\simeq V_{\sqrt{2}\Z}\otimes \pi\simeq \W_{-2}(\sll_5,f_{1,2^2}).
\end{align*}

\section{Exceptional \tW-algebras}\label{sec: exceptional levels}
In this section, we use our conjectures to explain the rationality of exceptional $\W$-algebras. 
Recall that $\W$-algebras at admissible exceptional levels have been proven to be rational recently.
\begin{theorem}[\cite{AvE22,McR21}]
	The simple $\W$-algebra $\W_k(\sll_N,f_{r,n^s})$ with $sn+r=N$ and $n>r\geq 0$ is rational at levels $$k=-h^\vee+\frac{p}{n},\qquad\gcd(p,n)=1,\,p\geq n$$
 where $h^\vee=N$ is the dual Coxeter number of $\sll_N$.
\end{theorem}
The case $r=0,s=1$ corresponds to the regular $\W$-algebras and was proven earlier by Arakawa \cite{Ar2}.
For convenience, we use the following notation:
\begin{align*}
    \W_{p,q}(\sll_N)=\W_{-N+{p}/{q}}(\sll_N,f_{N}).
\end{align*}
We show that, for appropriate admissible levels, the simple hook-type $\W$-algebras can be realized as extensions over a rational regular $\W$-algebra and an affine vertex algebra.
\begin{theorem}\label{thm:branching_rules_hook_type}
    For the admissible level $k=-h^\vee+p/n$, we have a conformal embedding 
    \begin{align*}
    \W_{k}(\sll_{n+m},f_{1^m,n})\hookleftarrow \W_{p,p-n}(\sll_{p-(n+m)})\otimes L_{k^\sharp}(\sll_m)\otimes \pi
    \end{align*}
    where $k^\sharp=-m+(p-n)/n$.
\end{theorem}
\proof
Let $C^k_{\mathrm{univ}}(\sll_{n+m},f_{1^m,n})$ be the universal affine coset of $\W^k(\sll_{n+m},f_{1^m,n})$, that is, the specialization of the universal coset $C^{\mathbf{k}}_{\mathrm{univ}}(\sll_{n+m},f_{1^m,n})\subset \W^{\mathbf{k}}(\sll_{n+m},f_{1^m,n})$ defined similarly over the polynomial ring $R=\C[\mathbf{k}]$ (see \S \ref{sec: W_infty-algebra}). 
Let $L_{\text{sug}}$ be the conformal vector of the affine vertex subalgebra $V^{\mathbf{k}^\sharp}(\gl_m)$ obtained after the localization to the quotient field $\mathbb{F}=\C(\mathbf{k})$.
Then the universal coset is characterized as the grading zero subalgebra: 
\begin{align*}
    C^{\mathbf{k}}_{\mathrm{univ}}(\sll_{n+m},f_{1^m,n})= \mathrm{Ker}\left(L_{\text{sug},0}\colon \W^{\mathbf{k}}(\sll_{n+m},f_{1^m,n})\rightarrow \W^{\mathbf{k}}(\sll_{n+m},f_{1^m,n})\otimes_{R}\mathbb{F} \right).
\end{align*}
On the other hand, $\W_{k}(\sll_{n+m},f_{1^m,n})$ at $k=-h^\vee+\frac{p}{n}$ is a Kazhdan--Lusztig object with respect to the affine vertex subalgebra $V^{k^\sharp}(\gl_m)$. Since $k^\sharp=-m+(p-n)/n$ satisfies $k^\sharp+m \in \mathbb{R}_{>0}$, the conformal weights $\Delta_\lambda=\frac{(\lambda,\lambda+\rho)}{2(k^\sharp+h^\vee)}$ of the Weyl modules $\mathbb{V}^{k^\sharp}_\lambda$ with highest weight $\lambda\in P_+$ is positive and $\Delta_\lambda=0$ iff $\lambda=0$. Hence the affine coset $C_{k}(\sll_{n+m},f_{1^m,n})$ of $\W_{k}(\sll_{n+m},f_{1^m,n})$ is again characterized as the grading zero subalgebra of the conformal vector $L_{\text{sug}}$. Therefore, $C_k(\sll_{n+m},f_{1^m,n})$ is a quotient of $C^k_{\mathrm{univ}}(\sll_{n+m},f_{1^m,n})$. By \cite[Cor 6.5]{CL1}, its (unique) simple quotient is $\W_{-h^\vee+\frac{p-n}{p}}(\sll_{p-(n+m)},f_{p-(n+m)})$, which is deduced from the coincidence of the parameters
\begin{align*}
    &c=-\frac{(1+m+n-p)(mn+n^2-p-np)(-n+mn+n^2+p-np)}{(n-p)p},\\
    &\lambda=-\frac{(n-p)p}{(2+m+n-p)(mn+n^2-2p-np)(-2n+mn+n^2+2p-np)}.
\end{align*}
Now, we compare the asymptotic growth for the both sides of the conformal embedding 
\begin{align*}
    \W_{k}(\sll_{n+m},f_{1^m,n})\hookleftarrow C_k(\sll_{n+m},f_{1^m,n})\otimes \overline{V}_{k^\sharp}(\sll_m)\otimes \pi
\end{align*}
where $\overline{V}_{k^\sharp}(\sll_m)$ is a quotient of $V^{k^\sharp}(\sll_m)$. By setting $g(V)$ to be the asymptotic growth of the vertex operator algebra $V$, 
it follows from \cite[Cor 3.8, Prop 4.2]{AvEM} that 
the asymptotic growth 
$$G=g(C_k(\sll_{n+m},f_{1^m,n}))+g(\overline{V}_{k^\sharp}(\sll_m))+g(\pi)$$
has an upper bound 
\begin{align*}
    G\leq g(\W_{k}(\sll_{n+m},f_{1^m,n}))=n-1+m^2+2m-\frac{(n+m)((n+m)^2-1)}{pn}
\end{align*}
and a lower bound 
\begin{align*}
     G&\geq g(\W_{p,p-n}(\sll_{p-(n+m)}))+g(L_{k^\sharp}(\sll_m))+g(\pi)\\
      &=\left((p-(n+m)-1)-\frac{(p-(n+m))((p-(n+m))^2-1)}{(p-n)p} \right)+\left(1-\frac{m}{(p-n)n}\right)(m^2-1)+1\\
      &=n-1+m^2+2m-\frac{(n+m)((n+m)^2-1)}{pn}=g(\W_{k}(\sll_{n+m},f_{1^m,n})).
\end{align*}
Therefore, we have
\begin{align*}
    g(C_k(\sll_{n+m},f_{n,1^m}))=g(\W_{p,p-n}(\sll_{p-(n+m)})),\quad g(\overline{V}_{k^\sharp}(\sll_m))=g(L_{k^\sharp}(\sll_m)).
\end{align*}
The latter implies $\overline{V}_{k^\sharp}(\sll_m)=L_{k^\sharp}(\sll_m)$ by \cite[Thm 3.9]{AvEM}. Then the complete reducibility of $L_{k^\sharp}(\sll_m)$-modules in the Kazhdan--Lusztig category \cite{Ar} implies that $C_k(\sll_{n+m},f_{1^m,n})$ is simple and thus 
\begin{align*}
C_k(\sll_{n+m},f_{1^m,n})\simeq \W_{p,p-n}(\sll_{p-(n+m)}).
\end{align*}
This completes the proof.
\endproof
Since $\W_{k}(\sll_{n+m},f_{1^m,n})$ is completely reducible as a module over $\W_{p,p-n}(\sll_{p-(n+m)})\otimes L_{k^\sharp}(\sll_m)$ by \cite{Ar2,Ar}, it decomposes into 
\begin{align*}
\W_{k}(\sll_{n+m},f_{1^m,n})\simeq \bigoplus_{(\lambda,\mu)\in A} \mathbb{L}_{p,p-n}(\lambda)\otimes L_{k^\sharp}(\mu)\otimes \mathscr{V}_{a_{\lambda,\mu}+L}.
\end{align*}
Here the sum runs over a finite set $A$ which parameterizes the simple $\W_{p,p-n}(\sll_{p-(n+m)})$-modules $\mathbb{L}_{p,p-n}(\lambda)$ and $L_{k^\sharp}(\sll_m)$-modules $L_{k^\sharp}(\mu)$ appearing in the decomposition of multiplicities, and $\mathscr{V}_L$ is the Heisenberg vertex algebra $\pi$ or a rank one lattice vertex subalgebra extending $\pi$.
We apply $H_{f_{1^{m-n'},n'}}$ with $n'\leq n,m$ to $\W_{k}(\sll_{n+m},f_{1^m,n})$ and obtain
\begin{align*}
    H_{f_{1^{m-n'},n'}}(\W_{k}(\sll_{n+m},f_{1^m,n}))\simeq \bigoplus_{(\lambda,\mu)\in A}\mathbb{L}_{p,p-n}(\lambda)\otimes H_{f_{1^{m-n'},n'}}(L_{k^\sharp}(\mu))\otimes \mathscr{V}_{a_{\lambda,\mu}+L}.
\end{align*}
Since 
$$H_{f_{1^{m-n'},n'}}(L_{k^\sharp}(\sll_m))\simeq \W_{k^\sharp}(\sll_m,f_{1^{m-n'},n'})$$
by \cite[Thm 7.8]{AvE22}, we may iterate the hook-type BRST reductions. 
As a by-product of Conjecture~\ref{conj:successive_reductions}, we obtain the following conformal embedding for exceptional $\W$-algebras of type $\sll_N$.
\begin{conjecture}\label{conj:branching_rules}
Let $N\in\Z_{>0}$ and $k=-N+p/n$ to be admissible. 
\begin{enumerate}[wide, labelindent=0pt]
    \item The $\W$-algebra $\W_k(\sll_N,f_{1^r,n^s})$, with $sn+r=N$ and $r\geq 0$, has a conformal embedding
\begin{equation*}
\W_k(\sll_N,f_{1^{r},n^s})\hookleftarrow\bigotimes_{\ell=0}^{s}\W_{{p-\ell n},{p-(\ell+1)n}}(\sll_{p-N})\otimes L_{-r+{(p-sn)}/{n}}(\sll_{r})\otimes\pi^{\otimes (s-\delta_{r,0})}.
\end{equation*}
    \item The exceptional $\W$-algebra $\W_k(\sll_N,f_{r,n^s})$, with $sn+r=N$ and $n>r\geq 0$, has a conformal embedding
\begin{equation*}
\W_k(\sll_N,f_{r,n^s})\hookleftarrow\bigotimes_{\ell=0}^{s}\W_{{p-\ell n},{p-(\ell+1)n}}(\sll_{p-N})\otimes \W_{-r+{(p-sn)}/{n}}(\sll_{r})\otimes\pi^{\otimes (s-\delta_{r,0})}.
\end{equation*}
Here we set $L_\bullet(\sll_i)=\C$ for $i=0,1$.
\end{enumerate}
\end{conjecture}

In particular, by Theorems~\ref{thm: the case 2,2}, \ref{thm: the case 3,2}, \ref{thm: the case 2,2,1}, the conjecture is true in the following cases.
\begin{theorem}\label{nice embedding} We have the following conformal embeddings.
\begin{enumerate}[wide, labelindent=0pt]
        \item For $k=-4+p/2$ $(p=5,7,9,\dots)$,
        \begin{align*}
            \W_{-4+{p}/{2}}(\sll_4,f_{2,2})\hookleftarrow\W_{p,p-2}(\sll_{p-4})\otimes \W_{p-2,p-4}(\sll_{p-4})\otimes\pi.
        \end{align*}
        \item For $k=-5+p/3$ $(p\geq5,\ (p,3)=1)$,
        \begin{align*}
            \W_{-5+{p}/{3}}(\sll_5,f_{2,3})\hookleftarrow\W_{p,p-3}(\sll_{p-5})\otimes \W_{p-3,3}(\sll_{2})\otimes\pi.
        \end{align*}
        \item For $k=-5+p/2$ $(p=5,7,9,\dots)$,
        \begin{align*}
           \W_{-5+{p}/{2}}(\sll_5,f_{1,2^2})&\hookleftarrow\W_{p,p-2}(\sll_{p-5})\otimes \W_{-3+{(p-2)}/{2}}(\sll_3,f_{1,2})\otimes\pi\\
            &\hookleftarrow\W_{p,p-2}(\sll_{p-5})\otimes \W_{p-2,p-4}(\sll_{p-5})\otimes\pi^{\otimes 2}.
        \end{align*}
    \end{enumerate}
\end{theorem}

Combining with the inverse Hamiltonian reduction embeddings presented in \S\ref{sec: Non-trivial examples in low ranks}, we observe that they are preserved when considering simple quotients at these levels.
\begin{corollary}\label{iHR for simple quotients}
We have the following isomorphism of vertex algebras.
\begin{enumerate}[wide, labelindent=0pt]
\item For $k=-4+p/2$ $(p=5,7,9,\dots)$, $\W_k(\sll_4,f_{1^2,2})\simeq (\W_k(\sll_4,f_{2,2})\otimes \Pi^{\frac{1}{2}}[0])^{\SL_2}$.
\item For $k=-5+p/3$ $(p\geq5,\ (p,3)=1)$, $\W_k(\sll_5,f_{1^2,3})\simeq (\W_k(\sll_5,f_{2,3})\otimes \Pi^{\frac{1}{2}}[0])^{\SL_2}$.
\item For $k=-5+p/2$ $(p=5,7,9,\dots)$, $\W_k(\sll_5,f_{1^3,2})\simeq (\W_k(\sll_5,f_{1,2^2})\otimes \Pi^{\frac{1}{3}}[0]\otimes \beta\gamma)^{\SL_3}$.
\end{enumerate}
\end{corollary}

The second part of Conjecture \ref{conj:branching_rules} indicates that
exceptional $\W$-algebras are extensions of rational regular $\W$-algebras.
It also suggests relations between module categories.

The conformal embedding 
\begin{align*}
    \W_{-(p-5)+\frac{p-3}{p}}(\sll_{p-5})\otimes \W_{-2+\frac{p-3}{3}}(\sll_{2})\otimes \pi\hookrightarrow \W_{-5+{p}/{3}}(\sll_5,f_{2,3})
\end{align*}
in (2) is closely related to the level-rank duality (see, e.g.\ \cite[Appendix B]{CGNS}):
\begin{align}\label{affine level-rank duality}
    L_{m}(\sll_n)\otimes L_{n}(\sll_m)\otimes \pi\hookrightarrow V_{\Z^{nm}},
\end{align}
with $(n,m)=(p-5,2)$ in our case. 
Indeed, by \cite[Rem. 10.7]{ACF}, the Kazhdan--Lusztig category of $\W_{-(p-5)+\frac{p-3}{p}}(\sll_{p-5})$ (resp.\ $ \W_{-2+\frac{p-3}{3}}(\sll_{2})$) is a quotient of the braided tensor category of admissible representations over $L_{-(p-5)+\frac{p-3}{p}}(\sll_{p-5})$ (resp.\ $ L_{-2+\frac{p-3}{3}}(\sll_{2})$), which has the same fusion rules as  $L_{2}(\sll_{p-5})$ (resp.\ $ L_{p-5}(\sll_{2})$). Note that the last level $2$ (resp.\ $p-5$) is obtained by replacing the denominator of the level $-(p-5)+\frac{p-3}{p}$ (resp.\ $-2+\frac{p-3}{3}$) with one.
The embedding \eqref{affine level-rank duality} induces the decomposition  
\begin{align*}
V_{\Z^{nm}}\simeq \bigoplus_{a\in \Z_{nm}}\left(\bigoplus_{\substack{\lambda\in \widehat{\mathrm{P}}_+^m(n)\\\Proj(\lambda)=a}} L_m(\sll_n,\lambda)\otimes {L}_n(\sll_m,\sigma^{\frac{a-\ell(\lambda)}{n}}(\lambda^t)) \right)\otimes V_{\frac{a}{\ssqrt{nm}}+\ssqrt{nm}\Z}.
\end{align*} 
We refer to \cite{CGNS} for the notation. 

For the simple $\W$-algebra $\W_k(\sll_N,f_\lambda)$ at admissible levels $k=-N+p/q$, let us introduce the (simple) $\W_k(\sll_N,f_\lambda)$-modules
\begin{align*}
    \mathbf{L}_{p,q}^{\lambda}(\mu,\nu):=H_{f_\lambda,\nu}(L_k(\sll_N,\mu)),
\end{align*}
which are the BRST reduction of the admissible representations $L_k(\sll_N,\mu)$ with integral highest weights $\mu\in \widehat{\mathrm{P}}_+^{p-N}(N)$ twisted by weights $\nu\in \widehat{\mathrm{P}}(N)$. The Feigin--Frenkel duality 
$$\W_{p,q}(\sll_N)\simeq \W_{q,p}(\sll_N)$$
extends to 
$$\mathbf{L}_{p,q}^{(N)}(\mu,\nu)\simeq \mathbf{L}_{q,p}^{(N)}(\nu,\mu),$$
see \cite{ACF} for details. Now, we can formulate the conjecture on the decomposition in our case.
\begin{conjecture}\hspace{0mm}
\begin{enumerate}[wide, labelindent=0pt]
        \item For $k=-5+{p}/{3}$ $(p\geq7,\ (p,3)=1)$,
        \begin{align*}
            \W_{-5+{p}/{3}}(\sll_5,f_{2,3})\simeq \bigoplus_{a\in \Z_{2(p-5)}}\left(\bigoplus_{\substack{\lambda\in \widehat{\mathrm{P}}_+^{p-5}(2)\\\Proj(\lambda)=a}} \mathbf{L}_{p,p-3}^{(p-5)}(0,\sigma^{\frac{a-\ell(\lambda)}{2}}\lambda^t)\otimes \mathbf{L}_{p-3,3}^{(2)}(\lambda,0) \right)\otimes V_{\frac{5a}{\ssqrt{10(p-5)}}+\ssqrt{10(p-5)}\Z}.
        \end{align*}
        \item For $k=-5+{p}/{2}$ $(p=7,9,\dots)$,
         \begin{align*}
            \W_{-5+{p}/{2}}(\sll_5,f_{1,2^2})\simeq \bigoplus_{a\in \Z_{2(p-5)}}\left(\bigoplus_{\substack{\lambda\in \widehat{\mathrm{P}}_+^{p-5}(3)\\\Proj(\lambda)=a}} \mathbf{L}_{p,p-2}^{(p-5)}(0,\sigma^{\frac{a-\ell(\lambda)}{3}}\lambda^t)\otimes \mathbf{L}_{p-2,2}^{(1,2)}(\lambda,0) \right)\otimes V_{\frac{\sqrt{5}a}{\ssqrt{3(p-5)}}+\ssqrt{15(p-5)}\Z}.
        \end{align*}
    \end{enumerate}
\end{conjecture}


\appendix
\section{OPEs of lower rank \tW-algebras in type \tA}\label{appx1}
We collect all the OPE formulas of the $\W$-algebras of type $A$ with rank up to $4$.
Some are well-known for a long time, e.g. $\W^k(\sll_3,f_3)$ \cite{Za85} and $\W^k(\sll_3,f_{1,2})$ \cite{Ber,Pol} whereas some recently appear in the literature, see \cite{AMP, U0} for $\W^k(\sll_4,f_{2^2})$ and $\W^k(\sll_5,f_{1,2^2})$.
The OPEs here are presented from a relatively new perspective, that is, interpreted as examples of the webs of $\W$-algebras \cite{PR} built up from the truncations of the universal $\W_\infty$-algebra \cite{L}.
Accordingly, we avoid some special levels (mainly where affine vertex subalgebras of level 0 appear). The calculation of the OPEs is all based on constructing explicit strong generators following \cite[Thm 4.1]{KW04} and computing the OPEs among them.

We explain some notation commonly used in the following:
$$\{e,h,f\},\quad  \{J\},\quad \{L_1,L_2\}$$
are the generator of $V^\ell(\sll_2)$ at a certain level $\ell$, a Heisenberg field and mutually-commuting Virasoro fields. These three corresponding vertex algebras mutually commute. 

We take the conformal vector $L$ (which agrees with the one in \cite{KRW03} for the Dynkin grading) given by the sum 
$$L=L_{\mathrm{sug}}+\sum L_i$$
of the Sugawara vector $L_{\mathrm{sug}}$ for the affine vertex subalgebras and the Virasoro fields $L_1,L_2,\dots$. The tables in each case give the set of strong generators and their conformal weights with respect to $L$. 
For example, the $\W$-algebra $\W^k(\sll_3,f_{1,2})$ in \S \ref{2,1} has strong generators $J,v^\pm,L_1$ of conformal weight $1,3/2,2$, respectively. 

The fields $\Omega_{i,3},\dots,\Omega_{i,N}$ ($i=1,2$) have conformal weight $3,\dots,N$ and commute with $e,h,f,J$ and $L_j$ ($j\neq i$) if they exist. 
The fields $L_i,\Omega_{i,3}, \Omega_{i,4},\dots$ satisfy the same OPEs as the generators $L,W_{3}, W_4,\dots$ in the universal $\W_\infty$-algebra $\W_\infty[c,\lambda]$ (see Appendix \ref{Winfinity algebra}) up to normalization and specialization of parameters 
$$c=c_{\alpha+\beta,\beta}(k+a),\quad \lambda=\lambda_{\alpha+\beta,\beta}(k+a).$$
Such specializations appear in the affine cosets of hook-type $\W$-algebras, see \eqref{eq:c_parameter}-\eqref{eq:l_parameter}.
Accordingly, we express their OPEs as
\begin{align}\label{OPEs}
\lbr{W_p}{W_q}=P_{p,q}^{\alpha,\beta}(k+a)
\end{align}
including $L=W_2$ together with necessary data for normalization and truncation 
$$\Omega_{i,p}=n_{i,p} W_{i,p},\quad W_{i,N+M}=Q_{i,M}(W_{i,2},\dots,W_{i,N})\quad M=1,2,\dots$$
for $W_{i,N+M}$ appearing in \eqref{OPEs},
see \S \ref{4} for such an example. 

The fields $v_{i}^\pm$ are assigned for the remaining fields, which give extensions of the conformal embeddings
of the form
$$`` V^\ell(\sll_2)\otimes \pi^J \otimes \langle L_1,\Omega_{1,3},\dots \rangle \otimes \langle L_2,\Omega_{2,3},\dots \rangle\hookrightarrow \W^k(\g,f)".$$
If the sub-indexes are unnecessary, we always drop them off, e.g. $L_1=L$ in $\W^k(\sll_2)$, see \ref{2} below.

\subsection{Rank 1}
\subsubsection{}$\W^k(\sll_2,f_{2})$ at $k\notin \{-2\}$. \label{2}
\begin{fleqn}[\parindent]
\renewcommand{\arraystretch}{1.3}
\begin{align*}
{}\hspace{1cm}
    \begin{array}{|c|c|} \hline
      \Delta &2\\ \hline
      \text{gen.} & L\\ \hline
    \end{array}
    \hspace{2cm} \lbr{L}{L}=P_{2,2}^{2,0}(k).
\end{align*}
\renewcommand{\arraystretch}{1}
\end{fleqn}
\subsection{Rank 2}
\subsubsection{}$\W^k(\sll_3,f_{3})$ at $k\notin \{-3\}$.
\begin{fleqn}[\parindent]
\renewcommand{\arraystretch}{1.3}
\begin{align*}
{}\hspace{1cm}
    \begin{array}{|c|c|c|} \hline
      \Delta &2&3\\ \hline
      \text{gen.} & L & \Omega_3\\ \hline
    \end{array}
    \hspace{2cm} 
    \begin{array}{l}\lbr{W_p}{W_q}=P_{p,q}^{3,0}(k),\\ \displaystyle{W_3=n_k^{-1/2} \Omega_3},\end{array}
\end{align*}
\renewcommand{\arraystretch}{1}
\end{fleqn}
where
\begin{fleqn}[\parindent]
\begin{align*}
{}\hspace{1cm}
    &n_k=-\tfrac{1}{6}(k+3)^2(3k+4)(5k+12),\\
    &W_4=\tfrac{9(k+2)^2}{2(3k+4)(5k+12)} L'' - \tfrac{4(k+3)}{(3k+4)(5k+12)} L^2.
\end{align*}
\end{fleqn}
\subsubsection{}\label{2,1} {$\W^k(\sll_3,f_{1,2})$ at $k\notin\{-3,-\tfrac{3}{2}\}$.
\begin{fleqn}[\parindent]
\renewcommand{\arraystretch}{1.3}
\begin{align*}
{}\hspace{1cm}
    \begin{array}{|c|c|c|c|} \hline
      \Delta &1  & 3/2 & 2\\ \hline
      \text{gen.}        &J  & v^\pm  & L_1\\ \hline
    \end{array}
    \hspace{2cm} 
    \begin{array}{l}\lbr{L_1}{L_1}=P_{2,2}^{2,1}(k)\\ \lbr{J}{J}=\tfrac{1}{3}(2k+3)\lm{1} \end{array}
\end{align*}
\renewcommand{\arraystretch}{1}
\begin{align*}
{}\hspace{1cm}
    &\lbr{L_1}{v^\pm}=\tfrac{3(k+1)}{(3+2k)}v^\pm\lm{1}+v^\pm{}'\mp\tfrac{3}{(3+2k)}Jv^\pm,\qquad \lbr{J}{v^\pm}=\pm v^\pm, \qquad \lbr{v_\pm}{v_\pm}=0,\\
    &\lbr{v^+}{v^-}=(k+1)(2k+3)\lm{2}+3(k+1)J\lm{1}+(-(k+3)L_1+\tfrac{3(k+1)}{2}J'+\tfrac{9(k+1)}{2(3+2k)}J^2)
\end{align*}
\end{fleqn}
\subsection{Rank 3}
\subsubsection{} \label{4}
{$\W^k(\sll_4,f_{4})$ at $k\notin \{-4,-\tfrac{5}{2}, -\tfrac{10}{3}\}$}. 
\begin{fleqn}[\parindent]
\renewcommand{\arraystretch}{1.3}
\begin{align*}
{}\hspace{1cm}
\begin{array}{|c|c|c|c|} \hline
      \Delta &2  & 3 & 4\\ \hline
      \text{gen.}        &L  & \Omega_3  & \Omega_4\\ \hline
    \end{array}
    \hspace{2cm} 
    \begin{array}{l}\lbr{W_p}{W_q}=P_{p,q}^{4,0}(k),\\ 
     W_3=n_k^{-1/2}\Omega_3,\qquad W_4=n_k^{-1}\Omega_4,\end{array}
\end{align*}
\renewcommand{\arraystretch}{1}
\end{fleqn}
where 
\begin{fleqn}[\parindent]
\begin{align*}
{}\hspace{1cm}
    &n_k=-(4+k)^2(5+2k)(10+3k),\\
    &W_5=\tfrac{1}{2n_k^{3/2}(5+2k)(10+3k)}\big(3(86+59k+10k^2)\Omega_3''-16(4+k) L\Omega_3\big),\\
    &W_6=\frac{1}{n_k^2}\bigg(\tfrac{(4+k)^2(86+59k+10k^2)}{4}L''''+\tfrac{3(86+59k+10k^2)}{(5+2k)(10+3k)}\Omega_4''-\tfrac{8(4+k)}{(5+2k)(10+3k)}(L\Omega_4+3\Omega_3^2)\bigg).
\end{align*}
\end{fleqn}
\subsubsection{}{$\W^k(\sll_4,f_{1,3})$ at $k\notin\{-4, -\tfrac{8}{3}\}$}
\begin{fleqn}[\parindent]
\renewcommand{\arraystretch}{1.3}
\begin{align*}
{}\hspace{1cm}
    \begin{array}{|c|c|c|c|} \hline
      \Delta &1  & 2  & 3\\ \hline
      \text{gen.} &J  &L_1, v^\pm & \Omega_{1,3}\\ \hline
    \end{array}
    \hspace{2cm} 
    \begin{array}{l}\lbr{W_p}{W_q}=P_{p,q}^{3,1}(k),\\ \displaystyle{W_3=n_k^{-1/2} \Omega_{1,3}},\end{array}
\hspace{2cm}
\end{align*}
\renewcommand{\arraystretch}{1}
\end{fleqn}
where
\begin{fleqn}[\parindent]
\begin{align*}
{}\hspace{1cm}
    &n_k=-\tfrac{3(2+k)^2(4+k)^2(16+5k)}{2(8+3k)},\\
    &\lbr{J}{J}=\tfrac{1}{4}(3k+8)\lm{1},\qquad \lbr{J}{v^\pm}=\pm v^\pm,\qquad \lbr{v^+}{v^+}=\lbr{v^-}{v^-}=0\\
    &\lbr{L_1}{v^\pm}=\tfrac{2(7+3k)}{(8+3k)}v^\pm\lm{1}\mp\tfrac{4}{(8+3k)}Jv^\pm+v^\pm{}'\\
    &\lbr{\Omega_{1,3}}{v^\pm}=(k+4)\bigg(\pm\tfrac{2(2+k)(7+3k)(16+5k)}{(8+3k)^2}v^\pm\lm{2}+\tfrac{3(2+k)(16+5k)}{2(8+3k)}(\pm v^\pm{}'-\tfrac{4}{(8+3k)}Jv^\pm)\lm{1}\\
    &\hspace{2cm}+(\pm(3+k)v^\pm{}''-\tfrac{8(3+k)}{(8+3k)}Jv^\pm{}'- \tfrac{4(3+k)}{(8+3k)}J'v^\pm \pm\tfrac{16(3+k)}{(8+3k)^2}J^2v^\pm \mp\tfrac{2(4+k)}{(8+3k)}L_1v^\pm)\bigg)\\
    &\lbr{v^+}{v^-}=(2+k)\bigg((5+2k)(8+3k)\lm{3}+4(5+2k)J\lm{2}+(2(5+2k)(J'+\tfrac{4}{(8+3k)}J^2)-(4+k)L_1)\lm{1}\\
    &\hspace{2cm}+\frac{1}{2+k}\Omega_{1,3}-(4+k)(\tfrac{1}{2}L_1'+\tfrac{4}{(8+3k)}J L_1)+ \tfrac{2(5+2k)}{3}(J''+\tfrac{12}{(8+3k)}JJ'+\tfrac{16}{(8+3k)^2}J^3)\bigg).
\end{align*}
\end{fleqn}
where in $P_{3,3}^{3,1}$ the term $W_4=W_{3(1)}W_3$ is explicitly given by 
\begin{fleqn}[\parindent]
\begin{align*}
{}\hspace{1cm}
    W_4=&\tfrac{1}{(2+k)(16+5k)}\bigg( \tfrac{4}{3(2+k)(4+k)}((8+3k)(-2v^+v^-+\Omega_{1,3}')+8J \Omega_{1,3}\\
        &\hspace{1cm}+\tfrac{(8+3k)(2+k)}{2}L_1''-\tfrac{4(4+k)}{3}L_1^2-\tfrac{16}{3}((JL_1)'+\tfrac{4}{(8+3k)}L_1J^2)\\
        &\hspace{1cm}+\tfrac{4(2k+5)}{9(4+k)}((8+3k)J'''+12(J')^2+16JJ''+\tfrac{96}{(8+3k)}J'J^2+\tfrac{64}{(8+3k)^2}J^4)\bigg)
\end{align*}
\end{fleqn}

\subsubsection{}{$\W^k(\sll_4,f_{2,2})$ at $k\notin \{-2,-3,-4\}$} \label{2,2}
\begin{fleqn}[\parindent]
\renewcommand{\arraystretch}{1.3}
\begin{align*}
{}\hspace{1cm}
\begin{array}{|c|c|c|c|} \hline
      \Delta &1  & 2 \\ \hline
      \text{gen.} &J,v_1^\pm  & L_1,L_2, v_2^\pm\\ \hline
    \end{array}
    \hspace{2cm} 
    \begin{array}{l}
    \lbr{L_1}{L_1}=P_{2,2}^{2,2}(k),\quad \lbr{L_2}{L_2}=P_{2,2}^{2,0}(k+1),\\
    \lbr{J}{J}=(2+k)\lm{1},\end{array}
\end{align*}
\renewcommand{\arraystretch}{1}
\end{fleqn}
\begin{fleqn}[\parindent]
\begin{align*}
\hspace{1cm}
    &\lbr{J}{v_i^\pm}=\pm v_i^\pm,\\
    &\lbr{L_1}{v_1^\pm}=\Delta v_1^\pm\lm{1}\mp\tfrac{1}{(3+k)}v_2^\pm,\qquad\Delta:=\tfrac{(3+2k)(8+3k)}{4(2+k)(3+k)},\\
    &\lbr{L_1}{v_2^\pm}=\mp 2(3+k)\Delta v_1^\pm\lm{2}+(\Delta+1)v_2^\pm\lm{1}\pm\tfrac{(3+k)}{(2+k)^2}J^2v_1^\pm\mp\tfrac{2}{(2+k)}Jv_2^\pm\\
    &\hspace{2cm}\mp v_1^\pm L_2-\tfrac{(5+2k)}{(2+k)}Jv_1^\pm{}'-v_1^\pm J'+\tfrac{(5+2k)}{(3+k)}v_2^\pm{}'\pm\tfrac{(3+2k)(11+4k)}{8(3+k)}v_1^\pm{}'',\\
     &\lbr{L_2}{v_1^\pm}=-\tfrac{(3+2k)}{4(3+k)}v_1^\pm\lm{1}\pm\tfrac{1}{(3+k)}v_2^\pm\mp\tfrac{1}{(2+k)}Jv_1^\pm+v_1^\pm{}',\\
    &\lbr{L_2}{v_2^\pm}=-\tfrac{(3+2k)}{4(3+k)}v_2^\pm\lm{1}\mp\tfrac{(3+k)}{(2+k)^2}J^2v_1^\pm\pm\tfrac{1}{(2+k)}Jv_2^\pm\\
    &\hspace{2cm}\pm v_1^\pm L_2+\tfrac{(5+2k)}{(2+k)}Jv_1^\pm{}'+v_1^\pm J'-\tfrac{(2+k)}{(3+k)} v_2^\pm{}'\mp\tfrac{(3+2k)(11+4k)}{8(3+k)}v_1^\pm{}'',\\
    &\lbr{v_1^\pm}{v_1^\pm}=0,\qquad\lbr{v_2^\pm}{v_2^\pm}=-\tfrac{(4+k)(5+2k)(8+3k)}{4(2+k)^2}\left((v_1^\pm)^2+v_1^\pm v_1^\pm{}'\right),\\
    &\lbr{v_1^+}{v_1^-}=2(2+k)\lm{1}+J,\qquad\lbr{v_1^\pm}{v_2^\pm}=\mp\tfrac{(8+3k)}{2(2+k)}(v_1^\pm)^2,\\
    &\lbr{v_1^\pm}{v_2^\mp}=\pm(3+2k)(8+3k)\lm{2}+\tfrac{(3+2k)(8+3k)}{(2+k)}J\lm{1}\pm(2+k)L_2\\
    &\hspace{2cm}\mp(4+k)L_1\pm\tfrac{(7+3k)}{(2+k)}J^2\mp\tfrac{(4+k)}{2(2+k)}v_1^\pm v_1^\mp+(7+3k)J',\\
    &\lbr{v_2^+}{v_2^-}=\tfrac{(3+2k)(8+3k)(36+35k+8k^2)}{2(2+k)}(\lm{3}+J\lm{2})+(16+17k+4k^2)\left(L_2+\tfrac{(7+3k)}{(2+k)^2}(J^2+J')\right)\lm{1}\\
    &\hspace{2cm}-\tfrac{(4+k)}{(2+k)}\left((4+7k+2k^2)L_1+\tfrac{(8+9k+2k^2)}{4(2+k)}v_1^+v_1^-\right)\lm{1}+\tfrac{2(3+k)}{(2+k)}(2J^3+Jv_1^+v_1^-)\\
    &\hspace{2cm}+2(3+k)(2JL_2-\tfrac{(4+k)}{(2+k)}JL_1)-\tfrac{(8+3k)}{2(2+k)}(v_1^+v_2^-+v_1^-v_2^+)\\
    &\hspace{2cm}+\tfrac{(100+157k+77k^2+12k^3)}{(2+k)^2}JJ'+\tfrac{(80+84k+31k^2+4k^3)}{4(2+k)^2}v_1^+v_1^-{}'-\tfrac{(28+18k+3k^2)}{2(2+k)}v_1^-v_1^+{}'\\
    &\hspace{2cm}+\tfrac{1}{2}\left((16+17k+4k^2)L_2'-\tfrac{(4+k)(4+7k+2k^2)}{(2+k)}L_1'\right)+\tfrac{(288+620k+457k^2+142k^3+16k^4)}{4(2+k)^2}J''.
\end{align*}
\end{fleqn}

\subsubsection{}\label{2,1,1} {$\W^k(\sll_4,f_{1,1,2})$ at $k\notin \{-2,-3,-4\}$}.
\begin{fleqn}[\parindent]
\renewcommand{\arraystretch}{1.3}
\begin{align*}
{}\hspace{1cm}
    \begin{array}{|c|c|c|c|} \hline
      \Delta &1  & 3/2 & 2\\ \hline
      \text{gen.} &e,h,f,J  & v_1^{\pm}, v_2^{\pm}  & L_1\\ \hline
    \end{array}
    \hspace{2cm} 
    \begin{array}{l}\lbr{L_1}{L_1}=P_{2,2}^{2,2}(k),\\ 
     \langle e,h,f\rangle \simeq V^{k+1}(\sll_2),\quad \lbr{J}{J}=(2+k)\lm{1}, \end{array}
\end{align*}
\renewcommand{\arraystretch}{1}
\end{fleqn}
\begin{fleqn}[\parindent]
\begin{align}
{}\hspace{1cm}
   \label{eq: example of weyl module} &\Span\{v_1^+,v_2^+\}=\C^2\subset \weyl_{\varpi}^{k+1}\otimes \pi^J_1,\quad \Span\{v_1^-,v_2^-\}=\overline{\C}^2\subset \weyl_{\varpi}^{k+1}\otimes \pi^J_{-1},\\
   \nonumber &\lbr{L_1}{v_1^\pm}=\tfrac{(3+2k)(8+3k)}{4(2+k)(3+k)} v_1^\pm\lm{1}+\bigg(v_1^\pm{}'-\tfrac{1}{2(3+k)}(2e+h)v_1^\pm\mp\tfrac{1}{(2+k)}Jv_1^\pm\bigg),\\
    \nonumber  &\lbr{L_1}{v_2^\pm}=\tfrac{(3+2k)(8+3k)}{4(2+k)(3+k)} v_2^\pm\lm{1}+\bigg(v_2^\pm{}'-\tfrac{1}{2(3+k)}(2f-h)v_2^\pm\mp\tfrac{1}{(2+k)}Jv^\pm\bigg),\\
    \nonumber  &\lbr{v_i^\pm}{v_j^\pm}=0\quad (i,j=1,2),\\
    \nonumber  &\lbr{v_1^+}{v_1^-}=2(2+k)e\lm{1}+(2+k)e'+2Je,\qquad\lbr{v_2^+}{v_2^-}=-2(2+k)f\lm{1}-(2+k)f'-2Jf,\\
    \nonumber &\lbr{v_1^+}{v_2^-}=-2(1+k)(2+k)\lm{2}-((k+2)h+2(k+1)J)\lm{1}\\
    \nonumber &\hspace{2cm}+(4+k)L_1-\tfrac{2+k}{4(3+k)}h^2-hJ+\tfrac{1+k}{2+k}J^2-\tfrac{(2+k)^2}{2(3+k)}h'-(1+k)J'-\tfrac{2+k}{3+k}ef,\\
    \nonumber  &\lbr{v_2^+}{v_1^-}=2(1+k)(2+k)\lm{2}-((k+2)h-2(k+1)J)\lm{1}\\
    \nonumber  &\hspace{2cm}-(4+k)L_1+\tfrac{2+k}{4(3+k)}h^2-hJ-\tfrac{1+k}{2+k}J^2-\tfrac{(2+k)(4+k)}{2(3+k)}h'+(1+k)J'+\tfrac{2+k}{3+k}ef.
\end{align}
\end{fleqn}

\subsection{Rank 4}
\subsubsection{}{$\W^k(\sll_5,f_{5})$ at $k\notin\{-5,-\tfrac{18}{5},-\tfrac{30}{7}\}$}
\begin{fleqn}[\parindent]
\renewcommand{\arraystretch}{1.3}
\begin{align*}
{}\hspace{1cm}
\begin{array}{|c|c|c|c|c|} \hline
      \Delta &2  & 3 & 4 & 5\\ \hline
      \text{gen.}        &L  & \Omega_3  & \Omega_4 & \Omega_5\\ \hline
    \end{array}
    \hspace{2cm} 
    \begin{array}{l}\lbr{W_p}{W_q}=P_{p,q}^{5,0}(k),\\ 
     W_3=n_k^{-1/2}\Omega_3,\qquad W_4=n_k^{-1}\Omega_4,\qquad W_5=n_k^{-1/2}\Omega_5\end{array}
\end{align*}
\renewcommand{\arraystretch}{1}
\end{fleqn}
where 
\begin{fleqn}[\parindent]
\begin{align*}
{}\hspace{1cm}
    &n_k=-\tfrac{3}{10}(5+k)^2(18+5k)(30+7k),\\
    &W_6=\tfrac{(5+k)^2}{5n_k^2}\bigg(-12(5+k)^4 L^3+  62(5+k)L\Omega_4+ 156(5+k)\Omega_3^2+ 15(4+k)^2(5+k)^3 \left( (L')^2+\tfrac{11}{2}L''L\right)\\
    &\hspace{2cm}-\tfrac{3}{4}(6040+3056k+385k^2)\Omega_4'' -\tfrac{(5+k)^2(2736200 + 2784560 k + 1060492 k^2 + 179140 k^3 + 11325 k^4)}{40}L^{(4)}\bigg),\\
    &W_7=\tfrac{(5+k)^2}{5n_k^{5/2}}\bigg(62(5+k)\left(n_k L\Omega_5 + \tfrac{249}{31}\Omega_3\Omega_4\right) - 108(5+k)^4 L^2\Omega_3 + \tfrac{6}{5}(5+k)^3 (4624 + 2264 k + 279 k^2) L'\Omega_3'\\
    &\hspace{2cm}+ \tfrac{3}{10}(5+k)^3 (8908 + 5108 k + 693 k^2) L\Omega_3'' + \tfrac{9}{10}(5+k)^3 (14896 + 7256 k + 891 k^2)L''\Omega_3\\
    &\hspace{2cm}- \tfrac{1}{40}(5+k)^2 (49351620 + 50829408 k + 19561301 k^2 + 3333973 k^3 + 212355 k^4)\Omega_3\\
    &\hspace{2cm}- \tfrac{3}{4}n_k(10768 + 5480 k + 693 k^2) \Omega_5'' \bigg),\\
    &W_8=\tfrac{(5+k)^3}{25n_k^{3}}\bigg(-744(5+k)^6 L^4+ 3304 (5 + k)^3 \left(L^2\Omega_4 + \tfrac{804}{413} L\Omega_3^2\right)+ 3420 n_k \Omega_3\Omega_5 + 2490 \Omega_4^2 \\
    &\hspace{2cm} - 165 (5 + k)^2 (7120 + 3608 k + 455 k^2) L\Omega_4'' + 60 (5 + k)^5 (16400 + 8344 k + 1055 k^2) (L')^2L \\
    &\hspace{2cm}- 96 (5 + k)^2 (16055 + 8179 k + 1035 k^2) L'\Omega_4' - 540 (5 + k)^2 (7520 + 3832 k + 485 k^2) (\Omega_3')^2 \\
    &\hspace{2cm}+ 60 (5 + k)^5 (9440 + 4792 k + 605 k^2) L''L^2 - 6 (5 + k)^2 (99200 + 50728 k + 6435 k^2) L''\Omega_4 \\
    &\hspace{2cm}- \tfrac{1125}{8}(5 + k)^4 (4 + k)^2 (10768 + 5480 k + 693 k^2) (L'')^2 \\
    &\hspace{2cm}- 345 (5 + k)^2 (10708 + 5468 k + 693 k^2) \Omega_3''\Omega_3 \\
    &\hspace{2cm}- \tfrac{16}{5} (5 + k)^4 (14926850 + 15163220 k + 5765902 k^2 + 972700 k^3 + 61425 k^4) L^{(3)}L'\\
    &\hspace{2cm} - \tfrac{3}{4} (5 + k)^4 (37537600 + 38178400 k + 14532804 k^2 + 2453860 k^3 + 155075 k^4) L^{(4)}L \\
    &\hspace{2cm}+ \tfrac{9}{8}(5 + k) (31517800 + 31880080 k + 12076744 k^2 + 2030640 k^3 + 127875 k^4) \Omega_4^{(4)} \\
    &\hspace{2cm}+ \tfrac{(5 + k)^3 }{4800}(3284339216000 + 5030716348800 k + 
    3206045930240 k^2 + 1088119208832 k^3\\
    &\hspace{2cm} + 207431281760 k^4 + 21059173800 k^5 + 889552125 k^6)L^{(6)}\bigg).
\end{align*}
\end{fleqn}
\subsubsection{}{$\W^k(\sll_5,f_{1,4})$ at $k\notin\{-5,-\frac{15}{4}\}$}
\begin{fleqn}[\parindent]
\renewcommand{\arraystretch}{1.3}
\begin{align*}
{}\hspace{1cm}
\begin{array}{|c|c|c|c|c|c|} \hline
      \Delta &1 &2 &5/2 & 3 & 4\\ \hline
      \text{gen.} & J  &L_1 &v^\pm & \Omega_{1,3}  & \Omega_{1,4}\\ \hline
    \end{array}
    \hspace{2cm} 
    \begin{array}{l}\lbr{W_{1,p}}{W_{1,q}}=P_{p,q}^{4,1}(k)\\ 
     W_{1,3}=n_k^{-1/2}\Omega_{1,3},\qquad W_{1,4}=n_k^{-1}\Omega_{1,4}\\
     \lbr{J}{J}=\tfrac{1}{5} (15+4 k)\lm{1}\end{array}
\end{align*}
\renewcommand{\arraystretch}{1}
\end{fleqn}
\begin{align*}
&\lbr{J}{v^\pm}=v^\pm,\quad \lbr{L_1}{v_+}=\tfrac{5(7+2k)}{(15+4k)}v_+\lm{1}+v_+'\mp\tfrac{5}{15+4 k}Jv^\pm,\\
&\lbr{\Omega_{1,3}}{v^\pm}=\tfrac{(5+k)(13+4k)(25+6k)}{(15+4k)}\bigg(\pm \tfrac{5(7+2k)}{3(15+4k)}v^\pm\lm{2}+(\pm\tfrac{1}{2}v^\pm{}'-\tfrac{5}{2(15+4k)}Jv^\pm)\lm{1}\bigg)\\
&\hspace{2cm}\pm((4+k)(5+k)v^\pm{}''+\tfrac{25(4+k)(5+k)}{(15+4k)^2}J^2v^\pm-\tfrac{2(5+k)^2}{(15+4k)}L_1v^\pm) -\tfrac{5(4+k)(5+k)}{(15+4k)} (J'v^\pm+2Jv^\pm{}'),\\
&\lbr{\Omega_{1,4}}{v^\pm}=\tfrac{5(5+k)^2(13+4k)(25+6k)}{(15+4k)^2}\bigg(\tfrac{(7+2k)(115+62k+8k^2)}{2(15+4k)}v^\pm\lm{3}+\tfrac{(545+298k+40k^2)}{2}(v^\pm{}'\mp\tfrac{1}{(15+4k)}Jv^\pm)\lm{2}\bigg)\\
&\hspace{2cm}+\bigg(\tfrac{2(4+k)(5+k)^2(365+205 k+28 k^2)}{(15+4 k)}(v^\pm{}''\mp\tfrac{5}{(15+4 k)}(J'v^\pm+2Jv^\pm{}')+\tfrac{25}{(15+4 k)^2}J^2v^\pm)\\
&\hspace{2cm}-\tfrac{4 (5+k)^3 (40+27 k+4 k^2)}{(15+4k)^2}L_1v^\pm\bigg)\lm{1}
+(5+k)^2\bigg(-\tfrac{12}{15+4 k}\Omega_{1,3}v^\pm+\tfrac{2 (5+k) (5+2 k) (17+4 k)}{(15+4k)^2}L_1'v^\pm\\
&\hspace{2cm}+\tfrac{300 (4+k)^2}{(15+4 k)^2}(J^2v^\pm{}'+J'Jv^\pm)\mp\tfrac{4(4+k)(5+k)}{(15+4k)}(2L_1v^\pm{}'-\tfrac{10}{(15+4k)}JL_1v^\pm)\\
&\hspace{2cm}\mp\tfrac{20 (4+k)^2}{(15+4 k)}(3Jv^\pm{}''+3J'v^\pm{}'+J''v^\pm+\tfrac{25}{(15+4 k)^2}J^3v^\pm)+4 (4+k)^2v^\pm{}^{(3)}\bigg),\\
&\lbr{v^+}{v^-}=(3+k)(7+2k)(11+3k)((15+4k)\lm{4}+5J\lm{3})\\
&\hspace{2cm}+(3+k)(7+2k)\bigg(-(5+k)L_1+\tfrac{5}{2}(11+3k)(J'+\tfrac{5}{(15+4k)}J^2)\bigg)\lm{2}\\
&\\
&\hspace{2cm}+(k+3)\bigg(\Omega_{1,3}+\tfrac{5(7+2k)(11+3k)}{6}(J''+\tfrac{15}{(15+4k)}J'J+\tfrac{25}{(15+4k)^2}J^3)\\
&\hspace{6cm}-\tfrac{(5+k)(7+2k)}{2}(L_1'+\tfrac{10}{(15+4k)}JL_1)\bigg)\lm{1},\\
&\hspace{2cm}-\tfrac{1}{4(5+k)}\Omega_{1,4}+\tfrac{(k+3)}{2} \Omega_{1,3}'+\tfrac{5(3+k)}{(15+4k)}J\Omega_{1,3}\\
&\hspace{2cm}+\tfrac{(5+k)(7+2k)}{2}\bigg(\tfrac{(5+k)}{(15+4k)}L_1^2-\tfrac{3(7+2k)}{8} L_1'')-\tfrac{25(3+k)}{(15+4 k)^2}J^2L_1-\tfrac{5(3+k)}{(15+4k)}(JL_1)'\bigg)\\
&\hspace{2cm}+\tfrac{5(3+k)(7+2k)(11+3k)}{4}\bigg(\tfrac{1}{6}J'''+\tfrac{5}{(15+4k)}(\tfrac{2}{3}J''J+\tfrac{1}{2}(J')^2)+\tfrac{25}{(15+4k)^2}J'J^2+\tfrac{125}{6(15+4k)^3}J^4\bigg),\\
\end{align*}
\begin{align*}
\end{align*}
where 
\begin{align*}
    &n_k=-\tfrac{(5+k)^2(7+2k)(13+4k)(25+6k)}{2(15+4k)}\\
    &W_5=n_k^{-1/2}\bigg(\tfrac{1}{(7+2k)(13+4k)^2(25+6k)^2}\bigg((\tfrac{10(3+k)(11+3k)}{3(5+k)^2}(\tfrac{625}{(15+4k)^2}J^5+250J^2J''+50(15+4k)J'J''\\
    &\hspace{6cm}+25(15+4k)JJ'''+\tfrac{1250}{(15+4k)}J'J^3+375(J')^2J+(15+4k)^2J'''')\\
    &\hspace{2cm}+200L_1^2J-\tfrac{100(3+k)}{(5+k)}(\tfrac{50}{3(15+4k)}L_1J^3+10L_1J'J+\tfrac{2(15+4k)}{3}L_1J''+5L_1'J^2+(15+4k)L_1'J')\\
    &\hspace{2cm}+\tfrac{1}{(5+k)^2(7+2k)}(-80(15+4k)^2v^+v^-+1000(3+k)J^2\Omega_{1,3}+200(3+k)(15+4k)J'\Omega_{1,3}\\
    &\hspace{4cm}+200(3+k)(15+4k)J\Omega_{1,3}'-(15+4k)(10875+9168k+2572 k^2+240 k^3)\Omega_{1,3}'')\\
    &\hspace{2cm}+40(15+4k)L_1'L_1-\tfrac{5}{2}(5+k)(7+2k)(13+4k)^2(25+6k)^2\Omega_{1,4}'-\tfrac{5(15+4k)(745+410k+56k^2)}{6(5+k)}L_1'''\bigg)\\
    &\hspace{2cm}-\tfrac{75(15+4 k)}{(13+4 k)^2(25+6 k)^2(5+k)}JL_1''-\tfrac{5+k}{2(15+4k)}50J\Omega_{1,4}+\tfrac{16(15+4k)(27+8k)}{(5+k)(7+2k)}L_1\Omega_{1,3}\bigg)\\
    &W_6=\tfrac{-1}{(7+2k)}\bigg(\tfrac{10(3+k)(11+3k)}{(13+4k)^2(25+6k)^2}\bigg(\tfrac{15625}{3(15+4k)^3}J^6+\tfrac{12500}{(15+4k)^2 }J'J^4+\tfrac{8750}{3(15+4k)}J''J^3+1250JJ'J''+\tfrac{100}{3}(15+4k)(J'')^2\\
    &\hspace{2cm}+375J'''J^2+50(15+4k)J'''J'+25(15+4k)J''''J+\tfrac{2(15+4k)^2}{3}J^{(5)}+\tfrac{5625}{(15+4k)}J^2(J')^2+250(J')^3\bigg)\\
    &\hspace{2cm}+\tfrac{(5+k)}{(13+4k)^2(25+6k)^2}\bigg(-\tfrac{37500(3+k)}{(15+4k)^2}L_1J^4-80(5+k)L_1^3-6000(3+k)L_1'JJ'+\tfrac{3000 (5+k)}{(15+4 k)}J^2L_1^2\\
    &\hspace{4cm}-\tfrac{30000(3+k)}{(15+4k)}J^2J'L_1-1500(3+k)L_1(J')^2-3000(3+k)L_1J''J-\tfrac{15000(3+k)}{(15+4k)}J^3L_1'\\
    &\hspace{4cm}-100(3+k)(15+4k)L_1J'''-375 (33+10k)J^2L_1''-25(925+518k+72 k^2)JL_1'''\\
    &\hspace{4cm}-300(3+k)(15+4k)J'L_1''\bigg)\\
    &\hspace{2cm}+\tfrac{20(3+k)}{(7+2 k)(13+4 k)^2(25+6 k)^2}\bigg(\tfrac{1000}{(15+4k)}J^3\Omega_{1,3}+300JJ'\Omega_{1,3}+10 (15+4k)J''\Omega_{1,3}+300J^2\Omega_{1,3}'\\
    &\hspace{4cm}+30(15+4k)J\Omega_{1,3}''+(15+4k)^2\Omega_{1,3}'''+30(15+4k)J'\Omega_{1,3}'\bigg)\\
    &\hspace{2cm}-\tfrac{(5+k)^3(7+2k)}{(15+4k)}(2(49+16k)L_1\Omega_{1,4}+75J\Omega_{1,4}')+\tfrac{240(15+4 k)^2}{(7+2 )(13+4k)^2(25+6k)^2}(v^+{}'v^--v^+v^-{}')\\
    &\hspace{2cm}+\tfrac{(5+k)^2(7+2k)(21225+18196 k+5144 k^2+480 k^3)}{4 (15+4k)} \Omega_{1,4}''-\tfrac{48(15+4k)(27+8k)}{(5+k)(7+2k)(13+4k)^2(25+6k)^2}\Omega_{1,3}^2\\
    &\hspace{2cm}-\tfrac{375(5+k)^3(7+2k)}{(15+4k)^2}J^2\Omega_{1,4}+\tfrac{1200(5+k)^2}{(13+4k)^2(25+6k)^2}JL_1'L_1-\tfrac{300(3+k)(5+k)(15+4k)}{(13+4k)^2(25+6 k)^2}J''L_1'\\
    &\hspace{2cm}+\tfrac{10(5+k)(15+4k)}{(13+4 k)^2(25+6 k)^2}(2(19+4k)(L_1')^2+(59+14k)L_1''L_1)\\
    &\hspace{2cm}+\tfrac{(1290375+1971350 k+1179822 k^2+347044 k^3+50304 k^4+2880 k^5)}{6(13+4 k)^2(25+6k)^2} L_1^{(4)}-\tfrac{2400(15+4k)}{(7+2k)(13+4k)^2(25+6 k)^2}Jv^+v^-\bigg)
\end{align*}

\subsubsection{}{$\W^k(\sll_5,f_{2,3})$} at $k\notin\{ -\frac{10}{3},-4,-5\}$. \label{3,2}
\begin{fleqn}[\parindent]
\renewcommand{\arraystretch}{1.3}
\begin{align*}
{}\hspace{1cm}
\begin{array}{|c|c|c|c|c|c|c|c|} \hline
      \Delta &1 &3/2 &2 &5/2  & 3 \\ \hline
      \text{gen.} &J &v_1^\pm &L_1, L_2 & v_2^\pm & \Omega_{1,3}  \\ \hline
    \end{array}
    \hspace{2cm} 
    \begin{array}{l}\lbr{W_{1,p}}{W_{1,q}}=P_{p,q}^{3,2}(k),\quad \lbr{L_2}{L_2}=P_{2,2}^{2,0}(k+2), \\ 
     W_{1,3}=n_k^{-1/2}\Omega_3,\\
     \lbr{J}{J}=\tfrac{2}{5}(10+3k)\lm{1}, \end{array}
\end{align*}
\renewcommand{\arraystretch}{1}
\end{fleqn}
\begin{fleqn}[\parindent]
\begin{align*}
{}\hspace{1cm}
    &\lbr{J}{v_{i}^{\pm}}=\pm v_i^\pm\ (i=1,2),\\
    &\lbr{L_1}{v^\pm_1}=\Delta v^\pm_1\lm{1}\mp\tfrac{1}{k+4}v^\pm_2,\quad \Delta=\tfrac{3(3+k)(15+4k)}{2(4+k)(10+3k)},\\
    &\lbr{L_1}{v^\pm_2}=\mp2(k+4)\Delta v^\pm_1\lm{2}+(\Delta+1)v^\pm_2\lm{1}+2 v_2^\pm{}'-\tfrac{5(k+4)}{2(3k+10)}(J'v^\pm_1+2 Jv^\pm_1{}') \\
    &\hspace{2cm}\mp (L_2v^\pm_1-(k+4)v^\pm_1{}''+\tfrac{5}{3k+10}Jv_2^\pm-\tfrac{25(k+4)}{4(3k+10)^2}J^2v^\pm_1), \\
    &\lbr{L_2}{v^\pm_1}=-\tfrac{5+2k}{4(k+4)}v^\pm_1\lm{1}+v^\pm_1{}'\mp(\tfrac{5}{2(3k+10)}Jv^\pm_1-\tfrac{1}{k+4}v^\pm_2),\\
    &\lbr{L_2}{v^\pm_2}=-\tfrac{5+2k}{4(k+4)}v^\pm_2\lm{1}-v^\pm_2{}'+\tfrac{5(k+4)}{2(3k+10)}(J'v^\pm_1+2Jv^\pm_1{}')\\
    &\hspace{2cm}\pm(L_2v^\pm_1-(k+4)v^\pm_1{}''+\tfrac{5}{2(3k+10)}Jv_2^\pm-\tfrac{25(k+4)}{4(3k+10)^2}J^2v^\pm_1),\\
    &\\
    &\lbr{\Omega_{1,3}}{v^\pm_1}=\mp\tfrac{5(3+k)(5+k)(8+3k)(15+4k)}{2(10+3k)^2}v^\pm_1\lm{2}\mp\tfrac{5(5+k)(8+3k)}{2(10+3k)}v^\pm_2\lm{1}\pm\tfrac{2(5+k)^2}{(10+3k)}L_1v_1^\pm\\
    &\hspace{2cm}+(5+k)\bigg(2v^\pm_2{}'-\tfrac{5(4+k)}{2(10+3k)}J'v^\pm_1 -\tfrac{5(4+k)}{(10+3k)}Jv^\pm_1{}'\\
    &\hspace{4cm}\pm( (4+k)v^\pm_1{}'' - L_2v^\pm_1 + \tfrac{25(4+k)}{4(10+3k)^2}J^2v^\pm_1 - \tfrac{5}{(10+3k)}Jv^\pm_2)\bigg),\\
    &\lbr{\Omega _{1,3}}{v^\pm_2}=\tfrac{5(k+5)(3k+8)}{2(3k+10)^2}(3(k+3)(k+4)(4 k+15)v^\pm_1\lm{3}\mp(10 k^2+71k+125)v^\pm_2\lm{2})\\
    &\hspace{2cm}-\tfrac{3(k+4)(k+5)}{(3k+10)}\bigg(\mp\tfrac{5(k+4)(7k+20)}{4(3 k+10)}(J'v^\pm_1+2 Jv^\pm_1{}')-\tfrac{5(7k+20)}{2(3 k+10)}Jv^\pm_2+\tfrac{25(k+4)(7k+20)}{8(3k+10)^2}J^2v^\pm_1\\
    &\hspace{5cm}+\tfrac{1}{2}(7k+20)((k+4) v^\pm_1{}''-L_2v^\pm_1\pm2 v^\pm_2{}')+\tfrac{2}{3}(k+5)L_1v^+_1\bigg)\lm{1}\\
    &\hspace{2cm}\pm\tfrac{5(k+5)(k+4)^2}{(3 k+10)}(J''v^\pm_1 +J'v^\pm_1{}')-\tfrac{75(k+5)(k+4)^2}{2(3 k+10)^2}(JJ'v^\pm_1+J^2v^\pm_1{}'\mp\tfrac{5}{2}Jv^\pm_1{}'')\\
    &\hspace{2cm}+\tfrac{15(k+5)(k+4)}{2(3k+10)}(J'v^\pm_2+2Jv^\pm_2{}'\mp\tfrac{2}{3}L_2Jv^\pm_1)\pm\tfrac{2(k+5)^2}{3 k+10}L_1v^\pm_2-2 (k+5) (k+4)^2 v^\pm_1{}'''\\
    &\hspace{2cm}\pm\tfrac{125 (k+5) (k+4)^2}{4 (3 k+10)^3}J^3v^\pm_1{}-\tfrac{2 (k+5)^2 (k+4)}{3 k+10}L_1'v^\pm_1\mp(k+5)L_2v^\pm_2\\
    &\hspace{2cm}\mp\tfrac{75 (k+5) (k+4)}{4 (3 k+10)^2}J^2v^\pm_2+2(k+5)(k+4)(L_1'v^\pm_1+L_2v^\pm_1{}'\mp\tfrac{3}{2}v^\pm_2{}''),\\
    &\lbr{v^\pm_1}{v^\pm_1}=0,\quad \lbr{v^\pm_1}{v^\pm_2}=\mp\tfrac{15+4k}{10+3k} (v^\pm_1)^2,\quad \lbr{v^\pm_2}{v^\pm_2}=-\tfrac{(k+5)(3k+11)(4k+15)}{(3k+10)^2}((v^\pm_1)^2\lm{1}+v^\pm_1{}'v^\pm_1),\\
    &\lbr{v^+_1}{v^-_1}=2(k+3)(3k+10)\lm{2}+5(3+k)J\lm{1}+\tfrac{25(k+3)}{4(3k+10)}J^2+\tfrac{5}{2}(k+3)J'-(k+5) L_1+(k+3)L_2,\\
    &\lbr{v^+_1}{v^-_2}=6(k+3)^2(4k+15)\lm{3}+\tfrac{15(k+3)^2(4k+15)}{3k+10}J\lm{2}\\
    &\hspace{2cm}+(-\tfrac{(k+5)(2k+5)(3k+11)}{(3 k+10)}L_1+\tfrac{3(k+3)^2(4k+15)}{4(3k+10)}(4L_2+10J'+\tfrac{25}{(3k+10)}J^2))\lm{1}\\
    &\hspace{2cm}+\Omega_{1,3}-\tfrac{(k+5)}{(3k+10)}v^+_1v^-_1-\tfrac{1}{2}(k+5)(3k+10)L_1'+(k+3)(2k+7)L_2'\\
    &\hspace{2cm}-\tfrac{5(k+3)}{(3k+10)}\bigg((k+5)L_1-(2k+7)L_2)J-(2k+7)(\tfrac{(3k+10)}{3}J''+\tfrac{5}{2}J'J+\tfrac{25}{12(3k+10)}J^3)\bigg),\\
    &\lbr{v^+_2}{v^-_1}=6(k+3)^2(4k+15)\lm{3}+\tfrac{15(k+3)^2(4k+15)}{(3k+10)}J\lm{2}\\
    &\hspace{2cm}+(-\tfrac{(k+5)(2k+5)(3k+11)}{(3k+10)}L_1+\tfrac{3(4k+15)(k+3)^2}{4(3k+10)}(4L_2+10J'+\tfrac{25}{(3k+10)}J^2))\lm{1}\\
    &\hspace{2cm}-\Omega_{1,3}-\tfrac{(k+5)}{(3k+10)}v^+_1v^-_1-\tfrac{1}{2}(k+2)(k+5)L_1'+(k+3)(2k+7)L_2'\\
    &\hspace{2cm}-\tfrac{5(k+3)}{(3k+10)}\bigg(((k+5)L_1-(2k+7)L_2)J-(2k+7)(\tfrac{(3k+10)}{3}J''+\tfrac{5}{2}J'J+\tfrac{25}{12(3k+10)}J^3)\bigg),\\
    &\lbr{v^+_2}{v^-_2}=\tfrac{6(k+3)^2(4k+15)(15k^2+103k+175)}{(3k+10)}(\lm{4}+\tfrac{5}{2(3k+10)}J\lm{3})\\
    &\hspace{2cm}+\bigg(\tfrac{15(4k+15)(15k^2+103k+175)(k+3)^2}{4(3k+10)^2}(\tfrac{5}{(3k+10)}J^2+2J'+\tfrac{4}{5}L_2)\\
    &\hspace{6cm}-\tfrac{(k+5)(2k+5)(3k+11)(9k^2+59k+95)}{(3k+10)^2}L_1\bigg)\lm{2}\\
    &\hspace{2cm}+\bigg(-\tfrac{2 (3k^2+18k+25)}{(3k+10)}\Omega _{1,3}-\tfrac{(k+5)(3k^2+21k+35)}{(3k+10)^2} v^+_1v^-_1\\
    &\hspace{3cm}-\tfrac{(k+3)(k+5)(9k^2+58k+90)}{(3k+10)}(L_1'+\tfrac{5}{(3k+10)}J L_1)\\
    &\hspace{3cm}+\tfrac{(k+3)(2k+7)(15k^2+102k+170)}{(3k+10)}(L_2'+\tfrac{5}{(3k+10)}J L_2)\\
    &\hspace{3cm}+\tfrac{5(k+3)(2k+7)(15k^2+102 k+170)}{(3k+10)}(J''+\tfrac{5}{2(3k+10)}J J'+\tfrac{25}{12(3k+10)^2}J^3)\bigg)\lm{1}\\
    &\hspace{2cm}-\tfrac{(3k^2+18k+25)}{(3k+10)} \Omega _{1,3}' +\tfrac{5(k+4)}{(3k+10)}J( v^+_1 v^-_1- \Omega _{1,3})-\tfrac{(4k+15)}{(3k+10)}(v^+_2 v^-_1+v^-_2v^+_1)\\
    &\hspace{2cm}+\tfrac{(9k^3+100k^2+375k+475)}{(3k+10)^2} v^+_1 v^-_1{}'-\tfrac{(4k^2+32 k+65)}{(3k+10)} v^+_1{}'v^-_1\\
    &\hspace{2cm}+(k+3)(k+4)L_1^2-(k+4)(k+5)L_1L_2\\
    &\hspace{2cm}-\tfrac{5(k+3)(k+5)(9k^2+58k+90)}{2(3k+10)^2}J' L_1-\tfrac{5(k+5)(9k^3+82k^2+242k+230)}{2(3 k+10)^2}J L_1'\\
    &\hspace{2cm}-\tfrac{25(k+4)(k+5)(3k+8)}{4(3k+10)^2} L_1 J^2-\tfrac{(k+5)(9 k^4+102 k^3+395 k^2+555 k+125)}{2(3k+10)^2}L_1''\\
    &\hspace{2cm}+\tfrac{5(k+3)(2k+7)(15k^2+102k+170)}{2(3k+10)^2} J'L_2+\tfrac{5(k+3)(15k^3+153k^2+516k+575)}{(3k+10)^2}J L_2'\\
    &\hspace{2cm}+\tfrac{25(k+3)(k+4)(5k+17)}{2(3k+10)^2}L_2 J^2+\tfrac{(k+3)(54k^4+756 k^3+3956 k^2+9165 k+7925)}{2(3k+10)^2}L_2''\\
    &\hspace{2cm}+\tfrac{625(k+3)(k+4)(5k+16)}{48(3k+10)^3}J^4+\tfrac{125(k+3)(30k^3+303k^2+1010k+1110)}{8(3k+10)^3}J'J^2
    &\hspace{2cm}+\tfrac{25(k+3)(180k^4+2523k^3+13220k^2+30675k+26575)}{12(3k+10)^3}J J''\\
    &\hspace{2cm}+\tfrac{25(k+3)(45k^4+633 k^3+3328 k^2+7745k+6725)}{4(3 k+10)^3}(J')^2\\
    &\hspace{2cm}+\tfrac{5(k+3)(90k^4+1218 k^3+6152 k^2+13735 k+11425)}{12(3k+10)^2}J''',
\end{align*}
\end{fleqn}
\begin{align*}
\end{align*}
where 
\begin{align*}
    &n_k=\tfrac{5(2+k)(4+k)(5+k)^2(8+3k)}{2(10+3k)},\\
    &W_4=\tfrac{1}{n_k}(4(5+k)(v_1^+v_2^- +v_2^+v_1^--\Omega_{1,3}')-\tfrac{20(4+k)(5+k)}{(10+3k)}Jv_1^+v_1^-- \tfrac{20(5+k)}{(10+3k)}J\Omega_{1,3}  \\
    &+\tfrac{2(2+k)(5+k)^3}{(10+3k)}L_1^2+\tfrac{25(2+k)(5+k)^2}{(10+3k)^2}L_1J^2-2(3+k)(5+k)L_2^2-\tfrac{25(3+k)(5+k)}{(10+3k)}L_2J^2   \\
    &+ \tfrac{4(5+k)(50+24k+3k^2)}{(10+3k)}v^+_1{}'v^-_1 -4(4+k)(5+k)v^+_1v^-_1{}'+4(5+k)^2L_1L_2\\
    &- \tfrac{3}{4}(2+k)^2(5+k)^2L_1''-(3+k)(5+k)(13+4k)L_2'' \\
    &+\tfrac{10(5+k)^2}{(10+3k)}((2+k)JL_1'+2(3+k)J'L_1 )-10(3+k)(5+k)(JL_2'+\tfrac{2(7+2k)}{(10+3k)}J'L_2  ) \\
    &-25(3+k)(5+k)(\tfrac{1}{15}(10+3k)J'''+\tfrac{25(2+k)}{24(10+3k)^3}J^4+\tfrac{(16+5k)}{3(10+3k)}JJ'' +\tfrac{1}{2} (J')^2+\tfrac{5(3+k)}{(10+3k)^2}J^2J')).
\end{align*}

\subsubsection{}{$\W^k(\sll_5,f_{1,1,3})$ at $k\notin\{-\tfrac{8}{3},-2,-4,-5\}$.}\label{3,1,1}
\begin{fleqn}[\parindent]
\renewcommand{\arraystretch}{1.3}
\begin{align*}
{}\hspace{1cm}
\begin{array}{|c|c|c|c|} \hline
      \Delta &1 &2  & 3 \\ \hline
      \text{gen.} &e,h,f,J &L_1, v_1^\pm, v_2^\pm  & \Omega_{1,3} \\ \hline
    \end{array}
    \hspace{2cm} 
    \begin{array}{l}\lbr{W_p}{W_q}=P_{p,q}^{3,2}(k),\\ 
     W_3=n_k^{-1/2}\Omega_{1,3},\qquad W_4=n_k^{-1}\Omega_4,\qquad W_5=n_k^{-3/2}\Omega_5\end{array}
\end{align*}
\renewcommand{\arraystretch}{1}
\end{fleqn}
where $n_k=-\tfrac{5(2+k)(4+k)(5+k)^2(8+3k)}{2(10+3k)}$ and

\begin{fleqn}[\parindent] 
\begin{align*}
{}\hspace{1cm}
    &\langle e,h,f\rangle \simeq V^{k+2}(\sll_2),\qquad\lbr{J}{J}=\tfrac{2}{5}(10+3k)\lm{1},\\
    &\Span\{v_1^+,v_2^+\}=\C^2\subset \weyl_{\varpi}^{k+2}\otimes \pi^J_1,\qquad \Span\{v_1^-,v_2^-\}=\overline{\C}^2\subset \weyl_{\varpi}^{k+2}\otimes \pi^J_{-1},\\
    &\lbr{L_1}{v^\pm_1}=\Delta v^\pm_1\lm{1}+v^\pm_1{}'-\tfrac{1}{2(4+k)}h v^\pm_1-\tfrac{1}{(4+k)}e v^\pm_2\mp\tfrac{5}{2(10+3k)}J v^\pm_1,\quad \Delta=\tfrac{3(3+k)(15+4k)}{2(4+k)(10+3k)},\\
    &\lbr{L_1}{v^\pm_2}=\Delta v^\pm_2\lm{1}+v_2^\pm{}'+\tfrac{1}{2(4+k)}h v^\pm_2-\tfrac{1}{(4+k)}f v^\pm_1\mp\tfrac{5}{2(10+3k)}J v^\pm_2,\\
    &\lbr{\Omega_{1,3}}{v_1^\pm}=\tfrac{5(5+k)(8+3k)}{2(10+3k)}\bigg(\mp\tfrac{(3+k)(15+4 k)}{(10+3k)}v_{1}^\pm\lm{2}\mp(4+k) v_{1}^\pm{}'\lm{1}+\tfrac{5(4+k)}{2(10+3k)} Jv_{1}^\pm\lm{1}\\
    &\hspace{2cm}\pm(\tfrac{1}{2}hv_{1}^\pm + ev_{2}^\pm)\lm{1}\mp\tfrac{2(4+k)(10+3k)}{5(8+3k)} v_{1}^\pm{}''\pm\tfrac{(10+3k)}{5(8+3k)}\left(h'v_{1}^\pm+2hv_1^\pm{}'+4ev_2^\pm{}'\right)\\
    &\hspace{2cm}+\tfrac{(4+k)}{(8+3 k)}\left(J'v_1^\pm+2Jv_1^\pm{}'\right)\pm\tfrac{4(5+k)}{5(8+3k)}Lv_1^\pm\mp\tfrac{(10+3 k)}{10 (4+k)(8+3 k)}\left(h^2v_1^\pm+4efv_1^\pm\right)\\
    &\hspace{2cm}\pm\tfrac{2(3+k)(10+3 k)}{5 (4+k)(8+3 k)}e'v_2^\pm-\tfrac{1}{(8+3 k)}\left(hJv_1^\pm+2Jev_2^\pm\right)\mp\tfrac{5(4+k)}{2(8+3k)(10+3k)}J^2v_1^\pm\bigg),\\
    &\lbr{\Omega_{1,3}}{v_2^\pm}=\tfrac{5(5+k)(8+3k)}{2(10+3k)}\bigg(\mp\tfrac{(3+k)(15+4k)}{(10+3k)}v_2^\pm\lm{2}\mp(4+k)v_2^\pm{}'\lm{1}+\tfrac{5 (4+k)}{2(10+3k)}Jv_2^\pm\lm{1}\\
    &\hspace{2cm}\mp(\tfrac{1}{2}hv_2^\pm- fv_1^\pm)\lm{1}\mp\tfrac{2(4+k)(10+3k)}{5(8+3k)} v_2^\pm{}''\mp\tfrac{(10+3k)}{5(8+3k)}\left(\tfrac{(2+k)}{(4+k)}h'v_2^\pm+2hv_2^\pm{}'-4fv_1^\pm{}'\right)\\
    &\hspace{2cm}+\tfrac{(4+k)}{(8+3 k)}\left(J'v_2^\pm+2Jv_2^\pm{}'\right)\pm\tfrac{4(5+k)}{5(8+3k)}Lv_2^\pm\mp\tfrac{(10+3 k)}{10 (4+k)(8+3 k)}\left(h^2v_2^\pm+4efv_2^\pm\right)\\
    &\hspace{2cm}\pm\tfrac{2(3+k)(10+3 k)}{5 (4+k)(8+3 k)}f'v_1^\pm+\tfrac{1}{(8+3 k)}\left(hJv_2^\pm-2Jfv_1^\pm\right)\mp\tfrac{5 (4+k)}{2(8+3k)(10+3k)}J^2v_2^\pm\bigg),\\
    &\lbr{v_1^\pm}{v_1^\pm}=\lbr{v_2^\pm}{v_2^\pm}=\lbr{v_1^\pm}{v_2^\pm}=0\\
    &\lbr{v_1^+}{v_1^-}=2(3+k)(10+3k)e\lm{2} + (3+k)(5Je+(10+3k)e')\lm{1} +\tfrac{1}{2}(3+k)(7+2k)e''\\
    &\hspace{2cm} + \tfrac{(3+k)}{2(4+k)}\left(-(6+k)he'+(5+k)h'e\right) + \tfrac{5}{2}(3+k)(Je'+J'e) - (5+k)L_1e\\
    &\hspace{2cm} + \tfrac{(3+k)}{2(4+k)}\left(\tfrac{1}{2}h^2e+2efe\right) + \tfrac{25(3+k)}{4(10+3k)}J^2e,\\
    &\lbr{v_2^+}{v_2^-}=-2(3+k)(10+3k)f\lm{2} -(3+k)(5Jf+(10+3k)f')\lm{1} -\tfrac{1}{2}(3+k)(7+2k)f''\\
    &\hspace{2cm} + \tfrac{(3+k)}{2(4+k)}\left(-(6+k)hf'+(7+k)h'e\right) - \tfrac{5}{2}(3+k)(Jf'+J'f) + (5+k)L_1f\\
    &\hspace{2cm} - \tfrac{(3+k)}{2(4+k)}\left(\tfrac{1}{2}h^2f+2eff\right) - \tfrac{25(3+k)}{4(10+3k)}J^2f,\\
    &\lbr{v_1^+}{v_2^-}=-2(2+k)(3+k)(10+3k)\lm{3}-(3+k)\left(5(2+k)J+(10+3k)h\right)\lm{2}\\
    &\hspace{2cm} + (2+k)(5+k)L\lm{1} - \tfrac{(3+k)(10+3k)}{2(4+k)}\left((3+k)h' + \tfrac{1}{2}h^2 + 2ef\right)\lm{1}\\
    &\hspace{2cm} - \tfrac{5}{2}(3+k)\left(hJ + (2+k)J' + \tfrac{5(2+k)}{2(10+3k)}J^2\right)\lm{1} + \Omega_{1,3} + \tfrac{1}{2}(2+k)(5+k)L'\\
    &\hspace{2cm} - \tfrac{(3+k)^3}{2(4+k)}h'' - \tfrac{5}{6}(2+k)(3+k)J'' + \tfrac{(5+k)}{2}Lh + \tfrac{5(2+k)(5+k)}{2(10+3k)}LJ\\
    &\hspace{2cm} - \tfrac{(3+k)^2}{4(4+k)}\left(3 h'h + 5 h'J + 4ef'\right) - \tfrac{5}{4}(3+k) J'h - \tfrac{25(2+k)(3+k)}{4(10+3k)}J'J - \tfrac{(3+k)(7+2k)}{(4+k)} e'f \\
    &\hspace{2cm} - \tfrac{25(3+k)}{8(10+3k)}hJ^2 - \tfrac{125(2+k)(3+k)}{24(10+3k)^2}J^3 - \tfrac{(3+k)}{8(4+k)}\left(h^3 + 5h^2J + 4hef + 20Jef\right),\\
    &\lbr{v_2^+}{v_1^-}=2(2+k)(3+k)(10+3k)\lm{3}+(3+k)\left(5(2+k)J-(10+3k)h\right)\lm{2}\\
    &\hspace{2cm} - (2+k)(5+k)L\lm{1} + \tfrac{(3+k)(10+3k)}{2(4+k)}\left(-(5+k)h' + \tfrac{1}{2}h^2 + 2ef\right)\lm{1}\\
    &\hspace{2cm} + \tfrac{5}{2}(3+k)\left(-hJ + (2+k)J' + \tfrac{5(2+k)}{2(10+3k)}J^2\right)\lm{1} - \Omega_{1,3} - \tfrac{1}{2}(2+k)(5+k)L'\\
    &\hspace{2cm} - \tfrac{(19+9k+k^2)(3+k)}{2(4+k)}h'' + \tfrac{5}{6}(2+k)(3+k)J'' + \tfrac{(5+k)}{2}Lh - \tfrac{5(2+k)(5+k)}{2(10+3k)}LJ\\
    &\hspace{2cm} + \tfrac{(3+k)}{4(4+k)}\left((11+3k)h'h - 5(5+k)h'J + 4(3+k)e'f\right) - \tfrac{5}{4}(3+k) J'h + \tfrac{25(2+k)(3+k)}{4(10+3k)}J'J\\
    &\hspace{2cm} + \tfrac{(3+k)(7+2k)}{(4+k)} ef' - \tfrac{25(3+k)}{8(10+3k)}hJ^2 + \tfrac{125(2+k)(3+k)}{24(10+3k)^2}J^3 - \tfrac{(3+k)}{8(4+k)}\left(h^3 - 5h^2J + 4hef - 20Jef\right).
\end{align*}
\end{fleqn}

\subsubsection{}{$\W^k(\sll_5,f_{1,2,2})$ at $k\notin \{-4,-5,-\frac{5}{2}\}$}\label{2,2,1}
\begin{fleqn}[\parindent]
\renewcommand{\arraystretch}{1.3}
\begin{align*}
{}\hspace{1cm}
\begin{array}{|c|c|c|c|c|} \hline
      \Delta &  1& 3/2 & 2\\ \hline
      \text{gen.} & e_*,h,f_*,J & v^\pm_1, v^\pm_2  & L_1,L_2, E,F \\ \hline
    \end{array}
    \hspace{2cm} 
    \begin{array}{l}\lbr{L_1}{L_1}=P_{2,2}^{2,3}(k),\quad \lbr{L_2}{L_2}=P_{2,2}^{2,1}(k+1)\\ \lbr{h}{h}=2(2k+5)\lm{1},\quad \lbr{J}{J}=\tfrac{2}{5}(2k+5)\lm{1}
    \end{array}
\end{align*}
\renewcommand{\arraystretch}{1}
\end{fleqn}
\begin{fleqn}[\parindent]
\begin{align*}
{}\hspace{1cm}
    &\langle e_*,h,f_*\rangle \simeq V^{(2k+5)}(\sll_2),\quad \Span\{E,H,F\}=\sll_2\subset \weyl_{2\varpi}^{2k+5}\otimes \pi^J,\\
    &\Span\{v_1^+,v_2^+\}=\C^2\subset \weyl_{\varpi}^{2k+5}\otimes \pi^J_1,\quad \Span\{v_1^-,v_2^-\}=\overline{\C}^2\subset \weyl_{\varpi}^{2k+5}\otimes \pi^J_{-1},\\
    &\lbr{L_1}{e_*}=\tfrac{1}{(4+k)}(\tfrac{(2+k) (10+3 k)}{ (5+2 k)} e_*\lm{1}+\tfrac{(10+3 k)}{2 }e_*'-\tfrac{(10+3 k)}{2(5+2 k)}he_*-E),\\
    &\lbr{L_1}{f_*}=\tfrac{1}{(4+k)}(\tfrac{(2+k) (10+3 k)}{ (5+2 k)} f_*\lm{1}+\tfrac{(10+3 k)}{2}f_*'+\tfrac{(10+3 k)}{2
    (5+2 k)}hf_*+F),\\
    & \lbr{L_2}{e_*}=\tfrac{1}{(4+k)}(-\tfrac{(2+k)^2}{ (5+2 k)}e_*\lm{1}-\tfrac{(2+k)}{2 } e_*'+\tfrac{(2+k)}{2(5+2
     k)}he_*+E),\\
    &\lbr{L_2}{f_*}=\tfrac{1}{(4+k)}(-\tfrac{(2+k)^2}{(5+2 k)}f_*\lm{1}- \tfrac{(2+k)}{2} f_*'-\tfrac{2+k}{2 (5+2 k)}hf_*-F),\\
    &\lbr{L_1}{v^+_2}=\lbr{L_1}{v^-_1}=0,\\
    &\lbr{L_1}{v^+_1}=\tfrac{1}{(4+k)}(\tfrac{(2+k) (10+3 k)}{(5+2 k)} v^+_1\lm{1}+(3+k) v^+_1{}'-\tfrac{5+k}{2(5+2 k)}h v^+_1-\tfrac{5(3+k)}{2(5+2 k)}J v^+_1+e_* v^+_2),\\
    &\lbr{L_1}{v^-_2}=\tfrac{1}{(4+k)}(\tfrac{(2+k) (10+3 k)}{(5+2 k)} v^-_2\lm{1}+(3+k) v^-_2{}'+\tfrac{(5+k)}{2(5+2 k)}
    hv^-_2+\tfrac{5(3+k)}{2(5+2 k)}Jv^-_2+f_*v^-_1),\\
    &\lbr{L_2}{v^+_1}=\tfrac{1}{k+4}(\tfrac{2 (2+k)}{(5+2 k)} v^+_1\lm{1}+ v^+_1{}'+\tfrac{1}{2(5+2k)}hv^+_1-\tfrac{5}{2(5+2 k)}J v^+_1-e_*v^+_2),\\
    &\lbr{L_2}{v^-_2}=\tfrac{1}{(4+k)}(\tfrac{2 (2+k)}{(5+2 k)} v^-_2\lm{1}+ v^-_2{}'-\tfrac{1}{2(5+2 k)}hv^-_2+\tfrac{5}{2 (5+2 k)}Jv^-_2-f_*v^-_1),\\
    &\lbr{L_2}{v^+_2}=\tfrac{1}{2(5+2 k)}(6 (2+k) v^+_2\lm{1}+hv^+_2-5Jv^+_2)+ v^+_2{}',\\
    &\lbr{L_2}{v^-_1}=\tfrac{1}{2(5+2 k)}(6 (2+k) v^-_1\lm{1}-hv^-_1+5Jv^-_1)+ v^-_1{}',\\
    &\lbr{L_1}{E}=\tfrac{(2+k)^2 (10+3 k)}{(4+k) (5+2 k)}e_*\lm{2}+\tfrac{(2+k)}{4(4+k)}((10+3 k)e_*'+\tfrac{2(15+4k)}{(5+2 k)}E-\tfrac{(10+3 k)}{(5+2 k)}he_*)\lm{1}\\
    &\hspace{2cm}+\tfrac{1}{2}E'-L_2e_*+\tfrac{1}{4+k}v^+_1v^-_1+\tfrac{(2+k)^2}{4 (4+k)}e_*''\\
    &\hspace{2cm}+\tfrac{1}{4(4+k) (5+2 k)}\bigg(-(2+k)((3+k) h+5 J)'e_*-(2+k)( (5+2 k) h+5 J)e_*'\\
    &\hspace{6cm}-2 ((5+k)h+5J)E+\tfrac{(2+k)}{2(5+2 k)} ((5+2 k) h^2e_*-25J^2e_*)\bigg),\\
    &\lbr{L_2}{E}=-\tfrac{(2+k)^2 (10+3 k)}{(4+k) (5+2 k)}e_*\lm{2}+(\tfrac{(2+k) (10+3k)}{4(4+k)}(-e_*'+\tfrac{1}{(5+2k)}he_*)+\tfrac{42+27k+4 k^2}{2 (4+k) (5+2 k)}E)\lm{1}\\
    &\hspace{2cm}+\tfrac{1}{2}E'+ L_2e_*-\tfrac{1}{4+k} v^+_1v^-_1-\tfrac{(2+k)^2}{4 (4+k)}e_*''\\
    &\hspace{2cm}+\tfrac{1}{4(5+2 k)(4+k)}\bigg((2+k)((3+k) h+5 J)'e_*+(2+k) ((5+2 k) h+5 J)e_*'\\
    &\hspace{6cm}-2 ((3+k)h-5 J)E-\tfrac{(2+k)}{2(5+2 k)}((5+2 k)h^2e_*-25J^2e_*)\bigg),\\
    &\lbr{L_1}{F}=-\tfrac{(2+k)^2 (10+3 k)}{(4+k) (5+2 k)}f_*\lm{2}+\tfrac{(2+k)}{4(4+k)}(-(10+3 k)f_*'+\tfrac{2(15+4 k)}{(5+2 k)}F-\tfrac{(10+3 k)}{(5+2 k)}hf_*)\lm{1}\\
    &\hspace{2cm}+\tfrac{6+k}{2(4+k)}F'+L_2f_*+\tfrac{1}{4+k}v^+_2v^-_2-\tfrac{(2+k)^2}{4(4+k)}f_*''\\
    &\hspace{2cm}+\tfrac{1}{4 (4+k) (5+2 k)}\bigg(2((5+k) h+5J)F-\tfrac{2+k}{2}(h^2f_*-\tfrac{25 }{(5+2 k)}J^2F)\\
    &\hspace{6cm}-(2+k)(((3+k)h-5 J)'f_*+((5+2 k)h-5 J)f_*')\bigg),\\
    &\lbr{L_2}{F}=\tfrac{(2+k)^2 (10+3 k)}{(4+k) (5+2 k)}f_*\lm{2}+(\tfrac{(2+k) (10+3 k)}{4 (4+k)}(f_*'+\tfrac{1}{(5+2k)}hf)+\tfrac{42+27k+4 k^2}{2 (4+k) (5+2 k)}F)\lm{1}\\
    &\hspace{2cm}+\tfrac{2+k}{2 (4+k)}F'-L_2f_*+\tfrac{(2+k)^2}{4 (4+k)}f_*''-\tfrac{1}{4+k}v^+_2 v^-_2\\
    &\hspace{2cm}+\tfrac{1}{4 (4+k) (5+2 k)}\bigg(2((3+k)h-5J)F+\tfrac{(2+k)}{2(5+2 k)}((5+2 k)h^2f_*-25J^2f_*)\\
    &\hspace{6cm}+(2+k)(((3+k) h-5 J)'f_*+((5+2 k)h-5J)f_*'))\bigg),\\
    &\\
    &\lbr{v^+_i}{v^+_j}=0,\quad \lbr{v^-_i}{v^-_j}=0\quad (i,j=1,2),\\
    &\lbr{v^+_1}{v^-_1}=(k+2)e\lm{1}+E+\tfrac{2+k}{2}e_*'+\tfrac{5(2+k)}{2(5+2 k)}Je_*,\\
    &\lbr{v^+_2}{v^-_2}=-(k+2)f_*\lm{1}-F-\tfrac{(k+2)}{2}f_*'-\tfrac{5(k+2)}{2(5+2k)}Jf_*,\\
    &\lbr{v^+_1}{v^-_2}=-(2+k)((5+2k)\lm{2}-\tfrac{1}{2}(h+5J)\lm{1}-\tfrac{1}{4}(h+5J)'-\tfrac{1}{8(5+2k)}(h+5J)^2)+(4+k)\mathbb{L}_2,\\
    &\lbr{v^+_2}{v^-_1}=(2+k)((5+2k)\lm{2}-\tfrac{1}{2}(h-5J)\lm{1}-\tfrac{1}{4}(h-5J)'+\tfrac{1}{8(5+2 k)}(h-5J)^2)-(4+k)L_2,\\
    &\\
    &\lbr{v^\pm_1}{E}=\mp\tfrac{(10+3k)}{2(5+2k)}e_*v^\pm_1,\quad \lbr{v^\pm_2}{F}=\mp\tfrac{(10+3 k)}{2 (5+2 k)}f_*v^\pm_2,\\
    &\lbr{v^\pm_2}{E}=\mp\tfrac{(2+k)(10+3 k)}{5+2 k}v^\pm_1\lm{1}\pm(\tfrac{k}{2 (5+2 k)}e_*v^\pm_2+\tfrac{1}{2}h v^\pm_1-\tfrac{5+5 k+k^2}{(5+2 k)} v^\pm_1{}')-\tfrac{5(3+k)}{2 (5+2 k)}Jv^\pm_1,\\
    &\lbr{v^\pm_1}{F}=\mp\tfrac{(2+k)(10+3k)}{(5+2k)}v^\pm_2\lm{1}\pm(\tfrac{k}{2 (5+2 k)}f_*v^\pm_1-\tfrac{1}{2}h v^\pm_2-\tfrac{5+5k+k^2}{5+2 k} v^\pm_2{}')-\tfrac{5(3+k)}{2 (5+2 k)}Jv^\pm_2,\\
    &\\
    &\lbr{E}{E}=\tfrac{(2+k)(10+3 k)}{2(5+2k)}(e^2\lm{1}+ee'),\quad \lbr{F}{F}=\tfrac{(2+k)(10+3 k)}{2 (5+2 k)}(f^2\lm{1}+ff'),\\
    &\lbr{E}{F}=-(2+k)^2(10+3k)\lm{3}-\tfrac{(2+k)^2(10+3k)}{5+2k}h\lm{2}\\
    &\hspace{2cm}+((3+k)(5+k)(L_1+L_2)-\tfrac{(10+12 k+3 k^2)}{4}(h'+\tfrac{1}{(5+2 k)}h^2)-\tfrac{30+24 k+5k^2}{2 (5+2 k)}e_*f_*)\lm{1}\\
     &\hspace{2cm}+v^+_1 v^-_2 + v^-_1 v^+_2 +\tfrac{(5+k)}{2}h(L_1+L_2)-\tfrac{5(5+k)}{2(5+2k)}J(L_1-L_2)+\tfrac{(2+k)(5+k)}{2}(L_1+L_2)')\\
    &\hspace{2cm}-\tfrac{k+2}{2}(e_*f_*'+\tfrac{(5+3 k)}{(5+2k)}e'_*f_*)+\tfrac{5}{2(5+2k)}Je_*f_*-\tfrac{1}{2}he_*f_*-\tfrac{(1+k)(2+k)}{4}h'''\\
    &\hspace{2cm}\tfrac{1}{4(5+2 k)}(25(2+k)JJ'- (5+9k+3k^2)hh'-5h'J-\tfrac{k}{2}h^3-\tfrac{5}{(5+2k)}h^2J+\tfrac{25(2+k)}{2(5+2 k)}hJ^2),
\end{align*}
\end{fleqn}
where
\begin{align*}
 &H=-(5+k)L_1-(3+k) L_2+\tfrac{1}{2}h'+\tfrac{1}{2 (5+2 k)}h^2-e_*f_*,\\
 &\mathbb{L}_1=L_2-\tfrac{1}{4+k}(L_1+L_2)-\tfrac{1}{2(4+k)}(\tfrac{1}{(5+2 k)}h^2-2e_*f_*+h'),\\
 &\mathbb{L}_2=L_1+\tfrac{1}{4+k}(L_1+L_2)+\tfrac{1}{2(4+k)}(\tfrac{1}{(5+2 k)}h^2-2e_*f_*+h').
\end{align*}
We note that $\{\mathbb{L}_1,\mathbb{L}_2\}$ is another pair of mutually-commuting Virasoro fields which also commute $h,J$. The central charges are again $c_{5,3}(k)$ and $c_{3,1}(k+1)$, respectively.

\subsubsection{}{$\W^k(\sll_5,f_{1,1,1,2})$ at $k\notin \{-4,-5,-\frac{5}{2}\}$}\label{2,1,1,1}
\begin{fleqn}[\parindent]
\renewcommand{\arraystretch}{1.3}
\begin{align*}
{}\hspace{1cm}
\begin{array}{|c|c|c|c|} \hline
      \Delta &1  & 3/2 & 2 \\ \hline
      \text{gen.} & \begin{array}{c}
           e_{i,j} (1\leq i\neq j \leq 3)  \\
           h_1, h_2, J
      \end{array} & \begin{array}{c}
           v_1^+, v_2^+, v_3^+   \\
           v_1^-, v_2^-, v_3^-
      \end{array}  & L_1 \\ \hline 
     \end{array}
    \hspace{1cm} 
    \begin{array}{l}\lbr{L_1}{L_1}=P_{2,2}^{2,3}(k),\\ 
     \langle e_{i,j}, h_1, h_2 \rangle \simeq V^{k+1}(\sll_3),\\
     \lbr{J}{J}=\tfrac{3}{5}(5+2k)\lm{1}\end{array}
\end{align*}
\renewcommand{\arraystretch}{1}
\end{fleqn}
\begin{fleqn}[\parindent]
\begin{align*}
{}\hspace{1cm}
&\Span\{v_1^+,v_2^+,v^+_3\}=\C^3\subset \weyl_{\varpi_1}^{k+1}\otimes \pi^J_1,\qquad \Span\{v_1^-,v_2^-,v_3^-\}=\overline{\C}^3\subset \weyl_{\varpi_2}^{k+1}\otimes \pi^J_{-1},\\
&\lbr{L_1}{v^+_1}=\Delta v^+_1\lm{1}+v^+_1{}'-\tfrac{5}{3(5+2 k)}Jv^+_1-\tfrac{1}{(4+k)}(\tfrac{2}{3} h_1v^+_1+\tfrac{1}{3}h_2v^+_1+e_{1,2}v^+_2+e_{1,3}v^+_3),\\
&\lbr{L_1}{v^+_2}=\Delta v^+_2\lm{1}+v^+_2{}'-\tfrac{5}{3(5+2 k)}Jv^+_2-\tfrac{1}{(4+k)}(-\tfrac{1}{3}h_1v^+_2+\tfrac{1}{3}h_2v^+_2+e_{2,1}v^+_1+e_{2,3}v^+_3),\\
&\lbr{L_1}{v^+_3}=\Delta v^+_3\lm{1}+v^+_3{}'-\tfrac{5}{3(5+2 k)}Jv^+_3-\tfrac{1}{(4+k)}(-\tfrac{1}{3}h_1v^+_3-\tfrac{2}{3}h_2v^+_3+e_{3,1}v^+_1+e_{3,2}v^+_2),\\
&\lbr{L_1}{v^-_1}=\Delta v^-_1\lm{1}+v^-_1{}'+\tfrac{5}{3(5+2 k)}Jv^-_1+\tfrac{1}{(4+k)}(+\tfrac{2}{3}h_1v^-_1+\tfrac{1}{3}h_2v^-_1+e_{2,1}v^-_2+e_{3,1}v^-_3),\\
&\lbr{L_1}{v^-_2}=\Delta v^-_2\lm{1}+v^-_2{}'+\tfrac{5}{3 (5+2 k)}Jv^-_2+\tfrac{1}{(4+k)}(-\tfrac{1}{3}h_1v^-_2+\tfrac{1}{3}h_2v^-_2+e_{1,2}v^-_1+e_{3,2}v^-_3),\\
&\lbr{L_1}{v^-_3}=\Delta v^-_3\lm{1}+v^-_3{}'+\tfrac{5}{3(5+2k)}Jv^-_3+\tfrac{1}{(4+k)}(-\tfrac{1}{3}h_1v^-_3-\tfrac{2}{3}h_2v^-_3+e_{2,3}v^-_2+e_{1,3}v^-_1),\\
&\\
&\lbr{v^+_i}{v^+_j}=\lbr{v^-_i}{v^-_j}=0\quad (i,j=1,2,3),\\
&\lbr{v^+_1}{v^-_1}=-(1+k)(5+2k)\lm{2}-(\tfrac{2}{3}(5+2k)h_1+\tfrac{1}{3}(5+2k)h_2+\tfrac{5}{3}(1+k)J)\lm{1}\\
&\hspace{2cm}+(k+5)(L_1+L_{\mathrm{sug}})-(\tfrac{7}{9}h_1^2+\tfrac{7}{9}h_1h_2+\tfrac{10}{9}h_1J+\tfrac{4}{9}h_2^2+\tfrac{5}{9}h_2J+\tfrac{10}{9}J^2)\\
&\hspace{2cm}-(2e_{1,2}e_{2,1}+2 e_{1,3}e_{3,1}+e_{2,3}e_{3,2})+\tfrac{(1-2 k)}{3} h_1'+\tfrac{2-k}{3}h_2'-\tfrac{5(1+k)}{6} J',\\
&\lbr{v^+_2}{v^-_2}=-(1+k)(5+2k)\lm{2}-(\tfrac{-1}{3}(5+2k)h_1+\tfrac{1}{3}(5+2k)h_2+\tfrac{5}{3}(1+k)J)\lm{1}\\
&\hspace{2cm}+(5+k)(L_1+L_{\mathrm{sug}})-(\tfrac{4}{9}h_1^2+\tfrac{1}{9}h_1h_2-\tfrac{5}{9}h_1J+\tfrac{4}{9}h_2^2+\tfrac{5}{9}h_2J+\tfrac{10}{9}J^2)\\
&\hspace{2cm}-(2e_{1,2}e_{2,1}+e_{1,3}e_{3,1}+2e_{2,3}e_{3,2})+\tfrac{7+k}{3}h_1'+\tfrac{2-k}{3}h_2'-\tfrac{5(1+k)}{6} J',\\
&\lbr{v^+_3}{v^-_3}=-(1+k) (5+2 k)\lm{2}-(\tfrac{-1}{3}(5+2 k) h_1-\tfrac{2}{3}(5+2k)h_2+\tfrac{5}{3}(1+k)J)\lm{1}\\
&\hspace{2cm}+(5+k)(L_1+L_{\mathrm{sug}})-(\tfrac{4}{9}h_1^2+\tfrac{7}{9} h_1h_2-\tfrac{5}{9}h_1J+\tfrac{7}{9}h_2^2-\tfrac{10}{9}h_2J+\tfrac{10}{9}J^2)\\
&\hspace{2cm}-(e_{1,2}e_{2,1}+2e_{1,3}e_{3,1}+2e_{2,3}e_{3,2})+\tfrac{(7+k)}{3} h_1'+\tfrac{(11+2k)}{3}h_2'-\tfrac{5(1+k)}{6}J'\\
&\lbr{v^+_1}{v^-_2}=-(2k+5)e_{1,2}\lm{1}-\tfrac{1}{3}h_1e_{1,2}-\tfrac{2}{3}h_2e_{1,2}-\tfrac{5}{3}Je_{1,2}-e_{1,3}e_{3,2}-(k+2)e_{1,2}',\\
&\lbr{v^+_1}{v^-_3}=-(2k+5)e_{1,3}\lm{1}-\tfrac{1}{3}h_1e_{1,3}+\tfrac{1}{3}h_2e_{1,3}-\tfrac{5}{3}Je_{1,3}-e_{1,2}e_{2,3}-(k+2)e_{1,3}',\\
&\lbr{v^+_2}{v^-_1}=-(2k+5)e_{2,1}\lm{1}-\tfrac{1}{3}h_1e_{2,1}-\tfrac{2}{3}h_2e_{2,1}-\tfrac{5}{3}Je_{2,1}-e_{2,3}e_{3,1}-(k+2)e_{2,1}',\\
&\lbr{v^+_2}{v^-_3}=-(2k+5)e_{2,3}\lm{1}+\tfrac{2}{3}h_1e_{2,3}+\tfrac{1}{3}h_2e_{2,3}-\tfrac{5}{3}Je_{2,3}-e_{1,3}e_{2,1}-(k+3) e_{2,3}',\\
&\lbr{v^+_3}{v^-_1}=-(2k+5)e_{3,1}\lm{1}-\tfrac{1}{3}h_1e_{3,1}+\tfrac{1}{3}h_2e_{3,1}-\tfrac{5}{3}Je_{3,1}- e_{2,1}e_{3,2}-(k+3) e_{3,1}',\\
&\lbr{v^+_3}{v^-_2}=-(2k+5)e_{3,2}\lm{1}+\tfrac{2}{3}h_1e_{3,2}+\tfrac{1}{3}h_2e_{3,2}-\tfrac{5}{3}Je_{3,2}-e_{1,2}e_{3,1}-(k+3)e_{3,2}',
\end{align*}
\end{fleqn}
where $\Delta=\tfrac{(2+k)(10+3 k)}{(4+k)(5+2k)}$ and $L_{\mathrm{sug}}$ is the Sugawara vector for $V^{k+1}(\sll_3)\otimes \pi^J$.

\section{OPEs of \tWinfcl-algebra} \label{Winfinity algebra}
We set 
\begin{align*}
    [W_{n\Lambda}W_m]=P_{n,m}\in \W_\infty[c,\lambda][\Lambda],\quad (2 \leq n\leq m)
\end{align*}
where $W_2=L$. 
For lower $n,m$, $P_{n,m}$ are determined as follows \cite{L}.
\begin{align*}
   &P_{2,2}=\frac{c}{2} \lm{3}+ 2 L \lm{1}+L',\quad P_{2,3}=3 W_3 \lm{1}+W_3'\\
   &\\
   &P_{3,3}=\frac{c}{3}\lm{5}+2L\lm{3}+L'\lm{2}+W_4\lm{1}+\frac{1}{2}W_4'- \frac{1}{12}L'''\\
   &P_{2,4}=3c\lm{5}+10L\lm{3}+3L'\lm{2}+4W_4\lm{1}+W_4'\\
   &\\
   &P_{2,5}=-5(-37+16(2+c)\lambda)W_3\lm{3}-(-55+16(2+c)\lambda)W_3'\lm{2}+5W_5\lm{1}+W_5'\\
   &P_{3,4}=-(-31+16(2+c)\lambda)W_3\lm{3}-\frac{8}{3}(-5+2(2+c)\lambda)W_3'\lm{2}+W_5\lm{1}\\
   &\hspace{2cm}+\frac{2}{15}(48\lambda LW_3'-72 \lambda L'W_3+3W_5'+(-5+2(-1+c) \lambda)W_3''')\\
   &\\
   &P_{2,6}=-13c(-55+16\lambda(2+c))\lm{7}+(2100-768\lambda(2+c))L\lm{5}+(770-224\lambda(2+c))L'\lm{4}\\
   &\hspace{2cm}+((660-80\lambda(13+5c))W_4+640\lambda L^2+(50+40 \lambda(-1+c))L'')\lm{3}\\
   &\hspace{2cm}+((195-12\lambda(17+7c))W_4'+192\lambda L'L+\frac{1}{6}(-65+4\lambda(31+17c))L''')\lm{2}\\
   &\hspace{2cm}+6W_6\lm{1}+W_6'\\
   &P_{3,5}=-c(-55+16\lambda(2+c))\lm{7}-\frac{4}{3}(-175+64\lambda(2+c))L\lm{5}+(110-32\lambda(2+c))L'\lm{4}\\
   &\hspace{2cm}+((95-16\lambda(11+4c))W_4+128\lambda L^2+(10+8\lambda(-1+c))L'')\lm{3}\\
   &\hspace{2cm}+(64\lambda L'L+(\frac{75}{2}-4 \lambda(13+5c))W_4'+\frac{1}{12}(-25+8\lambda(9+5c))L''')\lm{2}+W_6\lm{1}\\
   &\hspace{2cm}+\frac{1}{3}W_6'+\frac{32 \lambda}{3}LW_4'-\frac{64\lambda}{3}L'W_4-\frac{16 \lambda}{3}L'''L +(-\frac{5}{4}+\frac{2}{3}\lambda(1+c))W_4'''\\
   &\hspace{2cm}+(\frac{5}{72}-\frac{1}{45} \lambda (13 + 5 c)\big) \partial^5 L,\\
   &P_{4,4}=-\frac{1}{3}c(-139+16\lambda(2+c))\lm{7}-\frac{4}{3}(-125+32\lambda(2+c))L\lm{5}+(\frac{250}{3}-\frac{64}{3}\lambda(2+c))L'\lm{4}\\
   &\hspace{2cm}+((72-48\lambda(3+c))W_4+128\lambda L^2+(10+8\lambda(-1+c))L'')\lm{3}\\
   &\hspace{2cm}+(128\lambda L'L+(36-24\lambda (3+c))W_4'+\frac{1}{18}(-35+8\lambda(23+13c))L''')\lm{2}\\ 
   &\hspace{2cm}+(\frac{4}{5}W_6+\frac{64 \lambda}{5} LW_4-\frac{288 \lambda}{5}W_3^2+32\lambda L''L+16\lambda (L')^2\\ 
   &\hspace{3cm}+\frac{1}{15}(35-4\lambda(19+11c))W_4''+\frac{1}{90}(-5+4\lambda(7+23c))L'''')\lm{1}\\ 
   &\hspace{2cm}+(-\frac{2}{5}W_6'-\frac{32 \lambda}{5}LW_4'+\frac{288 \lambda}{5}W_3'W_3-\frac{32 \lambda}{5}L'W_4\\ 
   &\hspace{3cm}-\frac{16 \lambda}{3}L'''L +(\frac{11}{6}-\frac{16 \lambda}{15}-\frac{8 \lambda c}{15})W_4'''+( -\frac{1}{4}+\frac{8 \lambda}{25}) L''''')\\
   &\\
   &P_{2,7}=18 (4725-4784\lambda(2+c)+256\lambda^2(26+23c+5c^2))W_3\lm{5}\\ 
   &\hspace{2cm}+14(2225-1920\lambda(2+c)+64\lambda^2(34+31c+7c^2))W_3'\lm{4}\\
   &\hspace{2cm}+(-5(-357+8\lambda(97+31c))W_5-640\lambda (-35+8\lambda(2+c))LW_3\\  
   &\hspace{3cm}+\frac{5}{2}(805-8\lambda(19+27c)+128\lambda^2(6+5c+c^2))W_3'')\lm{3}\\
   &\hspace{2cm}+(-\frac{3}{5}(-875+32\lambda(39+14c))W_5'-\frac{64}{5} \lambda(-425+4\lambda(79+29c))LW_3\\ 
   &\hspace{2cm}+\frac{288}{5}\lambda(5+4\lambda(13+3c))L'W_3+(-\frac{875}{2}+152\lambda(5+3c)-\frac{32}{5}\lambda^2(-23+15c+8c^2))W_3''')\lm{2}\\ 
   &\hspace{2cm}+7 W_7\lm{1} +W_7'\\
   &P_{3,6}= 2 \big(4375 - 4656 \lambda (2 + c) + 256 \lambda^2 (26 + 23 c + 5 c^2)\big)W_3\lm{5} \\ 
   &\hspace{2cm}+4\big(975-920\lambda(2+c)+32\lambda^2(34+31c+7c^2)\big)W_3'\lm{4}\\
   &\hspace{2cm}+\bigg(\big(225-8\lambda(71+21c)\big)W_5-128\lambda \big(-29+8\lambda(2+c)\big) LW_3 \\ 
   &\hspace{3cm}+\big(\frac{665}{2}-4\lambda(53+41c)+64\lambda^2(6+5c+c^2)\big) W_3''\bigg)\lm{3}\\ 
   &\hspace{2cm}+\bigg(\big(84-\frac{4}{5}\lambda(193+63c)\big)W_5'-\frac{32}{15}\lambda\big(-505+4\lambda(107+37c)\big)LW_3' \\ 
   &\hspace{2cm}-\frac{48}{5}\lambda\big(-55+4\lambda(-9+c)\big)L'W_3 \\ 
   &\hspace{2cm}+\big(-70+\lambda(\frac{490}{3}+82c)-\frac{16}{15}\lambda^2(-29+20c+9c^2)\big)W_3'''\bigg)\lm{2}\\ 
   &\hspace{2cm}+W_7\lm{1}+\bigg(\frac{2}{7}W_7'+\frac{496 \lambda}{35} LW_5'-\frac{248 \lambda}{7}L'W_5 \\ 
   &\hspace{2cm}+\frac{192 \lambda}{7}W_3W_4'-\frac{256 \lambda}{7} W_3'W_4+\frac{1536 \lambda^2}{35}L'LW_3-\frac{1024 \lambda^2}{35}L^2W_3' \\ 
   &\hspace{2cm}+\frac{8}{35} \lambda \big(-455+4\lambda(135+41c)\big)L'''W_3-\frac{192}{35} \lambda \big(5+2\lambda (-3+c)\big)L''W_3'\\ 
   &\hspace{2cm}+\frac{12}{35} \lambda \big(95+8\lambda (-3+c)\big) L'W_3''+\frac{8}{105} \lambda \big(-455+8\lambda (-25+7c)\big)LW_3''' \\ 
   &\hspace{2cm}+\big(-2 + \frac{2}{35} \lambda (17 + 21 c)\big)W_5''' +\frac{1}{105}\big(175-\lambda (149+205c)+24\lambda^2(11+c^2)\big)W_3'''''\bigg)\\
 &P_{4,5}=\big(4950-4928 \lambda (2+c)+256 \lambda^2 (26+23c+5c^2) \big)W_3\lm{5} 
\\ &\hspace{2cm}+\frac{2}{3} \big(3625-3600\lambda (2+c)+128\lambda^2(34+31c+7c^2)\big)W_3\lm{4}
\\ &\hspace{2cm}+\bigg(\big(140-8\lambda (49+13c)\big) W_5-128 \lambda \big(-23+8\lambda(2+c)\big)LW_3  
\\ &\hspace{2cm}+\big(\frac{525}{2}-4\lambda (87+55c)+64\lambda^2 (6+5c+c^2)\big)W_3''\bigg)\lm{3}
\\ &\hspace{2cm}+\bigg(\big(64-\frac{16}{5}\lambda(51+14c)\big)W_5'-\frac{32}{15} \lambda \big(-485+16\lambda(34+11c)\big)LW_3' 
\\ &\hspace{2cm}-\frac{48}{5}\lambda \big(-145+16\lambda (2+3c)\big)L'W_3+\frac{1}{30} \big(-1575+40\lambda(127+43c)
\\ &\hspace{2cm}-256\lambda^2 (-4+3c+c^2)\big)W_3'''\bigg)\lm{2} + \bigg(\frac{2}{3}W_7+\frac{64 \lambda}{3}LW_5-128 \lambda W_3W_4
\\ &\hspace{2cm}-\frac{32}{5}\lambda \big(-95+2\lambda(65+19c)\big)L''W_3-\frac{32}{15} \lambda \big(-125+2\lambda (65+19c)\big)L'W_3' 
\\ &\hspace{2cm}-\frac{32}{15}\lambda \big(35+4\lambda(-25+c)\big)LW_3''+\big( \frac{5}{2}-\frac{8}{5} \lambda (5+2c)\big)W_5'' 
\\ &\hspace{2cm}+\frac{1}{180} \big(-175+80\lambda (33+c)-64\lambda^2(65-6c+c^2)\big)W_3''''\bigg)\lm{1}
\\ &\hspace{2cm}+\bigg(\frac{2}{7}W_7'+\frac{384 \lambda}{35}LW_5'+\frac{32 \lambda}{7}L'W_5+\frac{192 \lambda}{7}W_3W_4' 
\\ &\hspace{2cm}-\frac{1152 \lambda}{7}W_3'-\frac{9216 \lambda^2}{35}L'LW_3+\frac{6144 \lambda^2}{35}L^2W_3' 
\\ &\hspace{2cm}-\frac{8}{105} \lambda \big(-2345+8\lambda (389+145c)\big)L'''W_3 
\\ &\hspace{2cm}-\frac{32}{35} \lambda \big(-145+8\lambda (13+5c)\big)L''W_3' 
+ \frac{8}{35} \lambda \big(-295+8\lambda(13+5c)\big)L'W_3'' 
\\ &\hspace{2cm}+\frac{16}{35} \lambda (-245+8\lambda(11+7c)\big)LW_3''' 
+ \big(-\frac{17}{6}+\frac{4}{105}\lambda (43+35c) \big)W_5'''
\\ &\hspace{2cm}+\frac{1}{420} \big(1925-16 \lambda (-87+130c)+64\lambda^2 (-29+5c^2)\big)W_3'''''\bigg)
\end{align*}

\end{document}